\documentclass[12pt,a4paper]{article}
\usepackage{color}
\usepackage{amscd,amsmath,amssymb,amsfonts,times}
\usepackage{graphicx}
\usepackage{mathrsfs}
\usepackage[all]{xy}
    \newtheorem{thm}{Theorem}[section]
    
    \newtheorem{lem}[thm]{Lemma}
    \newtheorem{prop}[thm]{Proposition}
    \newtheorem{cor}[thm]{Corollary}

    \newtheorem{defn}[thm]{Definition}
    \newtheorem{exmp}[thm]{Example}
    \newtheorem{rem}[thm]{Remark}
\numberwithin{equation}{thm}

\newcommand{\qed}
{\mbox{}\nolinebreak$\square$\medbreak\par}
\newenvironment{pf}{\par\smallskip\noindent{\it Proof.}}{\hfill\qed\par\smallskip}
\newenvironment{pf*}[1]{\par\smallskip\noindent\emph{#1.}}{\hfill\qed\par\smallskip}
\date{}
\title{Chern class and Riemann-Roch theorem for cohomology theory without homotopy invariance}
\author{Masanori Asakura\footnote{Supported by Grant-in-Aid for Scientific research C-24540001}\; and Kanetomo Sato\footnote{Supported by Grant-in-Aid for Scientific research B-23340003}}
\begin{document}
\maketitle

\def\back{\sharp}
\def\push{!}

\def\A{{\sf A}}
\def\Ab{{\sf Ab}}
\def\ab{{\sf ab}}
\def\ACov{{\sf ACov}}
\def\alg{{\sf alg}}
\def\an{{\sf an}}
\def\Aut{{\sf Aut}}
\def\BGL{{\sf B}_{\star}{\sf GL}}
\def\BGLs{{\text{\it B}}_{\star}{\text{\it GL}}}
\def\BG{{\text{\it B}}_{\star}\cG}
\def\BQP{{\sf B}_\star {\text{\it QP}}}
\def\BQPs{{\text{\it B}}_\star{\mathscr {Q\hspace{-1pt}P}}}
\def\Br{{\sf Br}}
\def\can{{\sf can}}
\def\CaCl{{\sf Cl}}
\def\cCov{{\mathscr C}\hspace{-.7pt}o\hspace{-.1pt}v}
\def\cd{{\sf cd}}
\def\ceC{\check{C}}
\def\cecH{\check{\mathscr H}}
\def\ceH{\check{H}}
\def\ch{{\sf ch}}
\def\CH{{\sf CH}}
\def\Char{{\sf Char}}
\def\Cl{{\sf Cl}}
\def\cl{{\sf cl}}
\def\codim{{\sf codim}}
\def\Coker{{\sf Coker}}
\def\Coim{{\sf Coim}}
\def\cone{\text{\rm Cone}}
\def\cont{{\sf cont}}
\def\Cor{{\sf Cor}}
\def\Cov{{\sf Cov}}
\def\crys{{\sf crys}}
\def\Cyc{{\sf Cyc}}
\def\deg{{\sf deg}}
\def\det{{\sf det}}
\def\disc{{\sf disc}}
\def\dim{{\sf dim}}
\def\div{{\sf div}}
\def\Div{{\sf Div}}
\def\dR{{\sf dR}}
\def\DR{{\sf DR}}
\def\log{{\sf log}}
\def\dlog{d\hspace{.5pt}\log}
\def\End{{\sf End}}
\def\et{{\sf \acute{e}t}}
\def\Et{{\sf \acute{E}t}}
\def\exp{{\sf exp}}
\def\Ext{{\sf Ext}}
\def\FAb{{\sf FAb}}
\def\FEt{{\sf F\acute{E}t}}
\def\fetc{{\sf f\acute{e}tc}}
\def\Fib{{\sf Fib}}
\def\Fibs{{\text{\it Fib}}}
\def\Fil{{\sf Fil}}
\def\fppf{{\sf fppf}}
\def\Frac{{\sf Frac}}
\def\Frob{{\sf Frob}}
\def\FSep{{\sf FSep}}
\def\FSets{{\sf FSets}}
\def\fslc{{\sf lcfs}}
\def\gal{{\sf gal}}
\def\Gal{{\sf Gal}}
\def\geo{{\sf geo}}
\def\GL{{\sf GL}}
\def\GLs{{\text{\it GL}}}
\def\Gmod{G\text{-}{\sf mod}}
\def\gr{{\sf gr}}
\def\Gr{{\sf Gr}\text{-}}
\def\GrAb{{\sf Gr\text{-}Ab}}
\def\GrRing{{\sf Gr\text{-}Ring}}
\def\Gps{{\sf GrSch}}
\def\gys{{\sf Gys}}
\def\h{{\sf \hspace{1pt}h}}
\def\H{H}
\def\Hom{{\sf Hom}}
\def\sHom{{\mathscr H}\hspace{-2pt}om}
\def\sMor{{\mathscr M}\hspace{-2pt}or}
\def\hom{{\sf hom}}
\def\homap{[-]}
\def\homapp{[-]}
\def\id{\text{\rm id}}
\def\Image{{\sf Im}}
\def\Ind{{\sf Ind}}
\def\ker{{\sf Ker}}
\def\Ker{{\sf Ker}}
\def\LCI{LCI}
\def\lim{{\sf lim}}
\def\lin{{\sf lin}}
\def\logsyn{{\sf log}\text{-}{\sf syn}}
\def\M{\mathbf{M}}
\def\Map{{\sf Map}}
\def\max{{\sf max}}
\def\MHS{{\sf MHS}}
\def\MHM{{\sf MHM}}
\def\mod{\;{\sf mod}\;}
\def\Mor{{\sf Mor}}
\def\mor{{\sf mor}}
\def\Nat{{\sf Nat}}
\def\nbd{{\sf nbd}}
\def\ncd{{\sf ncd}}
\def\nil{{\sf nil}}
\def\Nis{{\sf{Nis}}}
\def\nr{{\sf nr}}
\def\NSwT{{\sf NSwT}}
\def\Ob{{\sf Ob}}
\def\OBQP{{\sf\Omega}\BQP}
\def\OBQPs{{\sf\Omega}\BQPs}
\def\op{{\sf op}}
\def\ord{{\sf ord}}
\def\Ouv{{\sf Ouv}}
\def\perf{{\sf perf}}
\def\Pic{{\sf Pic}}
\def\pr{{\sf pr}}
\def\pre{{\sf pre}}
\def\Proj{{\sf Proj}}
\def\qis{\text{\rm qis}}
\def\Qis{{\sf Qis}}
\def\qpr{{\sf qp}}
\def\QP{{\text{\it QP}}}
\def\QPs{{\mathscr {Q\hspace{-1pt}P}}}
\def\rank{{\sf rank}}
\def\rat{{\sf rat}}
\def\rdim{{\sf rel}\text{-}{\sf dim}}
\def\red{{\sf red}}
\def\Re{{\sf Re}}
\def\Reg{{\sf Reg}}
\def\reg{{\sf reg}}
\def\res{{\sf Res}}
\def\Res{{\sf Res}}
\def\Rings{{\sf Rings}}
\def\rk{{\sf rk}}
\def\Rmod{R\text{-}{\sf mod}}
\def\rSym{{\sf Sym}}
\def\sAb{{\sf sAb}}
\def\Sch{{\sf Sch}}
\def\scr{{\diamondsuit}}
\def\sdeg{{\sf sd}}
\def\Sets{{\sf Sets}}
\def\sgn{{\sf sgn}}
\def\sh{{\sf \hspace{1pt}sh}}
\def\SLCI{‹­LCI}
\def\sm{{\sf sm}}
\def\Sm{{\sf Sm}}
\def\sp{{\sf sp}}
\def\Spec{{\sf Spec}}
\def\ss{{\sf ss}}
\def\ssm{\smallsetminus}
\def\Sym{{\text{\it Sym}}}
\def\syn{{\sf syn}}
\def\taut{{\sf taut}}
\def\Td{{\sf td}}
\def\top{{\sfE}}
\def\Tor{{\sf Tor}}
\def\tor{{\sf tors}}
\def\tot{{\sf tot}}
\def\tr{{\sf tr}}
\def\Tr{{\sf Tr}}
\def\trdeg{{\sf tr.deg}}
\def\univ{{\sf univ}}
\def\uu{\boldsymbol{u}}
\def\vv{\boldsymbol{v}}
\def\val{{\sf val}}
\def\vdim{{\sf vdim}}
\def\Zar{{\sf Zar}}
\def\zet{{\sf Z\et}}

\def\neg{\hspace{-1pt}}

\def\II{I\hspace{-0.7pt}I}
\def\III{I\hspace{-0.7pt}I\hspace{-0.7pt}I}
\def\IV{I\hspace{-0.7pt}V}
\def\VII{V\hspace{-0.7pt}I\hspace{-0.7pt}I}

\def\bA{\mathbb A}
\def\bC{\mathbb C}
\def\bE{\mathbb E}
\def\bF{\mathbb F}
\def\bG{\mathbb G}
\def\bH{\mathbb H}
\def\bI{\mathbb I}
\def\bK{\mathbb K}
\def\bL{\mathbb L}
\def\bN{\mathbb N}
\def\bP{\mathbb P}
\def\bQ{\mathbb Q}
\def\bR{\mathbb R}
\def\bT{\mathbb T}
\def\bZ{\mathbb Z}

\def\cA{\mathscr A}
\def\cB{\mathscr B}
\def\cC{\mathscr C}
\def\cD{\mathscr D}
\def\cE{\mathscr E}
\def\cF{\mathscr F}
\def\cG{\mathscr G}
\def\cH{\mathscr H}
\def\cI{\mathscr I}
\def\cJ{\mathscr J}
\def\cK{\mathscr K}
\def\cL{\mathscr L}
\def\cM{\mathscr M}
\def\cN{\mathscr N}
\def\cO{\mathscr O}
\def\cP{\mathscr P}
\def\cQ{\mathscr Q}
\def\cR{\mathscr R}
\def\cS{\mathscr S}
\def\cT{\mathscr T}
\def\cU{\mathscr U}
\def\cV{\mathscr V}
\def\cX{\mathscr X}
\def\cY{\mathscr Y}
\def\cZ{\mathscr Z}

\def\fT{{\mathfrak T}}

\def\fra{\mathfrak a}
\def\frb{\mathfrak b}
\def\fc{\mathfrak c}
\def\fd{\mathfrak d}
\def\fe{\mathfrak e}
\def\frg{\mathfrak g}
\def\fk{\mathfrak k}
\def\fl{\mathfrak l}
\def\fm{\mathfrak m}
\def\fn{\mathfrak n}
\def\fp{\mathfrak p}
\def\fq{\mathfrak q}
\def\fv{{\sf v}}
\def\fw{{\sf w}}

\def\rC{{\text {\it C}}}

\def\sfA{{\sf A}}
\def\sfB{{\sf B}}
\def\sfC{{\sf C}}
\def\sfD{{\text{\it D}}}
\def\sfE{{\sf E}}
\def\sfF{{\sf F}}
\def\sfG{{\sf G}}
\def\sfH{{\sf H}}
\def\sfI{{\sf I}}
\def\sfJ{{\sf J}}
\def\sfK{{\sf K}}
\def\sfL{{\sf L}}
\def\sfM{{\sf M}}
\def\sfN{{\sf N}}
\def\sfO{{\sf O}}
\def\sfP{\text{\it P}}
\def\sfQ{{\sf Q}}
\def\sfR{{\sf R}}
\def\sfS{{\text{\it Shv}}}
\def\sfT{{\sf T}}
\def\sfU{{\sf U}}
\def\sfX{{\sf X}}
\def\sfZ{{\sf Z}}
\def\sOmega{{\sf \Omega}}
\def\sfa{{\sf a}}
\def\sfb{{\sf b}}
\def\sfc{{\sf c}}
\def\sfd{{\sf d}}
\def\sfe{{\sf e}}
\def\sff{{\sf f}}
\def\sffa{{\sf g}}
\def\sfk{\text{\it k}}
\def\sfh{\text{\it h}}
\def\sfm{{\sf m}}
\def\sfp{\text{\it P}}
\def\sfs{{\sf s}}
\def\sft{{\sf t}}
\def\sfu{{\sf u}}
\def\sfz{{\sf z}}

\def\k{\kappa}
\def\ep{\epsilon}
\def\G{\Gamma}
\def\lam{\lambda}
\def\D{\Delta}
\def\vD{\varDelta}
\def\rD{\Delta}
\def\vG{\varGamma}
\def\vL{\varLambda}
\def\vS{\varSigma}
\def\vT{\Theta}
\def\ve{\varepsilon}
\def\vare{\ve}
\def\ra{\rightarrow}
\def\lra{\longrightarrow}
\def\lla{\longleftarrow}
\def\Lra{\Longrightarrow}
\def\Llra{\Longleftrightarrow}
\def\Ra{\Rightarrow}
\def\La{\Leftarrow}
\def\hra{\hookrightarrow}
\def\lmt{\longmapsto}
\def\thra{\twoheadrightarrow}

\def\rmapo#1{\overset{#1}{\lra}}
\def\rmapu#1{\underset{#1}{\lra}}
\def\lmapu#1{\underset{#1}{\lla}}
\def\rmapou#1#2{\overset{#1}{\underset{#2}{\lra}}}
\def\lmapo#1{\overset{#1}{\lla}}

\def\bs#1{\boldsymbol{#1}}
\def\wt#1{\widetilde{#1}}
\def\wh#1{\widehat{#1}}
\def\spt{\sptilde}
\def\ol#1{\overline{#1}}
\def\ul#1{\underline{#1}}
\def\us#1#2{\underset{#1}{#2}}
\def\os#1#2{\overset{#1}{#2}}
\def\lim#1{\us{#1}{\varinjlim}}

\def\Gm{\bG_{\hspace{-1pt}{\sf m}}}
\def\Ga{\bG_{\hspace{-1pt}{\sf a}}}
\def\BGm{{\sf B}_{\star}\Gm}

\def\KM#1{\sfK^\sfM_{#1}}
\def\zlm{\bZ/\ell^m}
\def\zlmn{\bZ/\ell^{m+n}}
\def\zln{\bZ/\ell^n}
\def\zl{\bZ_\ell}
\def\zp{\bZ_p}
\def\Zlam{\bZ_{(\ell)}}
\def\Zlamm{\bZ_{(\lam)}}
\def\Zn{\bZ/n}
\def\Zm{\bZ/m}
\def\Zp{\bZ/p}
\def\Zpm{\bZ/p^m}
\def\Zpr{\bZ/p^r}
\def\Zpn{\bZ/p^n}
\def\ql{\bQ_\ell}
\def\qp{\bQ_p}
\def\qzl{\bQ_\ell/\bZ_\ell}
\def\qzp{\bQ_p/\bZ_p}
\def\qz{\bQ/\bZ}
\def\Wm{W_{\hspace{-2pt}m}}
\def\Wn{W_{\hspace{-2pt}n}}
\def\Wnn{W_{\hspace{-2pt}n+1}}
\def\witt#1#2#3{W_{\hspace{-2pt}#2}{\hspace{1pt}}\Omega_{#1}^{#3}}
\def\logwitt#1#2#3{W_{\hspace{-2pt}#2}{\hspace{1pt}}\Omega_{#1,\log}^{#3}}
\def\logwittt#1#2{W_{\hspace{-2pt}#1}{\hspace{1pt}}\Omega_\log^{#2}}

\def\isom{\hspace{9pt}{}^\sim\hspace{-16.5pt}\lra}
\def\lisom{\hspace{9pt}{}^\sim\hspace{-16.5pt}\lla}
\def\bw{\boldsymbol{\wedge}}
\def\psr#1{\hspace{.5pt}[\hspace{-1.57pt}[{#1}]\hspace{-1.58pt}]}

\def\fD{f^*\hspace{-1.5pt}D}
\def\fE{f^*\hspace{-1.5pt}E}
\def\fL{f^*\hspace{-1.5pt}L}
\def\fsfc{f^\back\sfc}
\def\gsfc{g^\back\hspace{-.2pt}\sfc}
\def\psfc{p^\back\hspace{-.2pt}\sfc}
\def\gE{g^*\hspace{-1.5pt}E}
\def\pE{p^*\hspace{-1.5pt}E}

\def\tA{\wt{A}}
\def\wtB{\wt{B}}
\def\whA{\wh{A}}
\def\tB{\text{\it B}}
\def\tC{\wt{C}}
\def\tD{\wt{D}}
\def\tH{\text{\it H}\hspace{1.2pt}}
\def\wtH{\wt{\text{\it H}}{\hspace{1.2pt}}}
\def\whH{\wh{\text{\it H}}{\hspace{1.2pt}}}
\def\tK{\text{\it K}}
\def\tM{\wt{M}}
\def\tN{\wt{N}}
\def\tR{\text{\it R}}
\def\tV{\wt{V}}
\def\tW{\wt{W}}
\def\tX{\wt{X}}
\def\tY{\wt{Y}}
\def\whY{\wh{Y}}
\def\tZ{\wt{Z}}

\def\bSpec{\bs{\Spec}}
\def\bProj{\bs{\Proj}}
\def\cHs{{\cH\!{\text{\it o}}_\bullet}}
\def\cHsn{{\cH\!{\text{\it o}}}}
\def\cHsu{{\cH\!{\text{\it o}}}}
\def\cHss{{\cH\!{\text{\it o}}}}
\def\sC{{\text{\it C}}}
\def\sfGa{{\sf \Gamma}}
\def\sfGab{{\sf \Gamma}_{\!\tB}}
\def\sfGad{{\sf \Gamma}_{\!\cD}}
\def\SS{\vD^\op{\text{\it Shv}}_\bullet}
\def\SSu{\vD^\op{\text{\it Shv}}}

\def\As{A_{\hspace{-0.2pt}\star}}
\def\Bs{B_{\hspace{-0.5pt}\star}}
\def\Ds{D_{\hspace{-0.7pt}\star}}
\def\Es{E_{\hspace{-0.5pt}\star}}
\def\Fs{F_{\hspace{-0.7pt}\star}}
\def\Gs{G_{\hspace{-0.7pt}\star}}
\def\Is{I_{\hspace{-0.5pt}\star}}
\def\Ls{L_{\hspace{-0.5pt}\star}}
\def\Ms{M_{\hspace{-0.5pt}\star}}
\def\Ns{N_{\hspace{-0.5pt}\star}}
\def\Vs{V_{\hspace{-0.7pt}\star}}
\def\Xs{X_{\hspace{-0.5pt}\star}}
\def\Ys{Y_{\hspace{-1pt}\star}}
\def\cEs{\cE_\star}
\def\cFs{\cF_{\hspace{-0.7pt}\star}}
\def\cJs{\cJ_{\hspace{-0.7pt}\star}}
\def\fEs{f^*\hspace{-1.5pt}\Es}
\def\fLs{f^*\hspace{-1.5pt}\Ls}
\def\pEs{p^*\hspace{-1.5pt}\Es}
\def\cFp{\cF_{\hspace{-1pt}p}}
\def\Us{U_{\hspace{-0.7pt}\star}}
\def\Zs{Z_{\hspace{-0.7pt}\star}}

\def\scbullet{\text{\small$\bullet$}}
\def\tstar{\text{\tiny$\bigstar$}}
\def\scstar{\text{\scriptsize$\bigstar$}}
\def\bb{\hspace{1pt}--\hspace{0pt}}

\def\rz{z}

\def\nega{\vspace{-8pt}}
\vspace{-1.3cm}
\begin{center}
{\large with Appendix B by Kei Hagihara}
\end{center}
\begin{abstract}
In this paper, we formulate axioms of certain graded cohomology theory and define higher Chern class maps following the method of Gillet \cite{Gi}. We will not include homotopy invariance nor purity in our axioms. It will turn out that the Riemann-Roch theorem without denominators holds for our higher Chern classes. We will give two applications of our Riemann-Roch results in \S\S\ref{sect11}\bb\ref{sect13}.
\end{abstract}
\par
\section{Introduction}
In his papers \cite{Gr} and \cite{Gr2}, Grothendieck defined Chern classes and characters
\[ \sfc_i : \tK_0(X) \lra \CH^i(X), \qquad \ch_X : \tK_0(X) \lra \CH^*(X)_\bQ \]
for a smooth variety $X$ over a field $k$, where $\tK_0(X)$ (resp.\ $\CH^i(X)$) denotes the Grothendieck group of vector bundles over $X$(resp.\ the Chow groups of algebraic cycles of codimension $i$ on $X$ modulo rational equivalence). Concerning the Chern character, he proposed the celebrated Grothendieck-Riemann-Roch theorem, which asserts that for a proper morphism $f : Y \to X$ of smooth varieties over $k$, the equality
\stepcounter{thm}
\begin{equation}\label{eq1-1a}
  \ch_X(f_*\alpha) \cdot \Td(T_X) = f_\push(\ch_Y(\alpha) \cdot \Td(T_Y))
\end{equation}
holds in $\CH^*(X)_\bQ$ for any $\alpha \in \tK_0(X)$. Here $\Td(T_X)$ denotes the Todd class of the tangent bundle $T_X$ of $X$, and $f_*$ (resp.\ $f_\push$) denotes the push-forward of Grothendieck groups (resp.\ Chow rings). One immediately recovers the classical Riemann-Roch theorem for a smooth complete curve $X$ of genus $g$ with canonical divisor $K$:
\[ \ell(D)-\ell(K-D)=\deg(D)-g+1 \quad \hbox{ for a divisor $D$ on $X$,} \]
by considering the case of the structure morphism $X \to \Spec(k)$.

In \cite{Gi}, Gillet introduced certain axioms on graded cohomology theory $\sfGa(*)$ on a big Zariski site $\cC_\Zar$ including homotopy invariance and purity. Concerning such cohomology theory, he developed the general framework of universal Chern classes and characters, which endows with the Chern classes and characters for higher {\it K}-groups
\begin{equation}\label{eq1-1b}
 \sfC_{i,j} : \tK_j(X) \lra \tH^{2i-j}(X,\sfGa(i)), \qquad \ch_X : \tK_*(X) \lra \wh{\tH}{}^*(X,\sfGa(\bullet))_\bQ,
\end{equation}
where $\tK_*(X)$ denotes the algebraic {\it K}-group \cite{Q} and $\wh{\tH}{}^*(X,\sfGa(\bullet))_\bQ$ denotes the direct product of the cohomology groups $\tH^j(X,\sfGa(i)) \otimes \bQ$ for all integers $i$ and $j$. He further extended the formula \eqref{eq1-1a} to this last Chern character.
It is almost forty years since Gillet's paper \cite{Gi} was published, and the \tK-theory of schemes has been much developed by the discovery of the framework of $\bA^1$-homotopy theory, e.g.\ \cite{MV}, \cite{R}, \cite{KY2}.
See also Soul\'e's paper \cite{Sou} for Adams Riemann-Roch for higher \tK-theory, and
also the introduction of \cite{Hd} for a beautiful exposition on the history of Riemann-Roch theorems.
In this paper, we give a new result in a different direction, that is, we extend Gillet's results partially to graded cohomology theories which do {\it not} satisfy homotopy invariance or purity.

\subsection{Setting and results}\label{sect1.1}
Let $\Sch$ be the category consisting of schemes which are separated, noetherian, universally catenary and finite-dimensional, and morphisms of schemes. Let $\cC$ be a subcategory of $\Sch$ satisfying the following two conditions:
\begin{enumerate}
\item[($*_1$)] {\it If $f : Y \to X$ is smooth with $X \in \Ob(\cC)$, then $Y \in \Ob(\cC)$ and $f \in \Mor(\cC)$.}
\item[($*_2$)] {\it If $f : Y \to X$ is a regular closed immersion with $X,Y \in \Ob(\cC)$, then $f \in \Mor(\cC)$.}
\end{enumerate}
We do {\it not} assume that $\cC$ is closed under fiber products.
Let $\sfGa(*)=\{\sfGa(n)\}_{n \in \bZ}$ be a family of cochain complexes of abelian sheaves on the big Zariski site $\cC_\Zar$. Our axioms of {\it admissible cohomology theory} consist mainly of the following three conditions (see Definition \ref{axiom1} below for details):
\begin{enumerate}
\item[(1)] {\it A morphism $\varrho : \Gm[-1] \lra \sfGa(1)$ is given in $\sfD(\cC_\Zar)$.} ({\it From this $\varrho$, one obtains the first Chern class $\sfc_1(L) \in H^2(X_\Zar,\sfGa(1))$ of a line bundle $L$ on $X \in \Ob(\cC)$}.)
\item[(2)] {\it A projective bundle formula for projective bundles in $\cC$.}
\item[(3)] {\it For a strict closed immersion $f : \Ys \hra \Xs$ of codimension $r$ of simplicial objects in $\cC$, push-forward morphisms
\[ f_\push : f_*\sfGa(i)_{\Ys} \lra \sfGa(i+r)_{\Xs}[2r] \]
are given in $\sfD((\Xs)_\Zar)$ and satisfy transitivity, projection formula and compatibility with the first Chern class. Here $\sfD((\Xs)_\Zar)$ denotes the derived category associated with the additive category of unbounded complexes of abelian sheaves on $(\Xs)_\Zar$.}
\end{enumerate}
The conditions (1)\bb(2) have been considered both by Gillet \cite{Gi} Definition 1.2 and Beilinson \cite{Be} \S2.3\,(a)\bb(f). On the other hand, the last condition (3) has been considered only for regular closed immersions of usual schemes in those literatures, which we will need to verify the Whitney sum formula for Chern classes of vector bundles over simplicial schemes, cf.\ \S\ref{sect4} below. See \S\ref{sect3} for fundamental and important examples of admissible cohomology theories. We will define Chern class and character \eqref{eq1-1b} for an admissible cohomology theory $\sfGa(*)$, following the method of Gillet \cite{Gi}.
\par
As for the compatibility of the above axioms (1)\bb(3), the axiom (2) is compatible with (1) in the sense that the first Chern class of a hyperplane has been used in formulating (2). The push-forward morphisms in (3) will be compatible with (1) by assumption. Moreover, we will prove the following compatibility assuming that $\sfGa(*)$ is an admissible cohomology theory, cf.\ Corollary \ref{cor7-2}. Let $E$ be a vector bundle of rank $r$ on $Y \in \Ob(\cC)$ and let $X:=\bP(E\oplus \bs{1})$ be the projective completion of $E$, cf.\ \eqref{eq0-1}. Then for the zero-section $f : Y \to X$, we have
\[ f_\push(1) = \sfc_r(Q) \quad \hbox{ in \;\; $\tH^{2r}(X_\Zar,\sfGa(r))$,} \]
where $1$ denotes the unity of $\tH^0(Y_\Zar,\sfGa(0))$ and $Q$ denotes the universal quotient bundle over $X$. This formula shows a compatibility between the axioms (2) and (3), and plays an important role in our results on Riemann-Roch theorems:
\addtocounter{thm}{-1}
\begin{thm}[${\bs\S}$\ref{sect9}, ${\bs\S}$\ref{sect10}]\label{thm1-1}
Let $\sfGa(*)$ be an admissible cohomology theory on $\cC$, and let $f : Y \to X$ be a projective morphism in $\cC$ with both $X$ and $Y$ regular.
\begin{enumerate}
\item[{\rm(}1{\rm)}]
Assume that $f$ satisfies the assumption {\rm($\#'$)} in Theorem \ref{thm10-1} below. Then the formula \eqref{eq1-1a} holds for $\sfGa(*)$-cohomology, i.e., the following diagram commutes{\rm:}
\[\xymatrix{
\tK_*(Y) \ar[r]^-{f_*} \ar[d]_{\ch_Y(-) \, \cup \, \Td(T_f)} & \tK_*(X) \ar[d]^{\ch_X} \\
\whH^*(Y_\Zar,\sfGa(\scbullet))_\bQ \ar[r]^-{f_\push} & \whH^*(X_\Zar,\sfGa(\scbullet))_\bQ.
}\]
Here $T_f$ denotes the virtual tangent bundle of $f$, cf.\ \S\ref{sect10} below, and $f_\push$ denotes the push-forward morphism that will be constructed in \S\ref{sect7} below.
\item[{\rm(}2{\rm)}]
Assume that $f$ is a {\rm(}regular{\rm)} closed immersion of pure codimension $r \ge 1$ and satisfies the assumption {\rm($\#$)} in Theorem \ref{thm9-1} below. Then the Riemann-Roch theorem without denominators holds for $\sfGa(*)$-cohomology, i.e., the following diagram commutes for any $i,j \ge 0${\rm:}
\[\xymatrix{\tK_j(Y) \ar[r]^-{f_*} \ar[d]_{\sfP_{i-r,Y/X,j}} & \tK_j(X) \ar[d]^{\sfC_{i,j,X}} \\
 \tH^{2(i-r)-j}(Y_\Zar,\sfGa(i-r)) \ar[r]^-{f_\push} & \tH^{2i-j}(X_\Zar,\sfGa(i)), }\]
where $\sfP_{i-r,Y/X,j}$ denotes a mapping class defined by a universal polynomial $\sfP_{i-r,r}$ and universal Chern classes, cf.\ \S\ref{sect9} below.
\end{enumerate}
\end{thm}
The Riemann-Roch theorem without denominators was first raised as a problem in \cite{SGA6} Expos\'e XVI \S3, and proved by Jouanolou and Baum-Fulton-MacPherson for $\tK_0$ (\cite{J} \S1, \cite{BFM} Chapter IV \S5) and by Gillet for the Chern class maps \eqref{eq1-1b} under the assumption that $\sfGa(*)$ satisfies homotopy invariance and purity (\cite{Gi} Theorem 3.1). Theorem \ref{thm1-1}\,(2) removes those assumptions.

\subsection{A logarithmic variant}\label{sect1.2}
Let $X$ be a regular scheme which belongs to $\Ob(\cC)$, and let $D$ be a simple normal crossing divisor on $X$ whose strata belong to $\Ob(\cC)$ (see Definition \ref{defn12-a}\,(2)).
As an application of Theorem \ref{thm1-1}, we will construct Chern class maps
\[ \sfc_{i,(X,D)} : \tK_0(U) \lra \tH^{2i}((X,D)_\Zar,\sfGa(i)) \qquad \hbox{($U:=X \ssm D$)} \]
in \S\ref{sect13}, where $\sfGa(*)=\sfGa(*)^\log$ denotes an admissible cohomology theory on the category of log pairs in $\cC$, cf.\ Definitions \ref{defn12-a} and \ref{defn12-1} below. If $D$ is empty, then $\sfc_{i,(X,\emptyset)}$ is the $i$-th Chern class map of $X$ with values in an admissible cohomology in the sense of \S\ref{sect1.1}. Following the idea of Somekawa in \cite{So} Chapter \II, we construct the map $\sfc_{i,(X,D)}$ by induction on the number of irreducible components of $D$. A key point is to prove the Riemann-Roch theorem without denominators, analogous to Theorem \ref{thm1-1}\,(2), for this new Chern class map $\sfc_{i,(X,D)}$ to proceed the induction step. We have to note that a tensor product formula (cf.\ (L4) in \S\ref{sect13}) is necessary for this argument. Because we consider only $\tK_0$, one can derive this formula easily from that for the case $D=\emptyset$ and the surjectivity of the map $\tK_0(X) \to \tK_0(U)$. 

\par\medskip
This paper is organized as follows. In \S\ref{sect2} we will formulate admissible cohomology theory, whose examples will be explained in \S\ref{sect3}. We will construct Chern classes of vector bundles, universal Chern class and character, higher Chern class and character, following the method of Grothendieck and Gillet in \S\S\ref{sect4}\bb\ref{sect6} below. The section \ref{sect7} will be devoted to extending push-forward morphisms to projective morphisms in $\cC$, which plays a key role in our proof of Riemann-Roch theorems. We will give an explicit construction of Jouanolou's universal polynomial in \S\ref{sect8} for the convenience of the reader. After those preliminaries, we will prove Riemann-Roch theorems in \S\S\ref{sect9}\bb\ref{sect10}, and give an application in \S\ref{sect11}. In \S\S\ref{sect12}\bb\ref{sect13}, we will formulate a logarithmic variant of admissible cohomology theory and construct a Chern class map on $\tK_0$ with values in the admissible cohomology with log poles, which is another application of the Riemann-Roch theorem in \S\ref{sect9}. 

\par\medskip
The authors express their gratitude to Kei Hagihara and Satoshi Mochizuki for valuable comments and suggestions.

\subsection{Notation and conventions}\label{notation}
In this paper, all schemes are assumed to be {\it separated}, {\it noetherian of finite dimension} and {\it universally catenary}. Unless indicated otherwise, all cohomology groups of schemes are taken over the Zariski topology.
\par
For a scheme $X$, a closed subset $Z \subset X$ and a cochain complex $\cF^\bullet$ of abelian sheaves on $X_\Zar$, we define the hypercohomology group $\tH^i_Z(X,\cF^\bullet)$ with support in $Z$ as the $i$-th cohomology group of the complex $\vG_Z(X,\cJ^\bullet)$, where $\cJ^\bullet$ denotes an injectively fibrant resolution of $\cF^\bullet$, cf.\ \S\ref{sectA.0} below. 
\par
{\it A projective morphism} $f : Y \to X$ of schemes means a morphism which factors as follows for some integer $n \ge 0$:
\[\xymatrix{  Y \; \ar@<-1pt>@{^{(}->}[r]^i & \bP^n \ar[r]^p & X, }\]
where $i$ is a closed immersion and $p$ is the natural projection, cf.\ \cite{Ha} p.\ 103. When $X$ is regular, a projective morphism $f : Y \to X$ in the sense of \cite{EGAII} 5.5.2 is projective in our sense by the existence of an ample family of line bundles over $X$, cf.\ \cite{SGA6} Expos\'e \II\, Corollaire 2.2.7.1.
\par
For a vector bundle $E$ over a scheme $X$, we define the projective bundle $\bP(E)$ as
\addtocounter{thm}{2}
\begin{equation}\label{eq0-1}
 \bP(E):=\bProj(\Sym_{\cO_X}^\bullet(\cE^\vee)),
\end{equation}
where $\cE$ denotes the locally free sheaf on $X$ represented by $E$ and $\cE^\vee$ means its dual sheaf over $\cO_X$. We define the tautological line bundle $L^\taut$ over $\bP(E)$ as follows:
\begin{equation}\label{eq0-2}
 L^\taut:=\bSpec(\Sym_{\cO_{\bP(E)}}^\bullet(\cO(-1))),
\end{equation}
where $\cO(-1)$ denotes the $\cO_{\bP(E)}$-dual of the twisting sheaf $\cO(1)$ of Serre. For a Cartier divisor $D$ on $X$, {\it the line bundle over $X$ associated to $D$} means the line bundle
\begin{equation}\label{eq0-3}
 \bSpec(\Sym_{\cO_X}^\bullet(\cO_X(-D))),
\end{equation}
which represents the invertible sheaf $\cO_X(D)$ on $X$.
\par
Let $\vD$ be the simplex category, whose objects are ordered finite sets
\[ [p]:=\{0,1,2,\dotsc,p\}  \qquad \hbox{($p \ge 0$)} \]
and whose morphisms are order-preserving maps.

\begin{defn}\label{def1-1}
{\rm Let $\cB$ be a category.
\begin{enumerate}
\item[(1)]
{\it A simplicial object} $\Xs$ in $\cB$ is a functor
\[ \Xs : \vD^\op \lra \cB. \]
{\it A morphism $f : \Ys \to \Xs$ of simplicial objects in $\cB$} is a natural transform of such contravariant functors.
\item[(2)]
For a simplicial object $\Xs$ in $\cB$ and a morphism $\alpha : [p] \to [q]$ in $\vD$, we often write 
\[ \alpha^X : X_q \lra X_p \qquad \hbox{($X_p:=\Xs([p])$)} \]
for $\Xs(\alpha)$, which is a morphism in $\cB$.
\end{enumerate}}
\end{defn}

\begin{defn}\label{def1-2}
{\rm Let $\Xs$ be a simplicial scheme. 
\begin{enumerate}
\item[(1)]
{\it A vector bundle} over $\Xs$ is a morphism $f : \Es \to \Xs$ of simplicial schemes such that $f_p : E_p \to X_p$ is a vector bundle for any $p \ge 0$ and such that the commutative diagram
\begin{equation}\notag
\xymatrix{ E_q \ar[r]^{f_q} \, \ar@<-1pt>[d]_{\alpha^E} & X_q \ar[d]^{\alpha^X} \\
 E_p \ar[r]^{f_p} \, & X_p}
\end{equation}
induces an isomorphism $E_q \cong \alpha^{X*}\!E_p:=E_p \times_{X_p} X_q$ of vector bundles over $X_q$ for any morphism $\alpha : [p] \to [q]$ in $\vD$\; (cf.\ \cite{Gi2} Example 1.1).
\item[(2)]
We say that a morphism $f : \Ys \to \Xs$ of simplicial schemes is {\it a closed immersion} if $f_p : Y_p \to X_p$ is a closed immersion for each $p \ge 0$.  We say that a closed immersion $f : \Ys \to \Xs$ is {\it exact}, if the diagram
\begin{equation}\label{eq1-1}
\xymatrix{ Y_q \ar[r]^{f_q} \, \ar@<-1pt>[d]_{\alpha^Y} & X_q \ar[d]^{\alpha^X} \\
 Y_p \ar[r]^{f_p} \, & X_p}
\end{equation}
is cartesian for any morphism $\alpha : [p] \to [q]$ in $\vD$. We say that a closed immersion $f : \Ys \to \Xs$ is {\it strict} if it is exact and $f_p : Y_p \to X_p$ is regular for each $p \ge 0$. {\it An effective Cartier divisor $\Xs$ on $\Ys$} is a strict closed immersion $\Xs \to \Ys$ of pure codimension $1$.
\end{enumerate}}
\end{defn}

\section{Admissible cohomology theory}\label{sect2}
The aim of this section is to formulate the axioms of admissible cohomology theory in Definition \ref{axiom1} below.
\begin{defn}[{{\it Graded cohomology theory}}]\label{axiom0}
{\rm Let $\cS$ be a site, and let $\sfGa(*)=\{\sfGa(i)\}_{i \in \bZ}$ be a family of complexes of abelian sheaves on $\cS$. We say that $\sfGa(*)$ is {\it a graded cohomology theory on $\cS$}, if it satisfies the following two conditions (cf.\ \cite{Gi} Definition 1.1):
\begin{enumerate}
\item[(a)]
$\sfGa(0)$ is concentrated in degrees $\ge 0$, and the $0$-th cohomology sheaf $\cH^0(\sfGa(0))$ is a sheaf of commutative rings with unity.
\item[(b)]
$\sfGa(*)$ is equipped with an associative and commutative product structure
\[ \sfGa(i) \otimes^{\bL} \sfGa(j) \lra \sfGa(i+j) \quad \hbox{ in } \;\; \sfD(\sfS^\ab(\cS)) \]
compatible with the product structure on $\cH^0(\sfGa(0))$ stated in (a).
\end{enumerate}}
\end{defn}

Now let $\cC$ be as in \S\ref{sect1.1}. For a simplicial object $\Xs$ in $\cC$, there is a natural restriction functor on the category of abelian sheaves
\[ \theta_{\!\Xs} : \sfS^\ab(\cC_\Zar) \lra \sfS^\ab((\Xs)_\Zar), \]
which sends a sheaf $\cF$ on $\cC_\Zar$ to the sheaf $(U \subset X_p) \mapsto \cF(U)$ on $(\Xs)_\Zar$.
This functor is exact and extends naturally to a triangulated functor on derived categories
\[ \theta_{\!\Xs} : \sfD(\cC_\Zar) \lra \sfD((\Xs)_\Zar). \]
\begin{defn}\label{def2-a}
{\rm Let $\sfGa(*)$ be a graded cohomology theory on $\cC_\Zar$ and let $\Xs$ be a simplicial object in $\cC$. For each $i \in \bZ$, we define a complex $\sfGa(i)_{\Xs}$ of abelian sheaves on $(\Xs)_\Zar$ by applying $\theta_{\Xs}$ to the complex $\sfGa(i)$. We will often omit the indication of the functor $\theta_{\!\Xs}$ in what follows.}
\end{defn}

\begin{defn}[{{\it First Chern class}}]\label{def2-b}
{\rm Let $\sfGa(*)$ be a graded cohomology theory on $\cC_\Zar$, and suppose that we are given a morphism
\[ \varrho : \cO^\times[-1] \lra \sfGa(1) \qquad \hbox{ in } \quad \sfD(\cC_\Zar), \]
where $\cO^\times$ means the abelian sheaf on $\cC_\Zar$ represented by the group scheme $\Gm$. Let $\Xs$ be a simplicial object in $\cC$, and let $\Ls$ be a line bundle over $\Xs$. There is a class $[\Ls] \in \tH^1(\Xs,\cO^\times)$ corresponding to $\Ls$ (cf.\ \cite{Gi2} Example 1.1).
We define {\it the first Chern class} $\sfc_1(\Ls) \in \tH^2(\Xs,\sfGa(1))$ as the value of $[\Ls]$ under the map
\[ \tH^1(\Xs,\cO^\times) \os{\varrho}\lra \tH^2(\Xs,\sfGa(1)). \]
}
\end{defn}
The first Chern classes are functorial in the following sense:
\begin{lem}\label{rem2-1}
Let $\sfGa(*)$ be a graded cohomology theory on $\cC_\Zar$, and suppose that we are given a morphism $\varrho : \cO^\times[-1] \to \sfGa(1)$ in $\sfD(\cC_\Zar)$. Then for a morphism $f : \Ys \to \Xs$ of simplicial objects in $\cC$ and a line bundle $\Ls$ over $\Xs$, we have
\[ \sfc_1(\fLs) = \fsfc_1(\Ls) \quad \hbox{ in } \;\; \tH^2(\Ys,\sfGa(1)). \]
Here $\fLs$ denotes $\Ls \times_{\Xs}\! \Ys$\,, the inverse image of $\Es$ by $f$.
\end{lem}
\begin{pf}
The assertion is obvious, because $\sfGa(1)$ and $\varrho$ are defined on the big site $\cC_\Zar$.
\end{pf}

\begin{defn}[{{\it Admissible cohomology theory}}]\label{axiom1}
{\rm We say that a graded cohomology theory $\sfGa(*)$ on $\cC_\Zar$ is {\it an admissible cohomology theory on $\cC$}, if it satisfies the axioms (1)\bb(3) below. Compare with \cite{Be} \S2.3\,(a)\bb(f), \cite{Gi} Definition 1.2.
\begin{enumerate}
\item[(0)]
For a scheme $X \in \Ob(\cC)$ and a dense open subset $U$ of $X$, the restriction map $\tH^0(X,\sfGa(0)) \to \tH^0(U,\sfGa(0))$ is injective.
\item[(1)]
({\it First Chern class}) \;
There exists a morphism
\[ \varrho : \cO^\times[-1] \lra \sfGa(1) \qquad \hbox{ in } \quad \sfD(\cC_\Zar). \]
\item[(2)]
({\it Projective bundle formula}) \;
For a scheme $X \in \Ob(\cC)$ and a vector bundle $E$ over $X$ of rank $r+1$, the morphism
\begin{align*}
 \gamma_E : \;& \bigoplus_{j=0}^r \ \sfGa(i-j)_X[-2j] \lra \tR p_*\sfGa(i)_{\bP(E)}, \quad
 (x_j)_{j=0}^r \mapsto \sum_{j=0}^r \  \xi^j \cup p^\back(x_j)
\end{align*}
is an isomorphism in $\sfD(X_\Zar)$. Here $p : \bP(E) \to X$ denotes the projective bundle associated with $E$, cf.\ \eqref{eq0-1}, and $\xi \in \tH^2(\bP(E),\sfGa(1))$ denotes the first Chern class of the tautological line bundle, cf.\ Definition \ref{def2-b}. See \eqref{eq0-2} for the definition of the tautological line bundle.
\item[(3)]
({\it Push-forward for strict closed immersions}) \;
For a strict closed immersion $f : \Ys \hra \Xs$ of simplicial objects in $\cC$ of pure codimension $r$, there are push-forward morphisms
\[ f_\push :  f_*\sfGa(i)_{\Ys} \lra \sfGa(i+r)_{\Xs}[2r] \qquad \hbox{($i \in \bZ$)} \]
in $\sfD((\Xs)_\Zar)$ which satisfy the following four properties.
\begin{enumerate}
\item[(3a)]
({\it Consistency with the first Chern class}) \;
When $r=1$, the push-forward map
\[ f_\push : \tH^0(\Ys, \sfGa(0)) \lra \tH^2(\Xs, \sfGa(1)) \]
sends $1$ to the first Chern class of the line bundle over $\Xs$ associated with $\Ys$, cf.\ Definition \ref{def2-b}.
\item[(3b)]
({\it Projection formula}) \;
The following diagram commutes in $\sfD((\Xs)_\Zar)${\rm:}
\[\hspace{-10pt} \xymatrix{ \sfGa(i)_{\Xs} \os{\bL}\otimes f_*\sfGa(j)_{\Ys} \ar[r]^-{\id \otimes f_\push} \ar[d]_{f^\back \otimes \id} & \sfGa(i)_{\Xs} \os{\bL}\otimes \sfGa(j+r)_{\Xs}[2r] \ar[rd]^-{\text{product}} \\
 f_*\sfGa(i)_{\Ys} \os{\bL}\otimes f_*\sfGa(j)_{\Ys} \ar[r]^-{\text{product}} &  f_*\sfGa(i+j)_{\Ys} \ar[r]^-{f_\push} &  \sfGa(i+j+r)_{\Xs}[2r]\, .} \]
\item[(3c)]
({\it Transitivity}) \;
For closed immersions $f : Y \hra X$ and $g : Z \hra Y$ of objects in $\cC$ of pure codimension $r$ and $r'$, respectively, the composite morphism
\begin{align*} (f \circ g)_*\sfGa(i)_{Z} & \cong f_*g_*\sfGa(i)_{Z} 
 \os{g_\push}\lra f_*\sfGa(i+r')_{Y}[2r']  \os{f_\push}\lra \sfGa(i+r+r')_{X}[2(r+r')] \end{align*}
agrees with $(f \circ g)_\push$\,.
\item[(3d)]
({\it Base-change property}) \;
Let 
\[\xymatrix{ Y' \; \ar@<-1pt>@{^{(}->}[r]^{f'} \ar@<-2.5pt>[d]_{g'} \ar@{}[rd]|{\square} & X' \ar[d]^g \\ Y \; \ar@<-1pt>@{^{(}->}[r]^f & X }\]
be a diagram in $\cC$ which is cartesian in the category of schemes and such that $f$ and $f'$ are strict closed immersions of pure codimension $r$.
Then the following diagram commutes in $\sfD(X_\Zar)$:
\def\uG{\ul{\vG}{\,}}
\[\xymatrix{ \tR (g \circ f')_*\sfGa(i)_{Y'} \ar[r]^-{f'_\push} & \tR g_*\sfGa(i+r)_{X'}[2r]
 \\ f_*\sfGa(i)_{Y} \ar[r]^-{f_\push} \ar[u]^{Rf_*(g'^\back)} & \sfGa(i+r)_{X}[2r] \ar[u]_{g^\back}\,. } \]
\end{enumerate}
\end{enumerate}}
\end{defn}

\begin{rem}\label{rem-axiom}
\begin{enumerate}
\item[{\rm(1)}]
In {\rm\cite{Be}} and {\rm\cite{Gi}}, the axiom {\rm(}3{\rm)} is considered only for usual schemes. In fact, if $\sfGa(*)$ satisfies homotopy invariance, then we need push-forward maps only for usual schemes to verify Theorem \ref{thm4-1}\,{\rm(}3{\rm)} below, cf.\ {\rm\cite{Le}} Part I, Chapter \III\, \S1.3.3.
\item[{\rm(2)}]
The axiom {\rm(}3{\rm)} implies the weak Gysin property in {\rm\cite{SCH}} p.\ 20 for simplicial schemes.
\item[{\rm(3)}]
We need the properties {\rm(}3c{\rm)} and {\rm(}3d{\rm)} in the axiom {\rm(}3{\rm)} only for usual schemes, to verify Theorem \ref{thm1-1}. A key step is to extend the push-forward morphisms to those for regular projective morphisms. See \S\ref{sect7} below for details.
\item[{\rm(4)}]
We do not assume that the push-forward map in the axiom {\rm(}3{\rm)} is an isomorphism, because the projective bundle formula and the full purity imply the homotopy invariance. Compare with the purity in the sense of {\rm\cite{Gi}} Definition 1.2\,{\rm(}vi{\rm)}.
\item[{\rm(5)}]
The axiom {\rm(}0{\rm)} is a technical one, but will be useful in our construction of push-forward morphisms for regular projective morphisms. See the proof of Lemma \ref{lem7-1}\,{\rm(}1{\rm)} below.
\end{enumerate}
\end{rem}

\section{Examples of admissible cohomology theory}\label{sect3}
We give several fundamental and important examples of $\cC$ and $\sfGa(*)$. The first four examples satisfy homotopy invariance and purity, while the others do neither of them.

\subsection{Motivic complex}\label{ex:motivic}
Let $k$ be a field, and let $\cC$ be the full subcategory of $\Sch/k$ consisting of schemes which are smooth separated of finite type over $k$. For $i \in \bZ$, we define $\sfGa(i)$ on $\cC_\Zar$ as follows:
\[ \sfGa(i) := \begin{cases} \bZ(i) \qquad & (i \ge 0), \\ 0 & (i < 0), \end{cases} \]
where $\bZ(i)=\rC\hspace{1pt}^\bullet(\bZ_\tr(\Gm^{\wedge i}))[-i]$ denotes the motivic complex of Suslin-Voevodsky \cite{SV} Definition 3.1, and $\rC\hspace{1pt}^\bullet$ denotes the singular complex construction due to Suslin.
We will prove that $\sfGa(*)$ is an admissible cohomology theory in Appendix \ref{appA} below, assuming that {\it $k$ admits the resolution of singularities} in the sense of \cite{SV} Definition 0.1.

\subsection{\'Etale Tate twist}\label{ex:etale}
Let $n$ be a positive integer, and let $\cC$ be the full subcategory of $\Sch$ consisting of regular schemes over $\Spec(\bZ[n^{-1}])$. For $i \in \bZ$, we define $\sfGa(i)$ on $\cC_\Zar$ as follows:
\[ \sfGa(i) := \begin{cases} \tR \vare_*\mu_n^{\otimes i} \qquad & (i \ge 0), \\ \tR \vare_* \big(\sHom(\mu_n^{\otimes (-i)},\Zn)\big) \qquad & (i < 0), \end{cases} \]
where $\mu_n$ denotes the \'etale sheaf of $n$-th roots of unity and $\vare : \cC_\et \to \cC_\Zar$ denotes the natural morphism of sites. Obviously $\sfGa(*)=\{\sfGa(i)\}_{i \in \bZ}$ is a graded cohomology theory on $\cC$. We define the morphism $\varrho$ in \ref{axiom1}\,(1) by the connecting morphism associated with the Kummer exact sequence on $\cC_\et$
\[ 0 \lra \sfGa(1) \lra \cO^\times \os{\times n}\lra \cO^\times \lra 0. \]
The property \ref{axiom1}\,(2) follows from the homotopy invariance (\cite{SGA4} Expos\'e XV, Th\'eor\`eme 2.1)  and the relative smooth purity.
\par
We check that $\sfGa(*)$ satisfies the axiom \ref{axiom1}\,(3), in what follows. Let $f : \Ys \hra \Xs$ be a strict closed immersion of pure codimension $r$ of simplicial schemes in $\cC$. We use Gabber's refined cycle class \cite{FG} Definition 1.1.2
\[ \cl_{X_0}(Y_0) \in \tH^{2r}_{Y_0}(X_0,\sfGa(r)). \]
By the spectral sequence (cf.\ \cite{C} Proposition \II.2)
\[ E_1^{a,b}=H^b_{Y_a}(X_a,\sfGa(r)) \Lra H^{a+b}_{\Ys}(\Xs,\sfGa(r)) \]
and the semi-purity in \cite{FG} \S8, we have
\[ H^{2r}_{\Ys}(\Xs,\sfGa(r)) \cong \ker(d_0^*-d_1^* : H^{2r}_{Y_0}(X_0,\sfGa(r)) \to H^{2r}_{Y_1}(X_1,\sfGa(r))). \]
By the functoriality in loc.\ cit.\ Proposition 1.1.3 and the assumption that the square \eqref{eq1-1} is cartesian, we have
\[ d_0^*\cl_{X_0}(Y_0)=\cl_{X_1}(Y_1)=d_1^*\cl_{X_0}(Y_0), \]
and consequently, $\cl_{X_0}(Y_0)$ belongs to $\ker(d_0^*-d_1^*)$. We thus obtain a cycle class
\[ \cl_{\Xs}(\Ys) \in H^{2r}_{\Ys}(\Xs,\sfGa(r)) \]
as the element corresponding to $\cl_{X_0}(Y_0)$. Since $f^*\sfGa(i)_{\Xs} \cong \sfGa(i)_{\Ys}$ on $(\Ys)_\Zar$\,, the cup product with $\cl_{\Xs}(\Ys)$ defines the desired push-forward morphism
\[\xymatrix{ f_\push :  f_*\sfGa(i)_{\Ys} \cong f_*f^*\sfGa(i)_{\Xs} \ar[rr]^-{\cl_{\Xs}(\Ys) \,\cup \, -} && \sfGa(i+r)_{\Xs}[2r]  \quad \hbox{ in } \;\; D((\Xs)_\Zar),}\]
which satisfies the properties \ref{axiom1}\,(3a), (3b). See loc.\ cit.\ Proposition 1.2.1 for (3c). The property (3d) follows from loc.\ cit.\ Proposition 1.1.3. Thus $\sfGa(*)$ is an admissible cohomology theory on $\cC$.

\subsection{Betti complex}\label{ex:betti}
Let $\cC$ be the full subcategory of $\Sch/\bC$ consisting of schemes which are smooth separated of finite type over $\bC$. Let $\cC_\an$ be the big analytic site associated with $\cC$. Let $A$ be a subring of $\bR$ with unity. For $i \in \bZ$, we define $\sfGa(i)$ on $\cC_\Zar$ as follows:
\[ \sfGa(i) := \tR \vare_* \big( (2\pi \sqrt{-1})^i A \big),  \]
where $\vare : \cC_\an \to \cC_\Zar$ denotes the natural morphism of sites. When $A=\bZ$, we define the morphism $\varrho=\varrho_\bZ$ in \ref{axiom1}\,(1) as the connecting morphism of the exponential exact sequence
\[ 0 \lra 2\pi \sqrt{-1} \cdot \bZ \lra \cO \os{\exp}\lra \cO^\times \lra 0.  \]
For a general $A$, we define $\varrho=\varrho_A$ as the composite morphism
\[ \varrho : \cO^\times[-1] \os{\varrho_\bZ}\lra 2\pi \sqrt{-1} \cdot \bZ \hra 2\pi \sqrt{-1} \cdot A. \]
The axioms \ref{axiom1}\,(2)\bb(4) can be checked in a similar way as for \S\ref{ex:etale}.

\subsection{Deligne-Beilinson complex}\label{ex:deligne}
Let $\cC$ be as in \S\ref{ex:betti}. Let $A$ be a subring of $\bR$ with unity. For $i \in \bZ$, we define $\sfGa(i)$ on $\cC_\Zar$ as follows:
\[ \sfGa(i) := \begin{cases} \sfGad(i) \qquad & (i \ge 0), \\ \tR \vare_* \big( (2\pi \sqrt{-1})^i A \big) \qquad & (i < 0), \end{cases} \]
where $\sfGad(i)$ denotes the Deligne-Beilinson complex on $\cC_\Zar$ in the sense of \cite{EV} Theorem 5.5 and $\vare : \cC_\an \to \cC_\Zar$ denotes the natural morphism of sites, cf.\ \S\ref{ex:betti}. See loc.\ cit., Theorem 5.5\,(b) and Proposition 8.5 (resp.\ \cite{Ja} \S3.2) for the axiom \ref{axiom1}\,(1) and (2) (resp.\ \ref{axiom1}\,(3) and (4)).

\subsection{Algebraic de Rham complex}\label{ex:deRham}
Let $k$ be a field, and let $\cC$ be the full subcategory of $\Sch/k$ consisting of schemes which are smooth separated of finite type over $k$. For $i \in \bZ$, we define $\sfGa(i)$ on $\cC_\Zar$ as the de Rham complex $\Omega^\bullet_{-/k}$ over $k$.
We see that $\sfGa(*)$ is an admissible cohomology theory on $\cC$, when we define $\varrho$ in \ref{axiom1}\,(1) by logarithmic differentials. See \cite{Ha1} Chapter \II\, \S2 for the axioms \ref{axiom1}\,(3).

\subsection{Logarithmic Hodge-Witt sheaf}\label{ex:loghw}
Let $p$ be a prime number, and let $\cC$ be the full subcategory of $\Sch$ consisting of regular schemes over $\Spec(\bF_p)$. Let $n$ be a positive integer. For $i \ge 0$, we define $\sfGa(i)$ on $\cC_\Zar$ as follows:
\[ \sfGa(i) := \begin{cases} \tR \vare_*\logwittt n i [-i] \qquad & (i \ge 0), \\ 0 \qquad & (i < 0), \end{cases} \]
where $\logwittt n i$ denotes the \'etale subsheaf of the logarithmic part of the Hodge-Witt sheaf $\witt {} n i$ on $\cC_\et$, cf.\ \cite{Il}, and $\vare : \cC_\et \to \cC_\Zar$ denotes the natural morphism of sites. Then $\sfGa(*)=\{\sfGa(i)\}_{i \in \bZ}$ is an admissible cohomology theory on $\cC$, which we are going to check. We define the morphism $\varrho$ in \ref{axiom1}\,(1) as the logarithmic differential map. See \cite{Gros} Chapter I\, Th\'eor\`eme 2.1.11 and \cite{Sh} Theorems 2.1, 2.2 for the axiom (2). To verify the axiom (3), we construct the push-forward map for a strict closed immersion $f : \Ys \hra \Xs$ of pure codimension $r$ of simplicial schemes in $\cC$. For $a \ge 0$, let $f_a : Y_a \hra X_a$ be the $a$-th factor of $f$. Note first that
\[ \tR^j f_a^!\sfGa(i+r)_{X_a} \cong \begin{cases} 0  \quad & \hbox{($j < i+2r$)} \\ \cH^i(\sfGa(i)_{Y_a}) \quad & \hbox{($j = i+2r$)} \end{cases} \]
by loc.\ cit.\ Theorem 3.2 and Corollary 3.4, which immediately implies
\[ \tR^j f^!\sfGa(i+r)_{\Xs} = 0 \quad \hbox{for } \;\; j < i+2r. \]
We denote the above isomorphism for $j=i+2r$ by $(f_a)_\push$\,. To show that the maps $(f_a)_\push$ for $a \ge 0$ give rise to an isomorphism
\[ \cH^i(\sfGa(i)_{\Ys}) \cong \tR^{i+2r}f^!\sfGa(i+r)_{\Xs}, \]
it is enough to check that the maps $(f_a)_\push$ are compatible with the simplicial structures of $\Xs$ and $\Ys$\,.
One can easily check this by the local description of $(f_a)_\push$ in \cite{Sh} p.\ 589 and the assumption that the square \eqref{eq1-1} is cartesian. Thus we obtain a morphism
\[\xymatrix{ f_\push : f_*\sfGa(i)_{\Ys} \ar[r] & \sfGa(i+r)_{\Xs}[2r]  \quad \hbox{ in } \;\; \sfD((\Xs)_\Zar). }\]
The properties (3a)\bb(3d) follows again from the local description in \cite{Sh} p.\ 589.

\subsection{$\bs{p}$-adic \'etale Tate twist}\label{ex:ptate}
Let $B$ be a Dedekind ring of mixed characteristics, and put $S:=\Spec(B)$. Let $p$ be a prime number which is not invertible on $S$, and let $\cC$ be the full subcategory of $\Sch/S$ consisting of regular schemes $X$ which are flat of finite type over $S$ and satisfy the following condition:
\begin{itemize}
\item {\it Let $B'$ be the integral closure of $B$ in $\vG(X,\cO_X)$. Then for any closed point $x$ on $\Spec(B')$ with $\ch(x)=p$, the fiber $X \times_{\Spec(B')} x$ is a reduced divisor with normal crossings on $X$.}
\end{itemize}
Fix a positive integer $n$. For $i \in \bZ$, we define
\[ \sfGa(i) := \tR \vare_*\fT_n(i), \]
where $\vare$ denotes the natural morphism of sites $\cC_\et \to \cC_\Zar$, and $\fT_n(i)$ denotes the $i$-th \'etale Tate twist with $\Zpn$-coefficients \cite{Sa2} Definition 3.5, a bounded complex of sheaves on $\cC_\et$. Then $\sfGa(*)= \{\sfGa(i)\}_{i \in \bZ}$ is an admissible cohomology theory on $\cC$. See \cite{Sa2} Theorem 4.1 and Proposition 5.5 for the axioms \ref{axiom1}\,(2), (3a)\bb(3c). The property (3d) follows from the construction of the push-forward morphisms given there and the corresponding property in \S\ref{ex:etale} above.


\section{Chern class of vector bundles}\label{sect4}
In this section we define Chern classes of vector bundles over simplicial schemes following the method of Grothendieck and Gillet (cf.\ \cite{Gr}\, p.\ 144 Theorem 1, \cite{Gi}\,Definition 2.10), and prove Theorem \ref{thm4-1} below. Let $\cC$ be as in \S\ref{sect1.1}, and let $\sfGa(*)$ be an admissible cohomology theory on $\cC$.
\begin{defn}[{{\it Chern class}}]\label{def4-1}
{\rm For a simplicial object $\Xs$ in $\cC$, we put
\[ \tH^{2*}(\Xs,\sfGa(*)) := \bigoplus_{i \ge 0} \ \tH^{2i}(\Xs,\sfGa(i)), \]
which is a commutative ring with unity by the axioms \ref{axiom0}\,(a), (b). For a vector bundle $\Es$ over $\Xs$, we define the {\it total Chern class of $\Es$}
\[ \sfc(\Es)=(\sfc_i(\Es))_{i \ge 0} \in \tH^{2*}(\Xs,\sfGa(*)) \]
as follows. Let $\Es$ be of rank $r$, and let $p$ be the natural projection $\bP(\Es) \to \Xs$. Let $\Ls^\taut$ be the tautological line bundle over $\bP(\Es)$ and put $\xi:=\sfc_1(\Ls^\taut) \in \tH^2(\bP(\Es),\sfGa(1))$.
There is an isomorphism
\begin{equation}\label{eq4-1}
 \bigoplus_{i=1}^r \ \tH^{2i}(\Xs,\sfGa(i)) \cong \tH^{2r}(\bP(\Es),\sfGa(r)), \quad (b_i)_{i=1}^r \mapsto \sum_{i=1}^r \ \xi^{r-i} \cup p^\back(b_i)
\end{equation}
by the axiom \ref{axiom1}\,(2), cf.\ \cite{Gi} Lemma 2.4. We define $\sfc_0(\Es):=1 \in \tH^0(\Xs,\sfGa(0))$ and define $(\sfc_1(\Es),\sfc_2(\Es),\dotsc,\sfc_r(\Es))$ as the unique solution $(\sfc_1,\sfc_2,\dotsc,\sfc_r)$ to the equation
\[ \xi^r + \xi^{r-1} \cup p^\back(\sfc_1) + \dotsb + \xi \cup p^\back(\sfc_{r-1}) + p^\back(\sfc_r) = 0 \]
in $\tH^{2r}(\bP(\Es),\sfGa(r))$. We define $\sfc_i(\Es):=0$ for $i > r$.}
\end{defn}
\begin{thm}\label{thm4-1}
Let $\Xs$ be a simplicial object in $\cC$.
\begin{enumerate}
\item[{\rm(1)}] {\rm(}Normalization{\rm)}\;
We have $\sfc_0(\Es)=1$ and $\sfc_i(\Es)=0$ for $i > \rank(\Es)$. If $\Es$ is a line bundle, then $\sfc_1(\Es)$ defined here agrees with the first Chern class in Definition \ref{def2-a}\,{\rm(}2{\rm)}.
\item[{\rm(2)}] {\rm(}Functoriality{\rm)}\;
For a morphism $f : \Ys \to \Xs$ of simplicial objects in $\cC$ and a vector bundle $\Es$ over $\Xs$, we have
\[ \sfc(\fEs)=\fsfc(\Es), \]
where $\fEs$ denotes $\Es \times_{\Xs}\! \Ys$\,, the inverse image of $\Es$ by $f$.
\item[{\rm(3)}] {\rm(}Whitney sum{\rm)}\;
For a short exact sequence $0 \to \Es' \to \Es \to \Es'' \to 0$ of vector bundles over $\Xs$, we have
\[ \sfc(\Es)=\sfc(\Es') \cup \sfc(\Es'') \quad \hbox{ in }\;\; \tH^{2*}(\Xs,\sfGa(*)). \]
\item[{\rm(4)}] {\rm(}Tensor product{\rm)}\;
For vector bundles $\Es$ and $\Es'$ over $\Xs$, we have
\begin{align*}
\hspace{15pt} \wt{\sfc}(\Es \otimes \Es') & = \wt{\sfc}(\Es) \,\scstar\; \wt{\sfc}(\Es') \\
 &\hspace{30pt} \hbox{ in }\;\; \wtH^{2*}(\Xs,\sfGa(*)):=\bZ \times \{ 1 \} \times I,
\end{align*}
where $\wt{\sfc}(\Es)$ denotes the augmented total Chern class $(\rk(\Es),\sfc(\Es))$, and $I$ denotes the positive part of the graded commutative ring $\tH^{2*}(\Xs,\sfGa(*))${\rm:}
\[ I = \bigoplus_{i > 0} \ \tH^{2i}(\Xs,\sfGa(i)). \]
We endowed $\wtH^{2*}(\Xs,\sfGa(*))$ with the $\lam$-ring structure associated with $\tH^{2*}(\Xs,\sfGa(*))$ {\rm(\cite{Gr2}} Chapter I \S3{\rm)}, and wrote $\scstar$ for its product structure.
\item[{\rm(5)}] The Chern classes $\sfc(\Es)$ are characterized by the properties {\rm(}1{\rm)}--{\rm(}3{\rm)}.
\end{enumerate}
\end{thm}

The properties (1) and (2) immediately follow from the definition of Chern classes and the functoriality in Remark \ref{rem2-1}. The property (4) follows from (3) and the splitting principle of vector bundles. The assertion (5) also follows from the splitting principle of vector bundles. The most important part of this theorem is the verification of the property (3), which is reduced to the Lemma \ref{lem4-1} below, by Grothendieck's arguments in \cite{Gr} Proof of Th\'eor\`eme 1.
In this paper, we cannot use the splitting bundle argument (cf. \cite{Le} Part I, Chapter \III\, \S1.3.3, \cite{ILO} Expos\'e XVI Proposition 1.5) because we do not assume homotopy invariance:

\begin{lem}\label{lem4-1}
Let $\Es$ be a vector bundle of rank $r$ on $\Xs$. Let $p : \bP(\Es) \to \Xs$ be the projective bundle associated with $\Es$, and let $L^\taut_\star$ be the tautological line bundle over $\bP(\Es)$. Suppose that we are given a filtration on $\Es$ by subbundles
\[ \Es=\Es^0 \supset \Es^1 \supset \dotsb \supset \Es^r = 0 \qquad \hbox{{\rm($r:=\rank(\Es)$)}}\]
such that the quotient $\Es^i/\Es^{i+1}$ is a line bundle for $0 \le i \le r-1$. Then we have
\[ \prod_{i=0}^{r-1} \ (\sfc_1(\Ls^\taut)+\psfc_1(\Es^i\hspace{1pt}/\Es^{i+1})) = 0
 \qquad \hbox{ in } \quad \tH^{2r}(\bP(\Es),\sfGa(r)). \]
\end{lem}
\begin{pf}
Put $\Fs:=\pEs \otimes \Ls^\taut$, and let $s : \bP(\Es) \to \Fs$ be the composite morphism
\[ s : \bP(\Es) \os{1}\lra \bA^1_{\bP(\Es)} \lra \pEs \otimes \Ls^\taut=\Fs \,, \]
where the central arrow is induced by the canonical inclusion $(\Ls^\taut)^\vee \hra \pEs$. Put
\[ \Fs^i:=\pEs^i \otimes \Ls^\taut, \;\; \Gs^i:=\Fs^i/\Fs^{i+1} \;\; \hbox{ and } \;\; \Vs^i:=s^{-1}(\Fs^i)  \;\; \hbox{ for } \;\; 0 \le i \le r. \]
Note that $\Fs^r=\Gs^r=0$ (as vector bundles over $\bP(\Es)$) and that $\Fs^{r-1}=\Gs^{r-1}$ is a line bundle over $\bP(\Es)$.
When $X_0$ is a point (and $E_0 \cong \bA^r_{X_0}$), the degree $0$ part of $s$ sends
\[ \begin{matrix} x=(b_1:b_2:\dotsb:b_r) & \longmapsto & (b_1,b_2,\dotsc, b_r) \otimes \bs{v} \\
      \qquad  \bigcap\hspace{-9pt}\bs{|}  & \phantom{\Big|} & \bigcap\hspace{-9pt}\bs{|} \;\;\; \\
  \bP(E_0) = \bP(\bA^r_{X_0}) \quad\; & \lra & F_0, \; \end{matrix}\]
where $\bs{v}$ denotes the dual vector of $(b_1,b_2,\dotsc,b_r) \in (L^\taut_0)^\vee_x$. By this local description of $s$, we see the following (where $\Xs$ is arbitrary):
\begin{itemize}
\item {\it $V^i_j$ is smooth over $X_j$ for any $i \le r-1$ and $j \ge 0$, and $V_j^r$ is empty for any $j \ge 0$. In particular, $\Vs^i$ is a simplicial object in $\cC$ for each $i=0,1,\dotsc, r$.}
\item {\it For each $i=0,1,\dotsc,r-1$, $\Vs^{i+1}$ is an effective Cartier divisor on $\Vs^i$ {\rm(}cf.\ {\rm\cite{EGAIV}}, Th\'eor\`eme 17.12.1{\rm)}, whose associated line bundle is isomorphic to $\Gs^i|_{\Vs^i}$\,. }
\end{itemize}
Now the assertion follows from these facts and a similar computations as in \cite{Gr} Proof of Lemma 2.
This completes the proof of Lemma \ref{lem4-1} and Theorem \ref{thm4-1}.
\end{pf}

For a simplicial scheme $\Xs$, let $\tK_0(\Xs)$ be the Grothendieck group of vector bundles over $\Xs$. As an immediate consequence of Theorem \ref{thm4-1}, we obtain the following corollary:
\begin{cor}\label{cor4-1}
For a simplicial object $\Xs$ in $\cC$, there exists a unique map
\[ \sfc=(\sfc_i)_{i \ge 0} : \tK_0(\Xs) \lra \tH^{2*}(\Xs,\sfGa(*)), \]
that satisfies the following four properties:
\begin{enumerate}
\item[{\rm(1)}]
We have $\sfc_0 \equiv 1$ {\rm(}constant{\rm)}. If $\alpha \in \tK_0(\Xs)$ is the class of a line bundle $\Ls$, then $\sfc_1(\alpha)$ agrees with $\sfc_1(\Ls)$.
\item[{\rm(2)}]
The map $\sfc$ is contravariantly functorial for morphisms of simplicial objects in $\cC$.
\item[{\rm(3)}]
For $\alpha,\beta \in \tK_0(\Xs)$, we have $\sfc(\alpha+\beta)=\sfc(\alpha) \cup \sfc(\beta)$.
\item[{\rm(4)}]
For  $\alpha,\beta \in \tK_0(\Xs)$, we have $\wt{\sfc}(\alpha \otimes \beta) = \wt{\sfc}(\alpha) \,\scstar\; \wt{\sfc}(\beta)$ in $\wtH^{2*}(\Xs,\sfGa(*))$. Here $\wt{\sfc}(\alpha)$ denotes the augmented total Chern class $(\rk(\alpha),\sfc(\alpha))$, and $\rk$ denotes the virtual rank function.
\end{enumerate}
\end{cor}

\section{Universal Chern class and character}\label{sect5}
Let the notation be as in \S\ref{sect4} and recall our convention on schemes we fixed in \S\ref{notation}.
We introduce here the following notation:
\begin{defn}\label{defn5-1}
{\rm Let $Y$ be a scheme, and let $\sfS(Y_\Zar)$ (resp.\ $\sfS_\bullet(Y_\Zar)$) be the category of sheaves of sets (resp.\ sheaves of pointed sets) on $Y_\Zar$.
\begin{enumerate}
\item[(1)]
We endow $\SS(Y_\Zar)$, the category of simplicial sheaves of pointed sets on $Y_\Zar$\,, with the Brown-Gersten model structure \cite{BG} Theorem 2, whose class of fibrations (resp.\ weak equivalences, cofibrations) are defined as that of global fibrations (resp.\ morphisms which induce topological weak equivalences on stalks, morphisms which have left lifting property with respect to all trivial fibrations). We write $\cHs(Y)$ for its associated homotopy category.
We will also use the unpointed version $\cHsu(Y)$ constructed from $\SSu(Y_\Zar)$, the category of simplicial sheaves of sets.
\item[(2)]
For a cochain complex $(\cF^\bullet,d^\bullet)$ of abelian sheaves on $Y_\Zar$ and an integer $j \in \bZ$, consider the following complex:
\[\dotsb \lra \cF^{j-2} \os{d^{j-2}}\lra \cF^{j-1} \os{d^{j-1}}\lra \ker(d^j: \cF^j \to \cF^{j+1}), \]
which we regard as a chain complex with the most right term placed in degree $0$. Taking the associated simplicial abelian sheaf to this complex (cf.\ \cite{GJ} p.\ 162), we obtain a simplicial abelian sheaf on $Y_\Zar$, which we denote by $\cK\!(\cF^\bullet,j)$.
\item[(3)]
Let $\BQP_{\!Y}$ be the nerve of the ${\text {\it Q}}$-category associated with the exact category ${\text {\it P}}_Y$ of locally free $\cO_Y$-modules. We call $\BQP_{\!Y}$ {\it the {\it K}-theory space} of $Y$.
Let $\BQPs_{\!Y}$ be the simplicial sheaf of pointed sets on $Y_\Zar$ associated with the presheaf
\[ U \subset Y \;\hbox{(open)} \longmapsto \BQP_U. \]
The $j$-th algebraic \tK -group $\tK_j(Y)$ of $Y$ is $\pi_{j+1}(\BQP_{\!Y})$ by definition \cite{Q}.
there is a natural map from $\tK_j(Y)$ to the generalized sheaf cohomology $H^{-j-1}(Y,\BQPs_{\!Y})$
 (see \S\ref{sect6} below), which is bijective if $Y$ is regular (cf.\ \cite{Q} \S7.1, \cite{BG} Theorem 5).
\end{enumerate}}
\end{defn}

Let $X$ be a scheme which belongs to $\Ob(\cC)$.
 We review the construction of universal Chern classes due to Gillet \cite{Gi} \S2 briefly,
 which will be complete in 3 steps.
\par\medskip\noindent
{\bf Step 1.}\;
Let $n$ be a non-negative integer, and let $\BGL_{n,X}$ be the classifying scheme of $\GL_{n,X}$, the general linear group scheme of degree $n$ over $X$. Applying the construction of Chern classes to the universal rank $n$ bundle $E^\univ_\star$ over $\BGL_{n,X}$, we obtain Chern classes
\[ \sfc_i(E^\univ_\star) \in \tH^{2i}(\BGL_{n,X},\sfGa(i)) \qquad \hbox{($i \ge 0$)} \]
which are called {\it the universal rank $n$ Chern classes}. Recall that there is a canonical map
\[ \alpha_{X,i,n} : \tH^{2i}(\BGL_{n,X},\sfGa(i)) \lra \Mor_{\cHsu(X)}(\BGLs_n(\cO_X),\cK\!(\sfGa(i)_X,2i)) \]
for $n \ge 1$ and $i \ge 0$, cf.\ \cite{Gi}, p.\ 221.
Here $\BGLs_n(\cO_X)$ denotes the simplicial sheaf of groups on $X$ represented by $\BGL_{n,X}$.
For $i \ge 0$, we define
\begin{equation}\notag 
 \sfc_{i,n}:=\alpha_{X,i,n}(\sfc_i(E^\univ_\star)) \in \Mor_{\cHsu(X)}(\BGLs_n(\cO_X),\cK\!(\sfGa(i)_X,2i))
\end{equation}
One can easily check that $\sfc_{i,n}$ for $i \ge 1$ is pointed, i.e., defines a morphism in $\cHs(X)$.
More precisely, we have used here the following well-known fact, cf.\ \cite{AGP} Corollary 4.4.7:
\begin{prop}\label{prop:pointed}
Let $\cF^\bullet$ be a cochain complex of abelian sheaves on $X_\Zar$. Then the canonical map forgetting base points
\[ \Mor_{\cHs(X)}(\cE_\star,\cK\!(\cF^\bullet,j)) \lra \Mor_{\cHsu(X)}(\cE_\star,\cK\!(\cF^\bullet,j)) \]
is injective for any $\cE_\star \in \Ob(\SS(X_\Zar))$ and $j \in \bZ$.
\end{prop}

\par\medskip\noindent
{\bf Step 2.}\;
The data $(\sfc_{i,n})_{n \ge 1}$ together with the stability (\cite{Gi} Theorem 1.12, Proposition 1.17)
yield a mapping class
\begin{equation}\notag 
 \sfc_{i,\infty} \in \Mor_{\cHs(X)}(\BGLs(\cO_X),\cK\!(\sfGa(i)_X,2i)) \qquad \hbox{($i \ge 1$)},
\end{equation}
loc.\ cit.\ p.\ 225. To proceed the construction of universal Chern classes,
 we recall here the following fact (compare with \cite{Sch} Warning 2.2.9):
\begin{prop}[{\bf cf.\ \cite{Gi} Proposition 2.15}]
There is a functorial isomorphism
\[ \OBQPs_{\!X} \cong \bZ \times \bZ_\infty \BGLs(\cO_X) \quad \hbox{ in } \;\; \cHs(X). \]
Here $\OBQPs_{\!X}$ denotes the loop space of $\BQPs_{\!X}$ {\rm(}cf.\ Definition \ref{defn5-1}\,{\rm(}3{\rm)}, {\rm\cite{GJ}} p.\ 33{\rm)},
which is a simplicial sheaf of pointed sets on $X_\Zar$ and whose base point is the constant loop at $0$.
\end{prop}
\par\medskip\noindent
{\bf Step 3.}\;
Finally we define {\it the $i$-th universal Chern class} $\sfC_i$ for $i \ge 1$ as the composite
\begin{equation}\label{eq9-1}
\begin{CD}
\sfC_i : \OBQPs_{\!X} @>{\simeq}>> \bZ \times \bZ_\infty \BGLs(\cO_X) \phantom{\cK\!(\sfGa(i)_X,2i) \cong} \\
  @>{\pr_2}>> \bZ_\infty \BGLs(\cO_X) \phantom{\bZ_\infty\cK\!(\sfGa(i)_X,2i) \cong}
 \\ @>{\bZ_\infty(\sfc_{i,\infty})}>> \bZ_\infty\cK\!(\sfGa(i)_X,2i) \cong  \cK\!(\sfGa(i)_X,2i)
\end{CD}
\end{equation}
in $\cHs(X)$, where the last isomorphism is obtained from the fact that simplicial abelian groups are $\bZ$-complete \cite{BK1} 4.2. For $i=0$, we define $\sfC_0 : \OBQPs_{\!X} \to \cK\!(\sfGa(0)_X,0)$ as the class of the constant map with value $1$ in $\cHsu(X)$ (not in $\cHs(X)$).
\par\medskip
We next review the universal augmented total Chern class and the universal Chern character, which will be useful later.
\begin{defn}[{\rm \cite{Gi}\,Definition 2.27, 2.34}]\label{def5-1}
{\rm Let $X$ be a scheme which belongs to $\Ob(\cC)$.
\begin{enumerate}
\item[(1)]
For $\cEs \in \Ob(\SS(X_\Zar))$, we define $\tH^{2i}(X,\cEs;\sfGa(i))$ (resp.\ $\tH^0(X,\cEs;\bZ)$) as the mapping class group $\Mor_{\cHs(X)}(\cEs,\cK\!(\sfGa(i)_X,2i))$ (resp.\ $\Mor_{\cHs(X)}(\cEs,\bZ)$).
\item[(2)]
We define {\it the universal augmented total Chern class} as
\[ \!\!\wt{\sfC}:=(\rk,\sfC_0,\sfC_1,\sfC_2,\dotsc) \in
\tH^0(X,\OBQPs_{\!X};\bZ) \times \{1\} \times \prod_{i \ge 1} \ \tH^{2i}(X,\OBQPs_{\!X};\sfGa(i)), \]
where $\rk$ denotes the class of the rank function $\OBQPs_{\!X} \to \bZ$.
\item[(3)]
(cf.\ \cite{Gr2} Chapter I (1.29))\; We define {\it the universal Chern character} as
\[ \ch := \rk + \eta \bigg(\log \bigg( 1+\sum_{i = 1}^\infty \ \sfC_i \bigg)\bigg) \in \prod_{i \ge 0} \ (\tH^{2i}(X,\OBQPs_{\!X};\sfGa(i)) \otimes \bQ), \]
where $\eta=(\eta_i)_{i \ge 0}$ denotes the graded additive endomorphism
\[ \eta_i(x_i) := \frac{(-1)^{i-1}}{(i-1)!} \cdot x_i \quad \hbox{($x_i \in \tH^{2i}(X,\OBQPs_{\!X};\sfGa(i)) \otimes \bQ$)}. \]
\end{enumerate}}
\end{defn}
For a morphism $f : Y \to X$ in $\cC$ and $i \ge 1$, there is a commutative diagram in $\cHs(Y)$
\begin{equation}\label{eq5-4}
\xymatrix{f^*\OBQPs_{\!X} \ar[d]_{\sfC_{i,X}} \ar[r]^-{f^\back} &  \ar[d]^{\sfC_{i,Y}} \OBQPs_{\!Y} \\
f^*\cK\!(\sfGa(i)_X,2i) \ar[r]^-{f^\back} & \cK\!(\sfGa(i)_Y,2i)}
\end{equation}
by Theorem \ref{thm4-1}\,(2) and the construction of $\sfC_i$. One can also check a similar commutativity for the universal Chern character.

\section{Chern class and character for higher $\bs{K}$-theory}\label{sect6}
\setcounter{equation}{0}
Let the notation be as in \S\ref{sect4}. Let $X$ be a scheme which belongs to $\Ob(\cC)$, and let $Z$ be a closed subset of $X$. For a simplicial sheaf of pointed sets $\cEs \in \Ob(\SS(X_\Zar))$ and a non-negative integer $j \ge 0$, we define
\[ \tH_Z^{-j}(X,\cEs):=\Mor_{\cHs(X)}(\cS^j_Z,\cEs), \]
where $\cS^j_Z$ denotes the constant sheaf (of simplicial pointed sets) on $Z$ associated with the singular simplicial set of the $j$-sphere. Note that there is a canonical isomorphism
\begin{equation}\label{lem6-1}
 \tH_Z^{-j}(X,\cK\!(\cF^\bullet,i)) \cong \tH^{i-j}_Z(X,\cF^\bullet)
\end{equation}
for a cochain complex $\cF^\bullet$ of abelian sheaves on $X_\Zar$ and $i \in \bZ$,
where the right hand side means the hypercohomology of $\cF^\bullet$ with support in $Z$, cf.\ \cite{BG} \S2 Proposition 2.
Let $\tK_j^Z(X)=\tK_j(X,X\ssm Z)$ be the $j$-th algebraic {\it K}-group of $X$ with support in $Z$. We define the Chern class map
\begin{equation}\notag 
 \sfC_{i,j,X}^Z : \tK_j^Z(X) \lra \tH^{2i-j}_Z(X,\sfGa(i))
\end{equation}
as the composite map
\[ \tK^Z_j(X) \lra \tH_Z^{-j}(X,\OBQPs_{\!X}) \os{\sfC_i}\lra \tH_Z^{-j}(X,\cK\!(\sfGa(i)_X,2i))
 \os{\eqref{lem6-1}} \cong \tH^{2i-j}_Z(X,\sfGa(i)). \]
Here $\sfC_i$ denotes the universal Chern class \eqref{eq9-1}. The map $\sfC_{i,0,X}^X$, i.e., the case $j=0$ and $Z=X$ agrees with $\sfc_i$ for $\Xs=X$ (constant simplicial scheme) mentioned in Corollary \ref{cor4-1}.
\begin{prop}\label{prop6-1}
\begin{enumerate}
\item[{\rm(1)}]
$\sfC_{i,j}^Z$ is contravariantly functorial in the pair $(X,Z)$, that is, for a morphism $f :X' \to X$ in $\cC$ and a closed subset $Z' \subset X'$ with $f^{-1}(Z) \subset Z'$, there is a commutative square
\[ \xymatrix{\tK_j^Z(X) \ar[r]^-{f^\back} \ar[d]_{\sfC_{i,j,X}^Z} & \tK_j^{Z'}(X') \ar[d]^{\sfC_{i,j,X'}^{Z'}} \\
 \tH^{2i-j}_Z(X,\sfGa(i)) \ar[r]^-{f^\back} & \tH^{2i-j}_{Z'}(X',\sfGa(i)). } \]
\item[{\rm(2)}]
$\sfC_{i,j,X}^Z$ is additive for $i,j >0$, and zero for $i=0$ and $j>0$.
\end{enumerate}
\end{prop}
\begin{pf}
(1) follows from the commutative diagram \eqref{eq5-4}. The additivity assertion of (2) follows from Theorem \ref{thm4-1}\,(3) and the arguments in \cite{Gi} Lemma 2.26. The last assertion follows from the fact that $\tH^{-j}_Z(X,\sfGa(0))$ is zero for $j > 0$, cf.\ Definition \ref{axiom0}\,(a).
\end{pf}
\nega
\begin{rem}
To prove Proposition \ref{prop6-1}\,{\rm(}2{\rm)}, we need the framework of Chern classes of representations {\rm\cite{Gi}} Definitions 2.1 and 2.10, which we omit in this paper because one can easily establish it under our setting by the same arguments as in loc.\ cit. 
\end{rem}
\begin{defn}\label{def6-3}
{\rm\begin{enumerate}
\item[(1)]
We define $\tK^Z_*(X)$ as the direct sum of $\tK^Z_i(X)$ with $i \ge 0$, and define $\whH^*_Z(X,\sfGa(\scbullet))_\bQ$ as the direct product of $\tH^i_Z(X,\sfGa(n))\otimes \bQ$ with $i,n \ge 0$.
\item[(2)]
We define the Chern character
\[ \ch_X^Z : \tK^Z_*(X) \lra \whH^*_Z(X,\sfGa(\scbullet))_\bQ \]
as the map induced by the universal Chern character defined in Definition \ref{def5-1}\,(3).
We often write $\ch_X$ for $\ch_X^X$.
\end{enumerate}}
\end{defn}

\begin{prop}\label{prop6-2}
\begin{enumerate}
\item[{\rm(1)}]
$\ch_X^Z$ is contravariantly functorial in the pair $(X,Z)$ in an analogous sense to Proposition \ref{prop6-1}\,(1).
\item[{\rm(2)}]
The Chern character $\ch_X^Z$ is a ring homomorphism.
\end{enumerate}
\end{prop}
\begin{pf}
(1) the commutative diagram \eqref{eq5-4}. 
The assertion (2) follows from Theorem \ref{thm4-1}\,(4), Proposition \ref{prop6-1}\,(2) and the arguments in \cite{Gi} Proposition 2.35.
\end{pf}

\section{Push-forward for projective morphisms}\label{sect7}
Let $\cC$ and $\sfGa(*)$ be as in \S\ref{sect4}. See \S\ref{notation} for the definition of projective morphisms.
\begin{defn}\label{def7-1}
{\rm Let $f : Y \to X$ be a projective morphism in $\cC$.
\begin{enumerate}
\item[(1)]
By taking a factorization
\[\xymatrix{ f : Y \;\ar@<-1pt>@{^{(}->}[r]^-g & \bP^m_X \ar[r] & X }\]
with $g$ a closed immersion, we define {\it the relative dimension} of $f$ as the integer $m-\codim(g)$. Because we deal with only universally catenary schemes, this number is independent of the factorization.
\item[(2)]
We say that $f$ is {\it regular} if $f$ has a factorization as above for which $g$ is a regular closed immersion. A regular projective morphism is a regular morphism in the sense of \cite{FL} p.\ 86.
\end{enumerate}}
\end{defn}
\par
Now let $\sfGa(*)$ be an admissible cohomology theory on $\cC$, and let $f : Y \to X$ be a regular projective morphism in $\cC$. The main aim of this section is to construct push-forward morphisms in $\sfD(X_\Zar)$
\begin{equation}\label{eq7-1}
 f_\push : \tR f_*\sfGa(i+r)_Y[2r] \lra \sfGa(i)_X \qquad \hbox{($i \in \bZ$)}
\end{equation}
and prove Theorem \ref{thm7-1} and Corollary \ref{cor7-2} below, where $r$ denotes the relative dimension of $f$. The results in this section will play key roles in the following sections. Taking a factorization
\[\xymatrix{ f : Y \;\ar@<-1pt>@{^{(}->}[r]^-g & \bP^m_X=:\bP^m \ar[r]^-p & X }\]
of $f$ such that $g$ is a regular closed immersion, we are going to define the push-forward morphism \eqref{eq7-1} as the composite
\[\xymatrix{ \tR f_*\sfGa(i+r)_Y[2r] \cong \tR p_*\tR g_*\sfGa(i+r)_Y[2r] \ar[r]^-{\tR p_*(g_\push)} &  \tR p_*\sfGa(i+m)_{\bP^m}[2m] \ar[r]^-{p_{\sharp}} & \sfGa(i)_X\,, }\]
where $p_\sharp$ denotes the composite of the isomorphism by the projective bundle formula and a projection
\[ p_\sharp : \tR p_*\sfGa(i+m)_{\bP^m}[2m] \cong \bigoplus_{j=0}^m \ \sfGa(i+j)_X[2j] \lra \sfGa(i)_X\,. \]
We have used the condition ($*_2$) on the category $\cC$ in \S\ref{sect1.1} to verify the existence of $g_\push$. 
\begin{thm}\label{thm7-1}
Let $f : Y \to X$ be a regular projective morphism in $\cC$.
\begin{enumerate}
\item[{\rm(1)}]
{\rm(}Well-definedness{\rm)} \;
$f_\push$ does not depend on the choice of a factorization of $f$. In particular, we have $f_\push=f_\sharp$ when $f$ is isomorphic to a natural projection $\bP^m_X \to X$.
\item[{\rm(2)}]
{\rm(}Projection formula{\rm)} \;
The following diagram commutes in $\sfD(X_\Zar)${\rm:}
\[\hspace{-10pt} \xymatrix{ \sfGa(i)_X \os{\bL}\otimes \tR f_*\sfGa(j+r)_Y[2r] \ar[r]^-{\id \otimes f_\push} \ar[d]_{f^\back \otimes \id} & \sfGa(i)_X \os{\bL}\otimes \sfGa(j)_X \ar[rd]^-{\text{product}} \\
 \tR f_*\sfGa(i)_Y \os{\bL}\otimes \tR f_*\sfGa(j+r)_Y[2r] \ar[r]^-{\text{product}} & \tR f_*\sfGa(i+j+r)_Y [2r] \ar[r]^-{f_\push} &  \sfGa(i+j)_X\, .} \]
\item[{\rm(3)}]
{\rm(}Transitivity{\rm)} \;
For another regular projective morphism $f' : Z \hra Y$ in $\cC$ of relative dimension $r'$, the composite morphism
\begin{align*} \tR (f \circ f')_*\sfGa(i+r+r')_Z[2(r+r')] & \cong \tR f_* \tR f'_*\sfGa(i+r+r')_Z [2(r+r')] \\
& \os{f'_\push}\lra \tR f_*\sfGa(i+r)_Y[2r]  \os{f_\push}\lra \sfGa(i)_X \end{align*}
agrees with $(f \circ f')_\push$\,.
\item[{\rm(4)}]
{\rm(}Base-change property{\rm)} \;
Let 
\[\xymatrix{ Y' \ar[r]^{f'} \ar[d]_\beta & X' \ar[d]^\alpha \\ Y \ar[r]^f & X }\]
be a commutative diagram in $\cC$, where $f$ and $f'$ are regular projective morphisms of relative dimension $r$, and $\alpha$ is a closed immersion. Let $U$ be an open subset of $Y$ for which the following square is cartesian{\rm:}
\[\xymatrix{ \beta^{-1}(U) \ar[rr]^-{f'|_{\beta^{-1}(U)}} \ar[d] \ar@{}[rrd]|{\square} && X' \ar[d]^\alpha \\ U \ar[rr]^-{f|_U} && X.} \]
Then for a closed subset $Z \subset Y$ contained in $U$, the diagram
\def\uG{\ul{\vG}{\,}}
\[\xymatrix{ \tR (\alpha \circ f')_*\tR\uG_{\beta^{-1}(Z)}(Y',\sfGa(i+r))[2r] \ar[r]^-{f'_\push} & \alpha_*\tR\uG_{\alpha^{-1}(f(Z))}(X',\sfGa(i))
 \\ \tR f_*\tR\uG_Z(Y,\sfGa(i+r))[2r] \ar[r]^-{f_\push} \ar[u]^{\beta^\back} & \tR\uG_{f(Z)}(X,\sfGa(i)) \ar[u]_{\alpha^\back} } \] commutes in $\sfD(X_\Zar)$.
\end{enumerate}
\end{thm}
\begin{rem}
\begin{enumerate}
\item[{\rm(1)}]
Applying Theorem \ref{thm7-1} to the example in \S\ref{ex:loghw}, we obtain push-forward morphisms of logarithmic Hodge-Witt sheaves
for projective morphisms of regular schemes over $\bF_p$, which satisfy the properties listed above. This result answers the problem raised in {\rm\cite{Sh}} Remark 5.5 affirmatively, and we obtain the same compatibility as in loc.\ cit.\ Theorem 5.4 for $\logwittt n q$ with $n \ge 2$ as well.
\item[{\rm(2)}]
We need the axiom \ref{axiom1}\,{\rm(}3{\rm)} only for regular closed immersions of usual schemes in $\cC$ to prove Theorem \ref{thm7-1}.
\end{enumerate}
\end{rem}

We first prepare the following lemmas:
\begin{lem}\label{lem7-0}
Let $R$ be a local ring and put $Y:=\Spec(R)$. Let $s : Y \to \bP^m_Y$ be a section of the natural projection $\bP^m_Y \to Y$. Then under a suitable choice of coordinates, there exists an affine open subset $\bA^m_Y \subset \bP^m_Y$ such that $s(Y) \subset \bA^m_Y$.
\end{lem}
\begin{pf}
Let $T_0,T_1,\dotsc,T_m$ be a set of homogeneous coordinates of $\bP^m$.
The global sections $a_i:=s^*(T_i) \in \vG(Y,s^*\cO(1))$ ($i=0,1,\dotsc,m$) generates the sheaf $s^*\cO(1)$, cf.\ \cite{Ha}\,II Theorem 7.1\,(a).
Since $R$ is local, we have $s^*\cO(1) \cong \cO_Y$ and $\vG(Y,s^*\cO(1)) \cong R$.
Under this identification, the ideal generated by $a_i$'s is $R$ itself and we have $a_{i_0} \in R^\times$ for some $i_0$, again by the assumption that $R$ is local.
Hence $s(Y)$ is contained in the complement of the hyperplane $\{T_{i_0} = 0\}$, which is the desired affine open subset.
\end{pf}

\begin{lem}\label{lem7-1}
Let $X$ be an object of $\cC$. Let $m$ and $n$ be non-negative integers, and put $\bP^m := \bP^m_X$, $\bP^n := \bP^n_X$ and $\bP^m \times \bP^n:=\bP^m \times_X \bP^n$.
\begin{enumerate}
\item[{\rm(1)}]
Let $p : \bP^m \to X$ be the natural projection, and let $s : X \to \bP^m$ be a section of $p$. Then the composite morphism in $\sfD(X_\Zar)$ \[ \sfGa(i)_X \os{s_\push}\lra Rp_*\sfGa(i+m)_{\bP^m}[2m] \os{p_\sharp}\lra \sfGa(i)_X \] agrees with the identity morphism of $\sfGa(i)_X$ for each $i \in \bZ$. Here $s_\push$ denotes the push-forward morphism in Definition \ref{axiom1}\,{\rm(}3{\rm)}.
\item[{\rm(2)}]
Put $N:=mn+m+n$. Let $p' : \bP^n \to X$ and $q : \bP^N:=\bP^N_X \to X$ be the natural projections. Let $\psi : \bP^m \times \bP^n \hra \bP^N$ be the Segre embedding, and let $\pi : \bP^m \times \bP^n \to \bP^m$ be the first projection. Then the diagram
\[\xymatrix{ \tR (p \times p')_*\sfGa(i+m+n)_{\bP^m \times \bP^n}[2(m+n)] \ar[r]^-{\psi_\push} \ar[d]_{\pi_\sharp} & \tR q_*\sfGa(i+N)_{\bP^N}[2N] \ar[d]^{q_\sharp} \\ \tR p_*\sfGa(i+m)_{\bP^m}[2m] \ar[r]^-{p_\sharp} & \sfGa(i)_X}\]
is commutative in $\sfD(X_\Zar)$ for each $i \in \bZ$.
\item[{\rm(3)}]
Let $g : Y \hra X$ be a regular closed immersion of codimension $c$ in $\cC$. Let $p_Y : \bP^m_Y \to Y$ be the natural projection, and let $g' : \bP^m_Y \hra \bP^m(=\bP^m_X)$ be the regular closed immersion induced by $g$. Then the diagram
\[ \xymatrix{ g_*\tR p_{Y*}\sfGa(i+m)_{\bP^m_Y}[2m] \ar[rr]^-{\tR p_*(g'_\push)} \ar[d]_{p_{Y\sharp}} && \tR p_*\sfGa(i+m+c)_{\bP^m}[2(m+c)] \ar[d]^{p_\sharp} \\ g_*\sfGa(i)_Y \ar[rr]^-{g_\push} && \sfGa(i+c)_X[2c]}\]
commutes in $\sfD(X_\Zar)$ for each $i \in \bZ$.
\item[{\rm(4)}]
Let $Y \in \Ob(\cC)$ be a scheme which is projective over $X$, and let $g : Y \hra \bP^m$ and $h : Y \hra \bP^n$ be regular closed immersions over $X$. Let $\varphi : Y \hra \bP^m \times \bP^n$ be the closed immersion induced by $g$ and $h$. Then $\varphi$ is a regular closed immersion, and the diagram
\[ \xymatrix{ g_*\sfGa(i)_Y \ar[rr]^-{\tR\pi_*(\varphi_\push) } \ar[rd]_{g_\push} && \tR\pi_*\sfGa(i+n+c)_{\bP^m \times \bP^n}[2(n+c)] \ar[ld]^{\pi_\sharp} \\ & \sfGa(i+c)_{\bP^m}[2c] } \]
commutes in $\sfD((\bP^m)_\Zar)$ for each $i \in \bZ$. Here $c$ denotes $\codim(g)$, and $\pi : \bP^m \times \bP^n \to \bP^m$ denotes the first projection.
\end{enumerate}
\end{lem}
\begin{pf*}{Proof of Lemma \ref{lem7-1}}
(1) We first show that
\begin{equation}\label{eq7-1-1}
 p_\sharp(s_\push(1))=1 \quad \hbox{ in } \;\; \tH^0(X,\sfGa(0)),
\end{equation}
where $1$ means the unity of $\tH^0(X,\sfGa(0))$. Replacing $X$ with an affine dense open subset if necessary, we may assume that $s : X \hra \bP^m$ is given by the homogeneous coordinate $(1\,{:}\,a_1\,{:}\,\dotsc\,{:}\,a_n)$ with $a_i:=s^\sharp(T_i/T_0) \in \vG(X,\cO_X)$ under a suitable choice of projective coordinates $(T_0\,{:}\,T_1\,{:}\,\dotsc\,{:}\,T_n)$ with $s(X) \subset \{T_0 \ne 0\}=\bA^n$, cf.\ Definition \ref{axiom1}\,(0), Lemma \ref{lem7-0}. Let $\xi \in \tH^2(\bP^m,\sfGa(1))$ be the first Chern class of the tautological line bundle over $\bP^m$. In order to show \eqref{eq7-1-1}, it is enough to check
\begin{equation}\label{eq7-1-2}
 s_\push(1)=\xi^m \quad \hbox{ in } \;\; \tH^{2m}(\bP^m,\sfGa(m)).
\end{equation}
If $m=1$, we have $s_\push(1)=\xi$ by Definition \ref{axiom1}\,(3a). For $m \geq 2$, we have $s_\push(1)=\xi^m$ by Definition \ref{axiom1}\,(3b), (3c) and induction on $m$. Thus we obtain \eqref{eq7-1-1}.

Lemma \ref{lem7-1}\,(1) follows from \eqref{eq7-1-1} and the following commutative diagram in $\sfD(X_\Zar)$:
\[\xymatrix{ \sfGa(i)_X \ar[r]^-{\id \otimes s_\push(1)} \ar[d]_{s_\push}
 & \sfGa(i)_X \os{\bL}\otimes Rp_*\sfGa(m)_{\bP^m}[2m] \ar[d]^{\id \otimes p_\sharp} \ar[ld]_{\text{($\star$)}} \\
Rp_*\sfGa(i+m)_{\bP^m}[2m] \ar[r]^-{p_\sharp} & \sfGa(i)_X,
}\]
where the arrow ($\star$) is $p^\back \otimes \id$ followed by product.
The left triangle commutes by the projection formula in Definition \ref{axiom1}\,(3b). On the other hand, the right triangle commutes as well because the projective bundle formula in Definition \ref{axiom1}\,(2) is compatible with the multiplication by $\sfGa(*)_X$. This completes the proof of Lemma \ref{lem7-1}\,(1).
\par
(2) Let $s : X \to \bP^m$ and $s': X \to \bP^n$ be the zero sections.
Let $\xi \in \tH^2(\bP^m,\sfGa(1))$ (resp.\ $\eta \in \tH^2(\bP^n,\sfGa(1))$, $\zeta \in \tH^2(\bP^N,\sfGa(1))$) be the first Chern class of the tautological line bundle over $\bP^m$ (resp.\ $\bP^n$, $\bP^N$). Let $\sigma : X \to \bP^m \times \bP^n$ be the morphism induced by $s$ and $s'$. By similar arguments as in the proof of \eqref{eq7-1-2}, we see that
{\allowdisplaybreaks
\begin{align*}
 &\sigma_\push(1) = \pi^\back (\xi^m) \cup \pi'^\back (\eta^n) \quad \hbox{ in } \;\; \tH^{2(m+n)}(\bP^m \times \bP^n,\sfGa(m+n)), \\
 &(\psi \circ \sigma)_\push(1) = \zeta^N \qquad \qquad\, \hbox{ in } \;\; \tH^{2N}(\bP^N,\sfGa(N)),
\end{align*}
}where $\pi'$ denotes the second projection $\bP^m \times \bP^n \to \bP^n$. By these facts and Definition \ref{axiom1}\,(3c), we obtain
\begin{equation}\label{eq7-1-3}
 \psi_\push(\pi^\back (\xi^m) \cup \pi'^\back (\eta^n))=\zeta^N \quad \hbox{ in } \;\; \tH^{2N}(\bP^N,\sfGa(N)). \end{equation}
On the other hand, if $0 \leq i \leq m$ and $0 \leq j \leq n$ with $i+j>0$, then we have
\begin{equation}\label{eq7-1-4}
 \psi_\push (\pi^\back (\xi^{m-i}) \cup \pi'^\back (\eta^{n-j}))
  = \zeta^{N-i-j} \cup q^\back(a_0) + \zeta^{N-i-j-1} \cup q^\back(a_1) + \dotsb+
   q^\back(a_{N-i-j})
\end{equation}
in $\tH^{2(N-i-j)}(\bP^N,\sfGa(N-i-j))$, for some $a_k = a_{ij,k} \in H^{2k}(X,\sfGa(k))$ ($k=0,1,\dotsc,N-i-j$)
 by the projective bundle formula in Definition \ref{axiom1}\,(2).
One can easily deduce the assertion of Lemma \ref{lem7-1}\,(2) from \eqref{eq7-1-3}, \eqref{eq7-1-4} and the projection formula in Definition \ref{axiom1}\,(3b).
\par
(3) By the axioms in Definition \ref{axiom1}\,(3), the following diagram commutes in $\sfD(X_\Zar)$:
\[ \xymatrix{ g_*\tR p_{Y*}\sfGa(i+m)_{\bP^m_Y}[2m] \ar[rr]^-{\tR p_*(g'_\push)} \ar@{=}[d]_{\wr} && \tR p_*\sfGa(i+m+c)_{\bP^m}[2(m+c)] \ar@{=}[d]^{\wr} \\ \bigoplus_{j=0}^m \ g_*\sfGa(i+j)_Y[2j] \ar[rr]^-{g_\push} && \bigoplus_{j=0}^m \ \sfGa(i+j+c)_X[2(j+c)],}\]
where the vertical isomorphisms follow from the projective bundle formula. The assertion follows from this fact.
\par
(4) The first assertion follows from \cite{SGA6} Expos\'e V\III\ Corollaire 1.3. As for the second assertion, replacing $\bP^m$ with $X$, we may assume that $m=0$ and that $Y$ is a closed subscheme of $X$ via $g$. Then decomposing $\varphi (=h)$ as $Y \hra \bP_Y^n \hra \bP^n_X$, we see that the assertion is reduced to the results in (1) and (3), by the transitivity in Definition \ref{axiom1}\,(3c).
\end{pf*}

\begin{pf*}{Proof of Theorem \ref{thm7-1}}
We write $\bP^m$ for the projective space $\bP^m_X$ over $X$ for simplicity. \par
(1) Suppose we are given two factorizations $Y \os{g}\hra \bP^m \os{p}\to X$ and $Y \os{h}\hra \bP^n \os{p'}\to X$ of $f$. There is a commutative diagram in $\cC$
\[ \xymatrix{Y \; \ar@<-1pt>@{^{(}->}[r]^-\varphi \ar[rd]_g & \bP^m \times \bP^n \; \ar@<-1pt>@{^{(}->}[r]^-\psi \ar[d] & \bP^{mn+m+n} \ar[d]^q \\ & \bP^m \ar[r]^p & X,}\]
where $\varphi$ denotes the closed immersion induced by $g$ and $h$, and $\psi$ denotes the Segre embedding. The right vertical arrow $q$ denotes the natural projection. Note that $\varphi$ and $\psi$ are both regular by
 \cite{SGA6} Expos\'e V\hspace{-1pt}\III\ Corollaire 1.3 and \cite{EGAIV} Th\'eor\`eme 17.12.1, respectively. By Lemma \ref{lem7-1}\,(2), (4) and the axioms in Definition \ref{axiom1}\,(3), we have
\[ p_\sharp \circ \tR p_*(g_\push) = q_\sharp \circ \tR q_*((\psi \circ \varphi)_\push) :  \tR f_*\sfGa(i+r)_Y[2r] \lra \sfGa(i)_X ,\]
which implies the assertion.
\par
(2) Fix a factorization $Y \os{g}\hra \bP^m \os{p}\to X$ of $f$. The projection formula holds for $g$ by the axiom in Definition \ref{axiom1}\,(3b), and holds for $p$ as well because the isomorphism of the projective bundle formula is compatible with the multiplication by $\sfGa(*)_X$, cf.\ Definition \ref{axiom1}\,(2). The assertion follows from these facts.
\par
(3) Taking factorizations $f : Y \os{g}\hra \bP^m \os{p}\to X$ and $f' : Z \os{h}\hra \bP^n \os{p'}\to Y$, one can easily deduce the assertion from Lemma \ref{lem7-1}\,(3), (4). The details are left to the reader.
\par
(4) Fix a factorization $Y \os{g}\hra \bP^m \os{p}\to X$ of $f$. There are cartesian squares
\[\xymatrix{ \beta^{-1}(U) \ar[rr] \ar[d] \ar@{}[rrd]|{\square} && \bP^m_{X'} \ar[d]  \ar[r]\ar@{}[rd]|{\square} & X' \ar[d]^\alpha \\ U \ar[rr]^-{g|_U} && \bP^m \ar[r] & X,} \]
where the horizontal arrows of the left square are (locally closed) immersions. The assertion holds for the right square and a closed subset $W \subset \bP^m$ by the definition of push-forward morphisms. On the other hand, the assertion holds for the left square and a closed subset $Z \subset Y$ contained in $U$ by excision and the base-change property in Definition \ref{axiom1}\,(3d). The assertion for $f$ follows from these facts and the transitivity established in (3).
\end{pf*}
\nega
\begin{cor}\label{cor7-1}
Let $Y \in \Ob(\cC)$ be a scheme which admits an ample family of invertible sheaves, and let $\pi : E \to Y$ be a vector bundle of rank $r+1$. Let $p : \bP(E) \to Y$ be the projective bundle associated with $E$, cf.\ \eqref{eq0-1}, which is projective in our sense by the assumption on $Y$. Then for each $i \in \bZ$, the map
$p_\push : \tR p_*\sfGa(i+r)_{\bP(E)}[2r] \to \sfGa(i)_Y$
agrees with the composite of the isomorphism by the projective bundle formula and a projection
\[ p_\sharp : \tR p_*\sfGa(i+r)_{\bP(E)}[2r] \cong \bigoplus_{j=0}^r \ \sfGa(i+j)_Y[2j] \lra \sfGa(i)_Y \quad \hbox{ in } \;\; \sfD(Y_\Zar). \]
\end{cor}
\nega
\begin{pf}
By a standard hyper-covering argument, we may assume that $E \cong \bA^{r+1}_Y$. Then the assertion follows from Theorem \ref{thm7-1}\,(1).
\end{pf}
\nega
\begin{cor}\label{cor7-2}
Let $Y \in \Ob(\cC)$ be a scheme which admits an ample family of invertible sheaves, and let $\pi : E \to Y$ be a vector bundle of rank $r$. Let $p : X:=\bP(E \oplus \bs{1}) \to Y$ be the projective completion of $E$, where $\bs{1}$ denotes the trivial line bundle over $Y$. Let $f : Y \hra X$ be the zero section of $p$, and let $Q$ be the universal quotient bundle $p^*(E \oplus \bs{1})/(L^\taut)^\vee$ on $X$. Then the Gysin map $f_\push$ sends the unity $1 \in \tH^0(Y,\sfGa(0))$ to $\sfc_r(Q) \in \tH^{2r}(X,\sfGa(r))$.
\end{cor}
\begin{pf}
Let $\xi:=\sfc_1(L^\taut) \in \tH^2(X,\sfGa(1))$ be the first Chern class of the tautological line bundle over $X$, and let $i_\infty : \bP(E) \hra \bP(E \oplus \bs{1})=X$ be the infinite hyperplane. The projective bundle formula in the axiom \ref{axiom1}\,(2) (for both  $X$ and $\bP(E)$) and the functoriality mentioned in Remark \ref{rem2-1} imply that the kernel of the pull-back map
\[ i_\infty^\back : \tH^{2r}(X,\sfGa(r)) \lra \tH^{2r}(\bP(E),\sfGa(r)) \]
is generated, over $\tH^0(Y,\sfGa(0))$, by the element
\begin{equation}\label{eq7-2}
 \xi^r + \xi^{r-1} \cup \psfc_1(E) + \xi^{r-2} \cup \psfc_2(E) + \dotsb  + \psfc_r(E) = \sfc_r(Q),
\end{equation}
where we have used Theorem \ref{thm4-1}\,(2), (3) for $\Xs=X$ to obtain the last equality. On the other hand, we have $i_\infty^\back \circ f_\push=0$, because $f_\push$ factors through $\tH^{2r}_Y(X,\sfGa(r))$ and $i_\infty^\back$ factors through $\tH^{2r}(X \ssm Y,\sfGa(r))$. Therefore $f_\push(1)$ belongs to $\ker(i_\infty^\back)$ and we have
\[ f_\push(1) = a \cdot \sfc_r(Q) \;\;\; \hbox{for some } \; a \in \tH^0(Y,\sfGa(0)). \]
It remains to check $a=1$. By the transitivity in Theorem \ref{thm7-1}\,(3), this claim is further reduced to showing $p_\push(\sfc_r(Q))=1$.
Finally this last equality follows from \eqref{eq7-2}, Theorem \ref{thm7-1}\,(2) and Corollary \ref{cor7-1} for $E \oplus \bs{1}$ over $Y$.
\end{pf}

\section{Construction of a universal polynomial}\label{sect8}
For an indeterminate $x$, we define
\[ \bZ[x]^\scr:=\{f(x) \in \bQ[x]\,|\, f(m) \in \bZ \, \hbox{ for any } \, m \in \bZ \}. \]
Let $n$ and $r$ be integers with $n \ge 0$ and $r \ge 1$. In this section, we construct a universal polynomial
\[ \sfp_{n,r}(t_0,t_1,\dotsc,t_n;u_1,u_2,\dotsc,u_r) \in \bZ[t_0]^\scr [t_1,\dotsc,t_n;u_1,u_2,\dotsc,u_r] \]
in an explicit way by modifying the polynomial of Fulton-Lang considered in \cite{FL} Chapter II \S4.
Our polynomial $\sfp_{n,r}$ agrees with the polynomial considered in \cite{Gr2} Chapter I Proposition 1.5 and \cite{J} \S1, up to signs of $u_j$'s with $j$ odd.
We will also provide Propositions \ref{prop8-1} and \ref{prop8-2} below concerning elementary properties of some power series related to this universal polynomial, which will be useful later in Theorem \ref{thm9-1} below.
\par
We start with indeterminates $\bs{a}=(a_1,a_2,\dotsc,a_n)$, $\bs{b}=(b_1,b_2,\dotsc,b_r)$ and a power series
\[ \sfF_{\!n,r}(\bs{a},\bs{b}):= \prod_{i=1}^n \; \prod_{j=0}^r \ \prod_{k_1<\dotsb<k_j} \ (1+a_i-b_{k_1}-b_{k_2}-\dotsb -b_{k_j})^{(-1)^j}, \]
which has constant term $1$ and is symmetric in $a_1,a_2,\dotsc,a_n$ and also in $b_1,b_2,\dotsc,b_r$. It is well-known that $\sfF_{\!n,r}(\bs{a},\bs{b})-1$ is divisible by $b_1 b_2 \dotsb b_r$ (cf.\ \cite{FL} p.\ 44).
We mention here a property of the power series $\sffa_r(\bs{b}):=\sfF_{1,r}(0,\bs{b})$. Let $s_i$ be the $i$-th elementary symmetric expression in $b_1,b_2,\dotsc,b_r$, and let $\sfG_r(t_1,t_2,\dotsc,t_r)$ be the power series satisfying \[ \sffa_r(\bs{b})=1+s_r \cdot \sfG_r(s_1,s_2,\dotsc,s_r). \]
\begin{prop}\label{prop8-1}
Let $\cC$ and $\sfGa(*)$ be as in \S\ref{sect4}, and let $\pi : E \to X$ be a vector bundle of rank $r$ with $X \in \Ob(\cC)$.
Put
\[ \lam_{-1}(E^\vee):= 1 - [E^\vee] + [{\sf\Lambda}^2 E^\vee] - \dotsb + (-1)^r [{\sf\Lambda}^rE^\vee] \in \tK_0(X). \]
Then in the complete cohomology ring $\whH^{2*}(X,\sfGa(*)):=\prod_{i \ge 0} \ \tH^{2i}(X,\sfGa(i))$, we have
\begin{equation}\label{eq8-1}
 \sum_{i \ge 0} \ \sfc_i(\lam_{-1}(E^\vee)) = 1+\sfc_r(E) \cup \sfG_r(\sfc_1(E),\sfc_2(E),\dotsc,\sfc_r(E)).
\end{equation}
\end{prop}
\begin{pf}
The assertion follows from Theorem \ref{thm4-1}\,(3) for $\Xs=X$ and \cite{FL} Chapter \II\, Proposition 4.1 for $e=\bs{1}$, the trivial line bundle over $X$.
\end{pf}
\par We next consider a power series
\begin{align*}
\sfJ_{n,r}(\bs{a},\bs{b}):= \sfF_{\!n,r}(\bs{a},\bs{b}) \cdot \sffa_r(\bs{b})^{-n},
\end{align*}
which is symmetric in $a_1,a_2,\dotsc,a_n$ and in $b_1,b_2,\dotsc,b_r$ as well.
\begin{lem}\label{lem8-1}
\begin{enumerate}
\item[{\rm(1)}] 
$\sfJ_{n,r}(\bs{a},\bs{b})-1$ is divisible by $b_1 b_2 \dotsb b_r$.
\item[{\rm(2)}] 
Let $\sigma_j$ be the $j$-th elementary symmetric expression in $a_1,a_2,\dotsc,a_n$. Then we have
\[ (1+\sigma_1+\sigma_2+\dotsb+\sigma_n)\bs{*}\sffa_r(\bs{b})=\sfJ_{n,r}(\bs{a},\bs{b}), \]
where $\bs{*}$ denotes the product of power series in the sense of {\rm\cite{Gr2}} Chapter I \S3 {\rm(}1.16{\rm)\bb(}1.17bis{\rm)} defined by regarding $a_j$ and $b_j$ as of degree $1$.
\end{enumerate}
\end{lem}
\begin{pf}
(1) follows from the fact that $\sfF_{\!n,r}(\bs{a},\bs{b})-1$ and $\sffa_r(\bs{b})-1$ are both divisible by $b_1 b_2 \dotsb b_r$.
The assertion (2) follows from the definition of the $\bs{*}$-product.
\end{pf}
\begin{defn}\label{def8-1}
{\rm Let $\sfh_{n,r}(t_0,\bs{a},\bs{b})$ be the homogeneous component of degree $n+r$, with respect to $\bs{a}$ and $\bs{b}$, of the power series
\[ \sfJ_{n,r}(\bs{a},\bs{b}) \cdot \sffa_r(\bs{b})^{t_0} =
  \sfF_{\!n,r}(\bs{a},\bs{b}) \cdot \sffa_r(\bs{b})^{t_0-n}
 \in \bZ[t_0]^\scr\psr{\bs{a},\bs{b}}, \]
where $t_0$ is of degree $0$, and $(1+x)^y$ means the binary power series
\[ \sum_{i \ge 0} \ \begin{pmatrix} y \\ i \end{pmatrix} x^i = 1+y\cdot x+\frac{y(y-1)}{2}\cdot x^2+\frac{y(y-1)(y-2)}{3!}\cdot x^3+\dotsb \in \bZ[y]^\scr\psr{x}. \]
By Lemma \ref{lem8-1}\,(1), $\sfh_{n,r}(t_0,\bs{a},\bs{b})$ is divisible by $b_1b_2 \dotsb b_r$. Finally, we define the desired polynomial $\sfp_{n,r}$ as that satisfying
\[  s_r \cdot \sfp_{n,r}(t_0,\sigma_1,\sigma_2,\dotsc,\sigma_n;s_1,s_2,\dotsc,s_r) = \sfh_{n,r}(t_0,\bs{a},\bs{b}), \]
where $\sigma_j$ is as in Lemma \ref{lem8-1}\,(2), and $s_j$ denotes the $j$-th elementary symmetric expression in $b_1,b_2,\dotsc,b_r$. Note that $\sfp_{n,r}(t_0,t_1,\dotsc,t_n;u_1,u_2,\dotsc,u_r)$ is weighted homogeneous of degree $n$, provided that $t_j$ and $u_j$ are of degree $j$.
}
\end{defn}
\begin{prop}\label{prop8-2}
Consider a series $1+\tau_1+\tau_2+\dotsb = 1+\sum_{i \ge 1} \ \tau_i$, and suppose that $\tau_j$ and $u_j$ are of degree $j$. Then we have
\begin{multline*}
 (1+u_r\cdot\sfG_r(u_1,\dotsc, u_r))^{t_0} \cdot 
 \big\{(1+\tau_1+\tau_2+\dotsb ) \bs{*} (1+u_r\cdot\sfG_r(u_1,\dotsc,u_r))\big\} \\
 = 1+ u_r \sum_{n \geq r} \ \sfp_{n-r,r}(t_0,\tau_1,\dotsc,\tau_{n-r};u_1,\dotsc,u_r)
\end{multline*}
Here $\bs{*}$ denotes the product considered in Lemma \ref{lem8-1}\,{\rm(}2{\rm)}.
Consequently, the weighted homogeneous component of degree $n$ of the left hand side is $u_r \cdot \sfp_{n-r,r}(t_0,\tau_1,\dotsc,\tau_{n-r};u_1,\dotsc,u_r)$ {\rm(}resp.\ zero{\rm)} for $n \ge r$ {\rm(}resp.\ for $1 \le n < r${\rm)}.
\end{prop}
\begin{pf}
The assertion is a consequence of Lemma \ref{lem8-1}. The details are straight-forward and left to the reader.
\end{pf}


\section{Riemann-Roch theorem without denominators}\label{sect9}
Let $\cC$ and $\sfGa(*)$ be as in \S\ref{sect4}, and let $\cR$ be the direct sum of $\cK(\sfGa(i),2i)$ with $i \ge 0$, which is a commutative graded ring object with unity in $\cHs(\cC_\Zar)$
 by the functoriality of $\cK(-,2i)$ and the assumption that $\sfGa(*)$ is a graded cohomology theory, cf.\ Definition \ref{axiom0}.
For $Z \in \Ob(\cC)$, let $\cR_Z$ be the restriction of $\cR$ onto $Z_\Zar$, which is a commutative graded ring object with unity in $\cHs(Z_\Zar)$.
\par
Let $f : Y \hra X$ be a regular closed immersion of codimension $r$ which belongs to $\cC$. Using the universal polynomial $\sfp_{n,r}$ constructed in \S\ref{sect8}, we define
\begin{align*}
\sfP_{\!n,Y/X}:=\sfp_{n,r}(\rk,\sfC_1,\sfC_2,\dotsc,\sfC_n;\sfc_1(N_{Y/X}),\sfc_2(N_{Y/X}),\dotsc,\sfc_r(N_{Y/X})) \qquad \\
\in \Mor_{\cHs(Y)}(\OBQPs_Y,\cK\!(\sfGa(n)_Y,2n)) \end{align*}
if $n \ge 0$. Here $\sfc_i(N_{Y/X}) \in \tH^{2i}(Y,\sfGa(i))$ denotes the $i$-th Chern class of the normal bundle $N_{Y/X}$, and we have taken the polynomial of Chern classes with respect to the ring structure on $\cR_Y$. See Definition \ref{def5-1}\,(2) and \eqref{eq9-1} for $\rk$ and $\sfC_i$, respectively. Note that $\sfP_{\!n,Y/X}$ is well-defined, because the rank function has integral values. We define $\sfP_{\!n,Y/X}$ as zero if $n < 0$. The main aim of this section is to prove a local version of Riemann-Roch theorem without denominators:
\begin{thm}\label{thm9-1}
Let $f : Y \hra X$ be as above, and assume that $X$ and $Y$ are both regular and that $Y$ has pure codimension $r \ge 1$ on $X$. Assume further the following condition{\rm:}
\begin{enumerate}
\item[{\rm($\#$)}]
The blow-up of $X \times \bP^1 := X \times_{\Spec(\bZ)} \bP^1_{\bZ}$ along $Y \times \{\infty\}$ belongs to $\Mor(\cC)$.
\end{enumerate}
Then the following diagram commutes in $\cHs(Y)$ for any $i \ge 1${\rm:}
\[\xymatrix{ \OBQPs_Y \ar@{=}[r]^-{f_*}_-\sim \ar[d]_{\sfP_{\!i-r,Y/X}} & \tR f^!\OBQPs_{\!X}
 \ar[d]^{\sfC_{i,X}^Y:=\tR f^!(\sfC_i)} \\
 \cK(\sfGa(i-r)_Y,2(i-r)) \ar[r]^-{f_\push} & \tR f^!\cK(\sfGa(i)_X,2i), }\]
where the upper horizontal arrow denotes the canonical isomorphism due to Quillen {\rm(\cite{Q}} \S7{\rm)}.
\end{thm}
\begin{rem}
Theorem \ref{thm9-1} is a generalization of a theorem of Gillet {\rm\cite{Gi}} Theorem 3.1. However, we have to note that his proof relies on an incorrect formula $\wt{\sfC}^Y\!(\cO_Y)$ $= j_\push(\lam_{-1}\wt{\sfC}(N))$ under the notation in loc.\ cit. Compare with \eqref{eq9-6} and \eqref{eq9-5} below.
\end{rem}
\begin{pf}
We prove Theorem \ref{thm9-1} in two steps by revising Gillet's arguments in \cite{Gi} \S3.
\par\medskip
{\bf Step 1}.
Assume that $f$ is isomorphic to the zero section of the projective completion $\pi : \bP(E \oplus \bs{1}) \to Y$ of a vector bundle $E \to Y$ of rank $r$, where $\bs{1}$ denotes the trivial line bundle over $Y$. In this step
we prove that the following diagram commutes in $\cHs(Y)$:
\begin{equation}\label{eq9-2}
\xymatrix{ \OBQPs_Y \ar[r]^-{f_*} \ar[d]_{\sfP_{\!i-r,Y/X}}
 & \tR\pi_*\OBQPs_{\!X} \ar[d]^{\sfC_{i,X}} \\
 \cK(\sfGa(i-r)_Y,2(i-r)) \ar[r]^-{f_\push}
 & \tR\pi_*\cK(\sfGa(i)_X,2i), }
\end{equation}
which is weaker than the theorem.
For $Z \in \Ob(\cC)$, let $\cR_Z$ be as we defined in the beginning of this section. Let $\wh{\cR_Z}{}^+$ be the product of $\cK(\sfGa(i)_Z,2i)$ with $i \ge 1$, and put 
\[ \vL_Z := \bZ \times \{1\} \times \wh{\cR_Z}{}^+ \in \Ob(\SS(Z_\Zar)), \]
which we endow with the ring structure associated with $\cR_Z$ to obtain a ring object in $\cHs(Z)$ (\cite{Gr2} Chapter I \S3). There is a commutative diagram in $\cHs(Y)$
\begin{equation}\label{eq9-3}
 \xymatrix{\OBQPs_Y \ar[rrd]_-{f_*} \ar[rr]^-{\pi^\back} && \tR\pi_*\OBQPs_{\!X} \ar[d]^-{?\, \bullet \,f_*(1_{\tK})}\ar[rr]^{\wt{\sfC}_X} && \tR\pi_*\vL_X \ar[d]^{? \tstar\, \wt{\sfC}_X(f_*(1_{\tK}))} \\
 && \tR\pi_*\OBQPs_{\!X}   \ar[rr]^{\wt{\sfC}_X} && \tR\pi_*\vL_X\,,}
\end{equation}
where $1_{\tK}$ denotes the unity of $\tK_0(Y)$, and $\scstar$ (resp.\ $\scbullet$) denotes the product structure on $\vL_X$ (resp.\ $\OBQPs_{\!X}$, cf.\ \cite{Gi} (2.31)). See Definition \ref{def5-1}\,(2) for $\wt{\sfC}_X$. The square commutes by Theorem \ref{thm4-1}\,(4) and loc.\ cit.\ Lemma 2.32. The triangle commutes by the projection formula for {\it K}-theory, cf.\ \cite{TT} Proposition 3.17. Let us remind here the following formulas:
{\allowdisplaybreaks
\begin{align}
 f_*(1_{\tK}) & = \lam_{-1}(Q^\vee) \quad \hbox{ in } \;\; \tK_0(X),
 & \hbox{(\cite{FL} Chapter V Lemma 6.2)}  \label{eq9-6} \\
 f_\push(1_\sfGa) & = \sfc_r(Q) \quad \hbox{ in } \;\; \tH^{2r}(X,\sfGa(r)),
 & \hbox{(Corollary \ref{cor7-2})} \label{eq9-5}
\end{align}
}where $Q$ denotes the universal quotient bundle $\pi^*(E \oplus \bs{1})/(L^\taut)^\vee$ on $X$ and $1_\sfGa$ denotes the unity of $\tH^0(Y,\sfGa(0))$. We have used the fact that a regular (noetherian and separated) scheme admits an ample family of invertible sheaves \cite{SGA6} Expos\'e \II\, Corollaire 2.2.7.1, in applying Corollary \ref{cor7-2}. By \eqref{eq9-3} and \eqref{eq9-6}, $\sfC_{i,X} \circ f_*$ agrees with the composite
\begin{align*}
\xymatrix{\OBQPs_Y \ar[r]^-{\pi^\back} & \tR\pi_*\OBQPs_{\!X} \ar[r]^-{\wt{\sfC}_X} & \tR\pi_*\vL_X \ar[rrr]^-{\tstar\,\wt{\sfC}_X(\lam_{-1}(Q^\vee))} &&& \tR\pi_*\vL_X} \\
\xymatrix{\ar[r]^-{\pr_i} & \tR\pi_*\cK(\sfGa(i)_X,2i)\,, }
\end{align*}
where $\pr_i$ denotes the natural projection. Hence we have
\[ \sfC_{i,X} \circ f_* = \pr_i \circ (\wt{\sfC}_X \bs{\star} \wt{\sfC}_X(\lam_{-1}(Q^\vee)) \circ \pi^\back \]
Noting that $\rk(\lam_{-1}(Q^\vee))=0$, we have
{\allowdisplaybreaks
\begin{align*}
\wt{\sfC}_X \bs{\star} \wt{\sfC}_X(\lam_{-1}(Q^\vee)) &
 =(\rk,\sfC_X) \bs{\star} (0,\sfC_X(\lam_{-1}(Q^\vee))) & \\
 & = \big(\rk,\sfC_X(\lam_{-1}(Q^\vee))^\rk \cup \{\sfC_X *\sfC_X(\lam_{-1}(Q^\vee))\}\big)& \hbox{(definition of $\bs{\star}$)} \\
 & = \bigg(\rk,1+\sfc_r(Q) \cup \sum_{j \ge r} \ \sfp_{\!j-r,r}(\rk,\sfC_X;Q) \bigg) &
 \hbox{(Prop.\ \ref{prop8-1}, \ref{prop8-2})}
\end{align*}
}where $\sfC_X$ denotes the total Chern class $1+ \sum_{j \ge 1} \, \sfC_{j,X}$, and we put
\begin{align*}
 \sfp_{\!j-r,r}(\rk,\sfC_X;Q):=\sfp_{\!j-r,r}(\rk,\sfC_{1,X},\dotsc,\sfC_{j-r,X};\sfc_1(Q),\dotsc,\sfc_r(Q)) \qquad \\
  \in \Mor_{\cHs(X)}(\OBQPs_X,\cK\!(\sfGa(j-r)_X,2(j-r))). \end{align*}
Thus we have
{\allowdisplaybreaks
\begin{align*}
\sfC_{i,X} \circ f_* & = (\sfc_r(Q) \cup \sfp_{\!i-r,r}(\rk,\sfC_X;Q)) \circ \pi^\back \\
 & =(f_\push(1_\sfGa) \cup \sfp_{\!i-r,r}(\rk,\sfC_X;Q)) \circ \pi^\back & \hbox{(by \eqref{eq9-5})} \\
 & = f_\push \circ f^\back\sfp_{\!i-r,r}(\rk,\sfC_X;Q) \circ f^\back \circ \pi^\back & \hbox{(the axiom \ref{axiom1}\,(3b))} \\
 & = f_\push \circ \sfP_{\!i-r,Y/X} & \hbox{(\eqref{eq5-4}, $f^*Q \cong N_{Y/X}$)}
\end{align*}
}and the diagram \eqref{eq9-2} commutes.
\par\medskip
{\bf Step 2.}
We prove the theorem using the result of Step 1 and deformation to normal bundle. Let $t : \Spec(\bZ) \to \bP^1_{\bZ}=:\bP^1$ be a morphism of schemes, and consider the following commutative diagram of schemes:
\[\xymatrix{ Y & \ar[l]_-{\;\,p} Y \times \bP^1 \, \ar@{}[rd]|{\square} \ar@<-1pt>@{^{(}->}[r]^-\xi \ar@/^7mm/[rr]^-{h:=f \times \id} & \sfM \ar[r]^-\varrho \ar@{}[rd]|{\square} & X \times \bP^1 \ar[r]^-{\pr_2} \ar@{}[rd]|{\square} & \bP^1 \\
 & \ar@{=}[lu] Y \ar[u] \, \ar@<-1pt>@{^{(}->}[r]^-{\xi_t} \ar@/_7mm/[rr]^-f & \sfM_t \ar[r]^-{\varrho_t} \ar[u] & X \ar[u]_{\phi_t} \ar[r] & \Spec(\bZ). \ar[u]_{t}  }\]
Here $\varrho$ denotes the blow-up of $X \times \bP^1$ along $Y \times \{\infty \}$, which is projective in our sense  because $X \times \bP^1$ is regular and admits an ample family of invertible sheaves. Note also that $\varrho$ is a morphism in $\cC$ by assumption. The arrow $\xi$ is a closed immersion induced by $h$, where we have used the fact that $Y \times \{ \infty \}$ is an effective Cartier divisor on $Y \times \bP^1$. The vertical arrows are morphisms induced by $t$. The arrow $\xi_t$ (resp.\ $\varrho_t$) is the base-change of $\xi$ (resp.\ $\varrho$). Note that we have
\[ \sfM_t \cong \begin{cases} X \quad & \hbox{($t \ne \infty$)} \\ \bP(N_{Y/X} \oplus \bs{1}) \cup \tX \quad &  \hbox{($t = \infty$)}, \end{cases} \]
where $\tX$ denotes the blow-up of $X$ along $Y$. In particular, $\xi_t$ (resp.\ $\varrho_t$) is identical to $f$ (resp.\ $\id_X$) when $t \ne \infty$. As for the case $t=\infty$, it is well-known that $\bP(N_{Y/X} \oplus \bs{1})$ meets $\tX$ along the infinite hyperplane $\bP(N_{Y/X})$ and that $\xi_\infty$ factors through the zero section
\[ s : Y \lra \bP(N_{Y/X} \oplus \bs{1})=:\bP_\infty \qquad \quad \hbox{($s(Y) \cap \bP(N_{Y/X}) = \emptyset$)} \]
(see e.g.\ \cite{FL} Chapter \IV \ \S5). Let $g$ be the restriction of $\varrho_\infty$ to $\bP_\infty$. We will prove the equality of morphisms
\begin{equation}\label{eq9-7}
 \sfC_{i,X}^Y \circ  f_* = g_\push \circ \sfC_{i,\,\bP_\infty}^Y \circ s_* :  \OBQPs_Y \lra \tR f^!\cK(\sfGa(i)_X,2i)
\end{equation}
in $\cHs(Y)$, where the right hand side means the composite morphism
\[ \OBQPs_Y \os{s_*}\lra \tR s^!\OBQPs_{\bP_\infty} \os{\sfC_n}\lra \tR s^!\cK(\sfGa(i)_{\bP_\infty},2i) \os{g_\push}\lra \tR f^!\cK(\sfGa(i)_X,2i). \]
We first check that \eqref{eq9-7} implies the theorem. Indeed, we have
{\allowdisplaybreaks
\begin{align*}
 \sfC_{i,X}^Y \circ f_*
 & = g_\push \circ \sfC_{i,\,\bP_\infty}^Y \circ s_* & \hbox{(by \eqref{eq9-7})} \\
 & = g_\push \circ s_\push \circ \sfP_{\!i-r,Y/\bP_\infty} & \hbox{(by \eqref{eq9-2})} \\
 & = f_\push \circ \sfP_{\!i-r,Y/X} & \hbox{(Theorem \ref{thm7-1}\,(3), $N_{Y/X} \cong N_{Y/\bP_\infty}$)}
\end{align*}
}as claimed.
\par
We prove \eqref{eq9-7} in what follows. Noting that $\sfM$ belongs to $\Ob(\cC)$ by the assumption ($\#$), consider composite morphisms
{\allowdisplaybreaks
\begin{align*}
\alpha : \OBQPs_{Y\times \bP^1} \os{\xi_*}\to \tR \xi^!\OBQPs_{\sfM}
 \os{\sfC_{i,\sfM}^{Y \times \bP^1}}\lra \tR\xi^!\cK(\sfGa(i)_{\sfM},2i)
 \os{\varrho_\push}\to \tR h^!\cK(\sfGa(i)_{X\times \bP^1},2i) , \\
\beta_t : \OBQPs_Y \os{p^\back} \to \tR p_*\OBQPs_{Y \times \bP^1}
 \os{\alpha}\to  \tR p_*\tR h^!\cK(\sfGa(i)_{X\times \bP^1},2i)
 \os{\phi_t^\back}\to \tR f^!\cK(\sfGa(i)_X,2i). \end{align*}
}One can easily check that
\[ \beta_0 = \sfC_{i,X}^Y \circ  f_* \quad \hbox{ and } \quad \beta_\infty = g_\push \circ \sfC_{i,\,\bP_\infty}^Y \circ s_* \]
by the functoriality of Chern class maps (cf.\ Proposition \ref{prop6-1}\,(1)) and the base-change property in Theorem \ref{thm7-1}\,(4). Moreover, we have $\phi_0^\back=\phi_\infty^\back$, i.e., $\beta_0=\beta_\infty$. Indeed, we have
\[ \tR p_*\tR h^!\cK(\sfGa(i)_{X\times \bP^1},2i) \cong \bigoplus_{j=0,1} \ \tR f^!\cK(\sfGa(i-j)_X,2(i-j)) \]
by the projective bundle formula, cf.\ the axiom \ref{axiom1}\,(2), and the pull-back of the tautological line bundle over $X\times \bP^1$ onto $X \times (\bP^1\ssm \{t\})$ is trivial for any $t \in \bP^1(\bZ)$, which imply that both $\phi_0^\back$ and $\phi_\infty^\back$ agree with the natural projection. Thus we have
\[  \sfC_{i,X}^Y \circ  f_* = \beta_0 = \beta_\infty = g_\push \circ \sfC_{i,\bP_\infty}^Y \circ s_*\,, \]
which completes the proof of the theorem.
\end{pf}
\begin{cor}\label{cor9-1}
Under the setting of Theorem \ref{thm9-1}, the diagram
\[\xymatrix{\tK_j(Y) \ar@{=}[r]^-{f_*}_-\sim \ar[d]_{\sfP_{i-r,Y/X,j}} & \tK_j^Y(X) \ar[d]^{\sfC_{i,j,X}^Y} \\
 \tH^{2(i-r)-j}(Y,\sfGa(i-r)) \ar[r]^-{f_\push} & \tH^{2i-j}_Y(X,\sfGa(i)) }\]
is commutative for any $i \ge 1$ and $j \ge 0$, where $\sfP_{i-r,Y/X,j}$ denotes the composite map
\[ \begin{CD} \tK_j(Y) @>{\sfP_{i-r,Y/X}}>> \tH^{-j}(Y,\cK\!(\sfGa(i-r)_Y,2(i-r))) \os{\eqref{lem6-1}}\cong \tH^{2(i-r)-j}(Y,\sfGa(i-r)). \end{CD}\]
\end{cor}

\section{Grothendieck-Riemann-Roch theorem}\label{sect10}
Let $\cC$ and $\sfGa(*)$ be as in \S\ref{sect4}. In this section, we prove Theorem \ref{thm10-1} below. Let $f : Y \to X$ be a projective morphism in $\cC$, and suppose that $X$ and $Y$ are regular. See Definition \ref{def6-3} for $\tK_*(Y)$ and $\whH^*(Y,\sfGa(\scbullet))_\bQ$. By taking a factorization
\[\xymatrix{ f : Y \;\ar@<-1pt>@{^{(}->}[r]^-g & \bP:=\bP^m_X \ar[r] & X }\]
with $g$ a closed immersion, we define {\it the virtual tangent bundle} $T_f$ of $f$ as
\[ T_f:=[g^*T_{\bP/X}]-[N_{Y/\bP}] \in \tK_0(Y), \]
which is independent of the factorization of $f$, cf.\ \cite{FL} Chapter V Proposition 7.1. Note that the relative dimension defined in \S\ref{sect7} is exactly the virtual rank of $T_f$. We further define {\it the Todd class} $\Td(T_f)$ $\in \whH^*(Y,\sfGa(\scbullet))_\bQ$ as $\Td(T_f) := \Td(g^*T_{\bP/X})/\Td(N_{Y/\bP})$,
which is independent of the factorization of $f$ as well.
\begin{thm}\label{thm10-1}
Let $f : Y \to X$ be a projective morphism in $\cC$ with both $X$ and $Y$ regular. Assume the following condition{\rm:}
\begin{enumerate}
\item[{\rm($\#'$)}]
$f : Y \to X$ is isomorphic to a projective space over $X$, or there exists a decomposition $Y \hra \bP^m_X \to X$ of $f$ such that the blow-up of $\bP^m_X \times_X \bP^1_X$ along $Y \times \{\infty\}$ belongs to $\Mor(\cC)$.
\end{enumerate}
Then the diagram
\[\xymatrix{
\tK_*(Y) \ar@{=}[r]^-{f_*}_-\sim \ar[d]_{\ch_Y(-) \, \cup \, \Td(T_f)} & \tK_*^{f(Y)}(X) \ar[d]^{\ch_X^{f(Y)}} \\
\whH^*(Y,\sfGa(\scbullet))_\bQ \ar[r]^-{f_\push} & \whH^*_{f(Y)}(X,\sfGa(\scbullet))_\bQ
}\]
is commutative, that is, for $\alpha \in \tK_*(Y)$ we have
\[ \ch_X^{f(Y)}(f_*\alpha)=f_\push (\ch_Y(\alpha) \cup \Td(T_f)) \quad \hbox{ in } \;\; \whH^*_{f(Y)}(X,\sfGa(\scbullet))_\bQ. \]
Here $f_\push$ denotes the push-forward morphism constructed in \S\ref{sect7}.
\end{thm}

\begin{pf}
When $Y$ is a projective space over $X$, then the assertion follows from the projective bundle formula (the axiom \ref{axiom1}\,(2)) and the arguments in \cite{FL} Chapter \II\, Theorem 2.2 (see also loc.\ cit.\ Chapter V\, Theorem 7.3). Hence by loc.\ cit.\ Chapter \II\, Theorem 1.1 and the same arguments as in Step 2 of the proof of Theorem \ref{thm9-1} (we need the assumption ($\#'$) here), we have only to check the commutativity of the diagram
\begin{equation}\label{eq10-1}
\xymatrix{
\tK_*(Y) \ar[r]^-{f_*} \ar[d]_{\ch_Y(-) \, \cup \, \Td(T_f)} & \tK_*(X) \ar[d]^{\ch_X} \\
\whH^*(Y,\sfGa(\scbullet))_\bQ \ar[r]^-{f_\push} & \whH^*(X,\sfGa(\scbullet))_\bQ
}\end{equation}
assuming that $X\cong \bP(E \oplus \bs{1})$, the projective completion of a vector bundle $E \to Y$ of rank $r$ and that $f$ is isomorphic to the zero section of  $\pi : \bP(E \oplus \bs{1}) \to Y$. Let $Q$ be the universal quotient bundle over $X$. Then we have
\begin{align*}
\ch(f_*(1_{\tK})) & = \ch(\lam_{-1}(Q^\vee)) & \hbox{(by \eqref{eq9-5})} \\
 & = \sfc_r(Q) \cup \Td(Q)^{-1} & \hbox{(\cite{FL} Chapter I Proposition 5.3)} \\
 & = f_\push(1_\sfGa \cup f^\back\Td(Q)^{-1}) & \hbox{(\eqref{eq9-6}, the axiom \ref{axiom1}\,(3b))} \\
 & = f_\push(\Td(T_f)) & \hbox{($f^\back[Q] = [N_{Y/X}] = -T_f$)}
\end{align*}
in $\whH^*(X,\sfGa(\scbullet))_\bQ$. Therefore the diagram \eqref{eq10-1} commutes by the arguments in loc.\ cit.\ Chapter \II\, Theorem 1.2 and the projection formula for \tK-theory and $\sfGa(*)$-cohomology.
\end{pf}


\section{Computation via $\bs{1}$-extension}\label{sect11}
\def\ssm{\smallsetminus}
\def\sing{{\sf sing}}
\def\Supp{{\sf Supp}}
As an application of the Riemann-Roch theorem without denominators, we compute Chern classes using $1$-extensions.

\par
Let $\cC$ be as in \S\ref{ex:betti}. Let $X \in \Ob(\cC)$ be a {\it proper smooth} variety over $\bC$ and let $Z \subset X$ be a reduced closed subscheme of pure codimension $r \ge 0$. Let $\sfGab(*)$ be the Betti complex with $A=\bQ$ (cf.\ \S\ref{ex:betti}) and let $\sfGad(*)$ be the Deligne-Beilinson complex with $A=\bQ$ (cf.\ \S\ref{ex:deligne}), which are both admissible cohomology theories on $\cC$. We are concerned with the Chern class maps
\begin{align*}
& \sfC^{\,\sing,Y}_{i,j} : \tK^Y_j(X) \lra \tH^{2i-j}_Y(X,\sfGab(i)), \\
& \sfC^{\cD,Y}_{i,j} : \tK^Y_j(X) \lra \tH^{2i-j}_Y(X,\sfGad(i)),
\end{align*}
where $Y$ is either $Z$ or $X$. We assume that $j \ge 1$ for simplicity (see Remark \ref{rem11-1} below for the case $j=0$), and put
\begin{align*}
 V & :=\Coker(H^{2i-j-1}_Z(X,\sfGab(i)) \lra H^{2i-j-1}(X,\sfGab(i))).
\end{align*}
We are going to compute the composite map
\begin{align}
\notag \varrho^\cD_{i,j} : \tK_j^Z(X) & \os{\sfC^{\cD,Z}_{i,j}}{\lra} \tH^{2i-j}_Z(X,\sfGad(i)) \lra \tH^{2i-j}(X,\sfGad(i)) \\
\notag & \; \os{\alpha_X}\cong \; \Hom_{\sfD(\MHS)}(\bQ,\tR\vG(X,\sfGab(i))[j]) \\
\label{eq11-1}\tag{$\spadesuit$} & \; \os{\tau}\cong \; \Ext^1_{\MHS}(\bQ,\tH^{2i-j-1}(X,\sfGab(i))) \lra \Ext^1_{\MHS}(\bQ,V)
\end{align}
in terms of the Chern class map of a regular dense open subset of $Z$, where $\MHS$ denotes the category of rational mixed Hodge structures and $\alpha_X=\alpha_{X,X}$ denotes the canonical isomorphism given in the following lemma. The isomorphism $\tau$ is obtained from the fact that $2$-extensions are trivial in $\MHS$ and the assumption that $j \ge 1$.
\begin{lem}\label{lem11-1}
Let $X \to \Spec(\bC)$ be as before, and let $T$ be a closed subscheme of $X$. Then for integers $i,j \ge 0$, there exists a canonical isomorphism
\[ \alpha_{X,T} : \tH^j_T(X,\sfGad(i)) \os{\simeq}\lra \Hom_{\sfD(\MHS)}(\bQ,\tR\vG_T(X,\sfGab(i))[j]) \]
fitting into a commutative diagram
\[\xymatrix{
\tH^j_T(X,\sfGad(i)) \ar[r] \ar@{=}[d]_{\alpha_{X,T}}^\wr & \tH^j(X,\sfGad(i)) \ar@{=}[d]_\wr^{\alpha_X=\alpha_{X,X}} \\
\Hom_{\sfD(\MHS)}(\bQ,\tR\vG_T(X,\sfGab(i))[j]) \ar[r] & \Hom_{\sfD(\MHS)}(\bQ,\tR\vG(X,\sfGab(i))[j]).}\]
\end{lem}
\begin{pf}
For an open subset $U \subset X$, let $\bQ(i)_{\cD,X,U}$ be the Deligne-Beilinson complex of $U$ on the analytic site $X_\an$, cf.\ \cite{EV} Definition 2.6. By the definition of $\sfGad(i)$, there is a natural homomorphism of complexes
\[ \beta_{X,U} : \tR\vG(X_\an,\bQ(i)_{\cD,X,U}) \lra \tR\vG(U_\Zar,\sfGad(i)), \]
which is a quasi-isomorphism by loc.\ cit.\ Lemma 2.8. Let $\MHM(X)$ be the category of mixed Hodge modules on $X$, cf.\ \cite{Sm}. For $\cM=(M, F^\bullet, K_\bQ) \in \Ob(\MHM(X))$, we define a complex $\cM^\dagger$ of abelian sheaves on $X_\an$ as
\[ \cM^\dagger := \cone(K_\bQ[-\dim X] \oplus F^0\DR_X(M) \to \DR_X(M))[-1], \]
where $\DR_X(M)$ and $F^0\DR_X(M)$ denote the complexes (on $X_\an$) with the most left term placed in degree $0$
\begin{align*}
\DR_X(M) & : M \lra M\otimes \Omega^1_X \lra \dotsb \lra M \otimes \Omega^q_X \lra \dotsb, \\
F^0\DR_X(M) & : F^0\! M \lra F^{-1}\! M\otimes \Omega^1_X \lra \dotsb \lra F^{-q}\! M \otimes \Omega^q_X \lra \dotsb,
\end{align*}
respectively. Because the assignment $H \mapsto H^\dagger$ is exact, this induces a functor
\[ (-)^\dagger : \sfD^b(\MHM(X)) \lra \sfD^b(X_\an). \]
Note also that for an open immersion $\psi : U \hra X$ there is a natural quasi-isomorphism of complexes (on $X_\an$)
\[ \bQ(i)_{\cD,X,U} \os{\qis}\lra (\tR \psi_*\bQ(i)_U)^\dagger, \]
where $\bQ(i)_U$ denotes the Hodge module associated to the constant sheaf $(2\pi\sqrt{-1})^i\bQ$ on $U_\an$. There is a diagram of quasi-isomorphisms of complexes
\[ \tR\vG(U_\Zar,\sfGad(i)) \os{\beta_{X,U}}{\us{\qis}\lla} \tR\vG(X_\an,\bQ(i)_{\cD,X,U}) \os{\qis}\lra \tR\vG(X_\an,(\tR \psi_*\bQ(i)_U)^\dagger). \]
Considering this diagram for $U=X$ and $X \ssm T$, we obtain an isomorphism
\[ \gamma_{X,T} : \tR\vG_T(X_\Zar,\sfGad(i)) \cong \tR\vG(X_\an,(\tR \phi^!\bQ(i)_X)^\dagger) \quad \hbox{ in } \;\; \sfD(\Ab), \]
where $\phi$ denotes the closed immersion $T \hra X$. Finally we define $\alpha_{X,T}$ as the composite
\begin{align*}
\tH^q_T(X_\Zar,\sfGad(i)) & \os{\gamma_{X,T}}{\,\;\cong\,\;} \tH^q(X_\an,(\tR \phi^!\bQ(i)_X)^\dagger) \\
& \,\;\cong\;\, \Hom_{\sfD(\MHM(X))}(\bQ_X,\tR \phi^!\bQ(i)_X[q]) \\
& \,\;\cong\;\, \Hom_{\sfD(\MHS)}(\bQ,\tR \vG_T(X_\an,\bQ(i)_X)[q]),
\end{align*}
which obviously fits into the commutative diagram in the lemma.
\end{pf}
\par
We return to the setting of the beginning of this section. Let $\xi_Z$ be an element of $\tK'_j(Z)$, and put
\[ \chi_Z:=\sfC^{\,\sing,Z}_{i,j}(f_*(\xi_Z)) \in H^{2i-j}_Z(X,\sfGab(i)), \]
where $f_*$ denotes the isomorphism $\tK'_j(Z) \cong \tK_j^Z(X)$. There is a localization exact sequence of rational Betti cohomology
\begin{equation}\label{eq11-3}
 0 \to V \to H^{2i-j-1}(X \ssm Z,\sfGab(i)) \os{\delta}\to H^{2i-j}_Z(X,\sfGab(i)) \os{\iota}\to H^{2i-j}(X,\sfGab(i)),
\end{equation}
where $\delta$ (resp.\ $\iota$) denotes the connecting map (resp.\ the canonical map). By the assumption that $j \ge 1$, we have $\iota(\chi_Z)=0$. Hence pulling-back this exact sequence by $\chi_Z$, that is, considering the fiber product in $\MHS$
\[ E:=\big\{(x,a) \in H^{2i-j-1}(X \ssm Z,\sfGab(i)) \times \bQ\,\big|\, \delta(x)=a \cdot \chi_Z \big\}, \]
we obtain a short exact sequence of rational mixed Hodge structures
\[ 0 \lra V \lra E \os{\pr_2}\lra \bQ \lra 0, \]
which we denote by $\eta_Z$. The following results computing $\sfC^\cD_{i,j}(\xi_Z)$ have been used in recent joint papers of Otsubo and the first author \cite{A1}, \cite{A2}:
\begin{thm}\label{thm11-1}
Assume that $i,j \ge 1$ {\rm(}see Remark \ref{rem11-1} below for the case $j=0${\rm)}. Then
\begin{enumerate}
\item[{\rm(1)}]
The map $\varrho^\cD_{i,j}$ in \eqref{eq11-1} sends $f_*(\xi_Z) \in \tK_j^Z(X)$ to the class of $\eta_Z$, up to the sign $(-1)^j$.
\item[{\rm(2)}]
Let $Z_\sing$ be the singular locus of $Z$, and put $Z^\circ:=Z \ssm Z_\sing$ and $X^\circ:=X \ssm Z_\sing$. Assume further that $2i-j < 2(r+1)$. Then the sequence of Betti cohomology 
\begin{equation}\label{eq11-2}
 0 \lra V \lra \tH^{2i-j-1}(X \ssm Z,\sfGab(i)) \os{\delta'}\lra \tH^{2i-j}_{Z^\circ}(X^\circ,\sfGab(i))
\end{equation}
is exact, where $\delta'$ is the composite of $\delta$ in \eqref{eq11-3} and a natural restriction map. Moreover, $\eta_Z$ is isomorphic to the pull-back of this exact sequence by the element
\[ g_\push(\sfP_{i-r,\,Z^\circ/X^\circ}(\xi_{Z^\circ})) \in \tH^{2i-j}_{Z^\circ}(X^\circ,\sfGab(i)), \]
where $g$ denotes the regular closed immersion $Z^\circ \hra X^\circ$ and $\xi_{Z^\circ}$ denotes the restriction of $\xi_Z$ to $\tK'_j(Z^\circ)=\tK_j(Z^\circ)$. See \S\ref{sect9} for $\sfP_{i-r,\,Z^\circ/X^\circ}$.
\end{enumerate}
\end{thm}
\begin{pf}
(1) Put $m:=2i-j$ for simplicity, and consider the following big diagram:
{\small\[\xymatrix{
\tK_j^Z(X) \ar[r]^{a_1} \ar[d]_{\sfC^{\cD,Z}_{i,j}}\ar@{}[rd]|{\rm (i)} & \tK_j(X) \ar[d]^{\sfC^{\cD}_{i,j}} \\
\tH^m_Z(X,\sfGad(i)) \ar[r]^{a_2} \ar@{=}[d]_{\alpha_{X,Z}}^\wr \ar@{}[rd]|{\rm (ii)} & \tH^m(X,\sfGad(i)) \ar@{=}[d]_\wr^{\alpha_X} \\
\Hom_{\sfD(\MHS)}(\bQ,\tR\vG_Z(X,\sfGab(i))[m]) \ar[r]^{a_3} \ar@{}[rd]|{\rm (iii)} & \Hom_{\sfD(\MHS)}(\bQ,\tR\vG(X,\sfGab(i))[m]) \\
\Hom_{\sfD(\MHS)}(\bQ,(\tau_{\le m}\tR\vG_Z(X,\sfGab(i))')[m]) \ar[r]^{a_4} \ar[d]_{\can} \ar@{=}[u]_{\wr}^{(*)} \ar@{}[rd]|{\rm (iv)} & \Hom_{\sfD(\MHS)}(\bQ,(\tau_{\le m-1}\tR\vG(X,\sfGab(i)))[m]) \ar[d]^\can \ar@{=}[u]^{\wr}_{(*')} \\
 \Hom_\MHS(\bQ,\Image(\delta)) \ar[r]^{a_5} & \Ext^1_{\MHS}(\bQ,V). }\]
}Here $\tau_{\le m}\tR\vG_Z(X,\sfGab(i))'$ denotes the complex
\[ \cone\big(\tau_{\le m}\tR\vG_Z(X,\sfGab(i)) \to \tH^m_Z(X,\sfGab(i))/\Image(\delta)[-m]\big)[-1], \]
The arrows $a_1$, $a_2$ are canonical maps, and $a_3$ is induced by the canonical morphism \[ \wt{\iota} : \tR\vG_Z(X,\sfGab(i)) \lra \tR\vG(X,\sfGab(i)). \]
The arrow $a_4$ denotes the morphism induced by $\wt{\iota}$, and $a_5$ denotes the connecting map associated with the short exact sequence in $\MHS$
\[ 0 \lra V \lra H^{m-1}(X \ssm Z,\sfGab(i)) \os{\delta}\lra \Image(\delta) \lra 0. \]
The squares (i) and (iii) commute obviously, and the square (ii) commutes by the construction of $\alpha_{X,Z}$ and $\alpha_X=\alpha_{X,X}$ in Lemma \ref{lem11-1}. The square (iv) commutes up to the sign $(-1)^m$, cf.\ \cite{J2} Lemma 9.5. The isomorphism $(*)$ (resp.\ $(*')$) in the left (resp.\ right) column follows from the fact that the Hodge $(0,0)$-part of $\tH^m_Z(X,\sfGab(i))$ lies in $\Image(\delta)$ (resp.\ the Hodge $(0,0)$-part of $\tH^m(X,\sfGab(i))$ is zero). Finally, the composite of the left vertical columns sends $f_*(\xi_Z)$ to $\chi_Z$, and  the composite of the right vertical columns agrees with the map $\varrho^\cD_{i,j}$ in question. The assertion follows from these facts.
\par
(2) Since we have $H^{2i-j}_{Z_\sing}(X,\sfGab(i))=0$ by the assumption on $i$ and $j$, the restriction map
\[  H^{2i-j}_Z(X,\sfGab(i)) \lra H^{2i-j}_{Z^\circ}(X^\circ,\sfGab(i))\] is injective, which implies the first assertion. The second assertion follows from Corollary \ref{cor9-1} for Betti cohomology.
\end{pf}
\begin{rem}\label{rem11-1}
As for the case $j=0$, one can easily modify the above arguments to obtain similar results for $\xi_Z \in \tK'_0(Z)$ satisfying $\sfC^{\,\sing}_{i,0}(f_*(\xi_Z))=0$ in $H^{2i}(X,\sfGab(i))$.
\end{rem}
\section{Admissible cohomology theory with log poles}\label{sect12}
Let $\cC$ be as in \S\ref{sect4}, let $\cC^\reg$ be the full subcategory of $\cC$ consisting of regular schemes. In this section we formulate a logarithmic variant of admissible cohomology theory, for which we will consider Chern class maps in the next section. We first fix the following terminology:
\begin{defn}\label{defn12-a}
{\rm
\begin{enumerate}
\item[(1)]
{\it A log pair} is a pair $(X,D)$ of a regular scheme $X$ and a simple normal crossing divisor $D$ on $X$ which may be empty. {\it A morphism $f : (Y,E) \to (X,D)$ of log pairs} is a morphism of schemes $f : Y \to X$ which satisfies $f^{-1}(D) \subset E$ as closed subsets of $Y$. {\it A closed immersion $f : (Y,E) \to (X,D)$ of log pairs} is a (regular) closed immersion $f : X \to Y$ with $f^{-1}(D) = E$ as closed subschemes of $Y$.
\item[(2)]
{\it A log pair in $\cC$} is a log pair $(X,D)$ such that $X$ belongs to $\Ob(\cC^\reg)$ and such that $D$ is empty or satisfies the following condition ($\star$):
\begin{enumerate}
\item[($\star$)]
Let $D_1,D_2,\dotsc,D_a$ be the irreducible components of $D$. Then the $n$-fold intersection $D_{i_1} \cap D_{i_2} \cap \dotsc \cap D_{i_n}$ is transversal and belongs to $\Ob(\cC^\reg)$, for any $1 \le n \le a$ and any $n$-tuple $(i_1, i_ 2, \dotsc, i_n)$ with $1 \le i_1 < i_ 2 < \dotsb < i_n \le a$.
\end{enumerate}
{\it A morphism of log pairs in $\cC$} is a morphism $f : (Y,E) \to (X,D)$ of log pairs such that $(X,D)$ and $(Y,D)$ are log pairs in $\cC$ and such that the underlying morphism $f : Y \to X$ belongs to $\Mor(\cC)$.
\item[(3)]
We define $\cC^\log$, {\it the category of log pairs in $\cC$}, to be the category whose objects are log pairs in $\cC$ and whose morphisms are those of log pairs in $\cC$. 
\item[(4)]
{\it An open subset of a log pair $(X,D)$} is a log pair of the form $(U,D \cap U)$, where $U$ is an open subset of $X$.
\item[(5)]
{\it An open covering of a log pair $(X,D)$} is a family $\{(X_i,D_i)\}_{i \in I}$ of open subsets of $(X,D)$ satisfying $\bigcup{}_{i \in I} \ X_i = X$.
\item[(6)]
({\it Big Zariski site})\;
We endow the category $\cC^\log$ with Zariski topology in such a way that a covering of a given object $(X,D)$ is an open covering $\{(X_i,D_i)\}_{i \in I}$ of $(X,D)$. We denote the resulting site by $(\cC^\log)_\Zar$.
\item[(7)]
({\it Small Zariski site})\;
For a log pair $(X,D)$ in $\cC$, let $\Ouv/(X,D)$ be the category of open subsets of $(X,D)$. We endow the category $\Ouv/(X,D)$ with Zariski topology in such a way that a covering of a given object $(U,D \cap U)$ is an open covering $\{(U_i,D \cap U_i)\}_{i \in I}$ of $(U,D \cap U)$. We denote the resulting site by $(X,D)_\Zar$.
\end{enumerate}}
\end{defn}
\noindent
We mention here some elementary facts, which will be used freely in what follows.
\begin{lem}\label{lem12-1}
Let $(X,D)$ be a log pair in $\cC$, and let $U$ be an open subset of $X$. Then{\rm:}
\begin{enumerate}
\item[{\rm(1)}]
The pair $(U,D \cap U)$ is a log pair in $\cC$.
\item[{\rm(2)}]
The natural morphism $(U,D \cap U) \to (X,D)$ is a morphism in $\cC^\log$.
\item[{\rm(3)}]
The natural morphism $(X,D) \to (X,\emptyset)$ induces an equivalence $(X,D)_\Zar \cong X_\Zar$.
\item[{\rm(4)}]
There is a natural functor $\theta_{(X,D)}:\sfS^\ab((\cC^\log)_\Zar) \to \sfS^\ab((X,D)_\Zar)$ induced by the restriction of topology.
\end{enumerate}
\end{lem}
\begin{pf}
(1) and (2) follows from our assumption ($*_1$) in \S\ref{sect1.1} on the category $\cC$. (3) is straight-forward and left to the reader. (4) follows from (1).
\end{pf}

We further fix the following terminology on simplicial schemes:

\begin{defn}
{\rm
\begin{enumerate}
\item[(1)]
{\it A simplicial log pair} in $\cC$ is a simplicial object in $\cC^\log$. {\it A morphism $f : (\Ys,\Es) \to (\Xs,\Ds)$ of simplicial log pairs in $\cC$} is a morphism of simplicial objects in $\cC^\log$.
\item[(2)]
{\it A closed immersion $f : (\Ys,\Es) \to (\Xs,\Ds)$ of simplicial log pairs in $\cC$} is a morphism of simplicial log pairs in $\cC$ such that the $p$-th factor $f_p: (Y_p,E_p) \to (X_p,D_p)$ is a closed immersion in the sense of Definition \ref{defn12-a}\,(1) for each $p \ge 1$.
\item[(3)]
We say that a closed immersion $f : (\Ys,\Es) \to (\Xs,\Ds)$ of simplicial log pairs in $\cC$ is {\it strict}, if the underlying closed immersion $\Ys \to \Xs$ is strict (i.e., exact) in the sense of Definition \ref{def1-2}\,(2).
\item[(4)]
({\it Zariski site})\;
For a simplicial log pair $(\Xs,\Ds)$ in $\cC$, let $\Ouv/(\Xs,\Ds)$ be the category whose object is an open subset $(U,D_p \cap U)$ of $(X_p,D_p)$ for some $p$, and whose morphism is a commutative square in $\cC^\log$
\[\xymatrix{ (V,D_q \cap V) \ar[d] \ar[r] & (U,D_p \cap U) \ar[d] &  \\ (X_q,D_q) \ar[r]^{\alpha^{(\Xs,\Ds)}} & (X_p,D_p) }\]
for some morphism $\alpha : [p] \to [q]$ in $\vD$.
 We endow the category $\Ouv/(\Xs,\Ds)$ with Zariski topology in such a way that a covering of a given object $(U,D_p \cap U) \to (X_p,D_p)$ is an open covering $\{(U_i,D_i)\}_{i \in I}$ of $(U,D_p \cap U)$. We denote the resulting site by $(\Xs,\Ds)_\Zar$.
\end{enumerate}}
\end{defn}
\noindent
The following lemma is a simplicial analogue of Lemma \ref{lem12-1}\,(3) and (4), whose proofs are left to the reader:
\begin{lem}
Let $(\Xs,\Ds)$ be a simplicial log pair in $\cC$.
\begin{enumerate}
\item[{\rm(1)}]
The natural morphism $(\Xs,\Ds) \to (\Xs,\emptyset)$ induces an equivalence $(\Xs,\Ds)_\Zar \cong (\Xs)_\Zar$.
\item[{\rm(2)}]
There is a natural functor $\theta_{(\Xs,\Ds)}:\sfS^\ab((\cC^\log)_\Zar) \to \sfS^\ab((\Xs,\Ds)_\Zar)$ induced by the restriction of topology.
\end{enumerate}
\end{lem}

\par We now formulate a logarithmic version of admissible cohomology theory.
\begin{defn}\label{defn12-1}
{\rm 
{\it An admissible cohomology theory on $\cC^\log$} is a graded cohomology theory $\sfGa(*)^\log=\{\sfGa(i)^\log\}_{i \in \bZ}$ on $(\cC^\log)_\Zar$ satisfying the following axioms. For a simplicial log pair $(\Xs,\Ds)$ in $\cC$, we denote $\theta_{(\Xs,\Ds)}(\sfGa(i)^\log)$ by $\sfGa(i)_{(\Xs,\Ds)}$, which we often regard as complexes of abelian sheaves on $(\Xs)_\Zar$.
\begin{enumerate}
\item[(0)]
For a log pair $(X,D) \in \Ob(\cC^\log)$ and a dense open subset $U$ of $X$, the restriction map $\tH^0(X_\Zar,\sfGa(0)_{(X,D)}) \to \tH^0(U_\Zar,\sfGa(0)_{(X,D)})$ is injective.
\item[(1)]
There exists a morphism
\[ \varrho^\log : \Gm^\log[-1] \lra \sfGa(1)^\log \qquad \hbox{ in } \quad \sfD((\cC^\log)_\Zar), \]
where $\Gm^\log$ denotes the sheaf on $(\cC^\log)_\Zar$ defined as $\Gm^\log(X,D):=\cO^\times(X \ssm D)$.
\item[(2)]
For a scheme $X \in \Ob(\cC^\reg)$ and a vector bundle $E$ on $X$, the complexes $\sfGa(i)_X:=\sfGa(i)_{(X,\emptyset)}$ and $\sfGa(i)_{\bP(E)}:=\sfGa(i)_{(\bP(E),\emptyset)}$ ($i \in \bZ$) satisfy the projective bundle formula in the sense of Definition \ref{axiom1}\,(2).
\item[(3)]
For a strict closed immersion $f : (\Ys,\Es) \hra (\Xs,\Ds)$ of simplicial log pairs in $\cC$ of codimension $r$, there are push-forward morphisms
\[ f_\push :  f_*\sfGa(i)_{(\Ys,\Es)} \lra \sfGa(i+r)_{(\Xs,\Ds)}[2r] \qquad \hbox{($i \in \bZ$)} \]
in $\sfD((\Xs)_\Zar)$ which satisfy the properties analogous to (3a)\bb(3d) in Definition \ref{axiom1}\,(3).
\item[(4)]
Let $(X,D)$ be a log pair in $\cC$, and let $Y$ be an irreducible component of $D$. Put $D' := \ol{D \ssm Y}$ (Zariski closure in $X$) and $E:=Y \cap D'$, and let
\[ f : (Y,E) \lra (X,D') , \qquad j : (X,D) \lra (X,D'), \]
be the natural morphisms in $\cC^\log$. Then there exists a distinguished triangle
\[\xymatrix{ f_*\sfGa(i-1)_{(Y,E)}[-2] \ar[r]^-{f_\push} & \sfGa(i)_{(X,D')} \ar[r]^-{j^\back} & \sfGa(i)_{(X,D)} \ar[r]^-\delta & f_*\sfGa(i-1)_{(X,D')}[-1] }\]
in $\sfD(X_\Zar)$ for any $i \in \bZ$.
\end{enumerate}}
\end{defn}

\begin{exmp}\label{exmp12-1}
{\rm For $\cC$ and $\sfGa(*)$ is as in \S\S\ref{ex:etale}\bb\ref{ex:deligne}, let $G(i)$ be the Godement resolution of $\sfGa(i)$ on $\cC_\Zar$. Then the assignment
\[ (X,D) \in \Ob(\cC^\log) \longmapsto j_*G(i)_{X \ssm D} \qquad \hbox{($i \in \bZ$)} \]
defines an admissible cohomology theory $\sfGa(*)^\log$ on $\cC^\log$, where $j$ denotes the open immersion  $X \ssm D \hra X$. For $\cC$ as in \S\S\ref{ex:deRham}\bb\ref
{ex:ptate}, we have an admissible cohomology theory $\sfGa(*)^\log$ on $\cC^\log$ as follows.
\begin{enumerate}
\item[(1)]
For $\cC$ as in \S\ref{ex:deRham}, we define $\sfGa(*)^\log$ on $\cC^\log$ as
\[ \sfGa(i)_{(X,D)}:=
\begin{cases}
 \Omega^\bullet_{X/k}(\log D) \qquad & \hbox{($i \ge 0$)} \\
 0  \qquad & \hbox{($i < 0$)},
\end{cases} \]
which is an admissible cohomology theory on $\cC^\log$.
\item[(2)]
For $\cC$ as in \S\ref{ex:loghw} and a positive integer $n$, we define $\sfGa(*)^\log$ on $\cC^\log$ as
\[ \sfGa(i)_{(X,D)}:=
\begin{cases}
 \vare_*G(i)^\bullet_{(X,D)} [-i] \qquad & \hbox{($i \ge 0$)} \\
 0 \qquad & \hbox{($i < 0$)},
\end{cases}
\]
where $\vare$ denotes the natural morphism of sites $X_\et \to X_\Zar$ and $G(i)^\bullet_{(X,D)}$ denotes the Godement resolution of the \'etale sheaf $\tau_*\logwitt {X \ssm D} n i$ with $\tau : X \ssm D \hra X$ the natural open immersion. Then $\sfGa(*)^\log$ is an admissible cohomology theory on $\cC^\log$. The axiom \ref{defn12-1}\,(4) follows from the purity theorem in \cite{Sh} Theorem 3.2.
\item[(3)]
For $\cC$ as in \S\ref{ex:ptate} and a positive integer $n$, we define $\sfGa(*)^\log$ on $\cC^\log$ as
\[ \sfGa(i)_{(X,D)}:= \vare_*G(i)^\bullet_{(X,D)} \qquad \hbox{($i \in \bZ$)}, \]
where $\vare$ denotes the natural morphism of sites $X_\et \to X_\Zar$, and $G_n(i)_{(X,D)}$ denotes the Godement resolution of the $i$-th \'etale Tate twist $\fT(i)_{(X,D)}$ with log poles along $D$ (\cite{Sa2} Corollary 3.9). Then $\sfGa(*)^\log$ is an admissible cohomology theory on $\cC^\log$. The axiom \ref{defn12-1}\,(4) is a variant of loc.\ cit.\ Theorem 3.12.
\end{enumerate}
}
\end{exmp}

In the rest of this section, let $\sfGa(*)^\log$ be an admissible cohomology theory on $\cC^\log$. The following fact is obvious:
\begin{lem}
We regard the category $\cC^\reg$ as a full subcategory $\cC^\log$ by sending a scheme $X \in \Ob(\cC^\reg)$ to the log pair $(X,\emptyset)$. Then the restriction of $\sfGa(*)^\log$ onto $(\cC^\reg)_\Zar$ is an admissible cohomology theory on $\cC^\reg$ in the sense of Definition \ref{axiom1}.
\end{lem}
The following proposition gives a logarithmic variant of the projective bundle formula in Definition \ref{axiom1}\,(2).
\begin{prop}\label{prop12-1}
Let $\sfGa(*)^\log$ be an admissible cohomology theory on $\cC^\log$.
Let $(X,D)$ be a log pair in $\cC$, and let $E$ be a vector bundle over $X$ of rank $r+1$. Let $\pi : \bP(E) \to X$ be the projective bundle associated with $E$, and let $D'$ be the divisor $\pi^{-1}(D)$ on $\bP(E)$. Let
\[ p : (\bP(E),D') \to (X,D) \qquad \hbox{{\rm(}resp.\ $\beta : (\bP(E),D')) \to (\bP(E),\emptyset)${\rm)}} \]
be the morphism of log pairs induced by $\pi$ {\rm(}resp.\ the natural morphism of log pairs{\rm)}.
Then we have
\begin{align*}
 \gamma_F : \;& \bigoplus_{j=0}^r \ \sfGa(i-j)_{(X,D)}[-2j] \cong \tR p_*\sfGa(i)_{(\bP(E),D')}, \quad
 (x_j)_{j=0}^r \mapsto \sum_{j=0}^r \  \beta^\back(\xi^j) \cup p^\back(x_j)
\end{align*}
in $\sfD(X_\Zar)$, where $\xi \in \tH^2(\bP(E),\sfGa(1)):=\tH^2(\bP(E)_\Zar,\sfGa(1)_{(\bP(E),\emptyset)})$ denotes the first Chern class of the tautological line bundle over $\bP(E)$.
\end{prop}
\begin{pf}
Note that $p$ and $\beta$ are morphisms in $\cC^\log$ by the assumption ($*_1$) in \S\ref{sect1.1} on the category $\cC$. If $D=\emptyset$, then the assertion follows from the axiom \ref{defn12-1}\,(2). One can easily deduce the general case from this case using the axioms \ref{defn12-1}\,(3), (4). The details are left to the reader as an exercise.
\end{pf}

We end this section with a logarithmic variant of Theorem \ref{thm7-1}, which will be used in the next section. We say that a morphism $f : (Y,E) \to (X,D)$ of log pairs is {\it projective}, if the underlying morphism $f : Y \to X$ is (regular) projective (cf.\ \S\ref{notation}, Definition \ref{def7-1}\,(2)) and we have $f^{-1}(D) = E$ as closed subschemes of $Y$. We define {\it the relative dimension} of a projective morphism $f : (Y,E) \to (X,D)$ as that of the underlying morphism $Y \to X$.
\begin{prop}\label{prop12-2}
Let $f : (Y,E) \to (X,D)$ be a projective morphism in $\cC^\log$. Then
\begin{enumerate}
\item[{\rm(1)}]
$f$ has a factorization
\[\xymatrix{ f : (Y,E) \;\ar@<-1pt>@{^{(}->}[r]^-g & (\bP^m_X,\bP^m_D) \ar[r]^p & (X,D)}\]
such that $g$ is a closed immersion in $\cC^\log$ and such that $p$ is the morphism induced by the projection $\bP^m_X \to X$. Furthermore, $p$ is a morphism in $\cC^\log$.
\item[{\rm(2)}]
For a decomposition of $f$ as in {\rm(}1{\rm)}, we define $f_\push$ as the composite morphism in $\sfD(X_\Zar)$
{\allowdisplaybreaks
\[\begin{CD} f_\push : \tR f_*\sfGa(i+r)_{(Y,E)}[2r]
 @= \tR p_*\tR g_*\sfGa(i+r)_{(Y,E)}[2r]  \phantom{AAAAAAii} \\
 @>{\tR p_*(g_\push)}>> \tR p_*\sfGa(i+m)_{(\bP^m_X,\bP^m_D)}[2m] \phantom{AAAAAAi} \\
 @<{\simeq}<{\text {\rm (Proposition \ref{prop12-1})}}< \bigoplus_{j=0}^m \ \sfGa(i+m-j)_{(X,D)}[2(m-j)] \\
 @>{\text {\rm projection}}>> \sfGa(i)_{(X,D)}\,, \phantom{AAAAAAAAAAAAAa} \end{CD}\]
}where $r$ denotes the relative dimension of $f$. Then $f_\push$ does not depend on the decomposition of $f$, and satisfies the properties analogous to those in Theorem \ref{thm7-1}\,{\rm(}2{\rm)}\bb{\rm(}4{\rm)}.
\end{enumerate}
\end{prop}
\begin{pf}
(1) The underlying morphism $Y \to X$ factors through a (regular) closed immersion $g$ followed by the projection $p$ as
\[ \xymatrix{ Y \;\ar@<-1pt>@{^{(}->}[r]^-g & \bP^m_X \ar[r]^p & X,} \]
both of which belongs to $\Mor(\cC)$ by the assumptions ($*_1$) and ($*_2$) on $\cC$ in \S\ref{sect1.1}. The assertions follow from this fact.
\par
(2) One can easily check the assertions by similar arguments as for \S\ref{sect7} and by endowing the projective spaces $\bP^n_X$ with the inverse image of $D \subset X$ as log poles.
\end{pf}

%

\section{Chern class for cohomology theory with log poles}\label{sect13}
The notation remains as in the previous section, and we assume the following technical but important condition on the category $\cC$ (note that the examples of $\cC$ in \S\ref{sect3} satisfy this condition):
\begin{enumerate}
\item[{\rm($*_3$)}]
{\it For a regular scheme $X \in \Ob(\cC)$ and a regular closed subscheme $Y \subset X$ which belongs to $\Ob(\cC)$, the blow-up of $X$ along $Y$ belongs to $\Mor(\cC)$.}
\end{enumerate}

Let $\sfGa(*)^\log$ be an admissible cohomology theory on $\cC^\log$. As an application of the Riemann-Roch theorem without denominators, we define Chern classes for log pairs with values in an admissible cohomology theory $\sfGa(*)^\log$ on $\cC^\log$.
%
%
Let $(X,D)$ be a log pair in $\cC$, and put $U := X \ssm D$. Let $j : U \hra X$ (resp.\ $f : D \hra X$) be the natural open (resp.\ closed) immersion. Put
\[ \tH^{2i}(X,D;\sfGa(i)):= \tH^{2i}(X_\Zar,\sfGa(i)_{(X,D)}) \;\; \hbox{ and } \;\; \tH^{2*}(X,D;\sfGa(*)):=\bigoplus_{i \ge 0} \ \tH^{2i}(X,D;\sfGa(i)). \]
We would like to construct a total Chern class map
\stepcounter{thm}
\begin{equation}\label{eq13-0}
  \sfc=(\sfc_i)_{i \ge 0} : \tK_0(U) \lra \tH^{2*}(X,D;\sfGa(*))
\end{equation}
which satisfy the following four properties:
\begin{enumerate}
\item[{\rm(L1)}]
{\it We have $\sfc_0 \equiv 1$ {\rm(}constant{\rm)}. If $\alpha \in \tK_0(U)$ is the class of a line bundle $L$ over $U$, then $\sfc_i(\alpha)=0$ for any $i \ge 2$, and $\sfc_1(\alpha)$ agrees with the value of the class of $L$ under the map}
\[ \varrho^\log : \Pic(U) \lra \tH^2(X,D;\sfGa(1)) \qquad \text{\it cf.\ Definition \ref{defn12-1}\,{\rm(}1{\rm)}}. \]
\item[{\rm(L2)}]
{\it The map $\sfc$ is contravariantly functorial for morphisms in $\cC^\log$.}
\item[{\rm(L3)}]
{\it For $\alpha,\beta \in \tK_0(U)$, we have
\[ \sfc(\alpha+\beta)=\sfc(\alpha) \cup \sfc(\beta) \quad \hbox{ in }\;\; \tH^{2*}(X,D;\sfGa(*)). \]}
\item[{\rm(L4)}]
{\it For $\alpha,\beta \in \tK_0(U)$, we have
\[ \wt{\sfc}(\alpha \otimes \beta) = \wt{\sfc}(\alpha) \,\scstar\; \wt{\sfc}(\beta) \quad \hbox{ in }\;\; \wtH^{2*}(X,D;\sfGa(*)).\; \]
Here $\wt{\sfc}(\alpha)$ denotes the augmented Chern class $(\rk(\alpha),\sfc(\alpha))$, and $\wtH^{2*}(X,D;\sfGa(*))$ denotes the $\lam$-ring defined in a similar way as for $\wtH^{2*}(\Xs,\sfGa(*))$ in Theorem \ref{thm4-1}\,{\rm(}4{\rm)}.}
\end{enumerate}
For $i \ge 1$, we define the $i$-th Chern class map $\sfc_i=\sfc_{i,(X,D)}$ by induction on the number $r$ of irreducible components on $D$. Assume first that $r=1$, i.e., $D$ is a regular integral closed subscheme of $X$, and consider the following diagram with exact rows:
\[\xymatrix{ \tK_0(D) \ar[r]^-{f_*} & \tK_0(X) \ar[r]^-{j^\back} \ar[d]_{\sfc_{i,X}} & \tK_0(U) \ar[r] & 0 \\
\tH^{2i-2}(D,\sfGa(i-1)) \ar[r]^-{f_\push} & \tH^{2i}(X,\sfGa(i)) \ar[r]^-{j^\back} & \tH^{2i}(X,D;\sfGa(i)), }\]
where the middle vertical arrow is the $i$-th Chern class map of $X$, cf.\ Corollary \ref{cor4-1}. For an element $\alpha \in \tK_0(U)$, we take a lift $\beta \in \tK_0(X)$ of $\alpha$ and define the $i$-th Chern class $\sfc_i(\alpha) \in \tH^{2i}(X,D;\sfGa(i))$ as
\[ \sfc_i(\alpha) := j^\back(\sfc_{i,X}(\beta)), \]
which is independent of the choice of $\beta$. Indeed, given another lift $\beta'$ of $\alpha$, we have $\beta=\beta'+f_*(\gamma)$ for some $\gamma \in \tK_0(D)$ and
\[ \sfc_{i,X}(\beta)=\sum_{k=0}^i \ \sfc_{i-k,X}(\beta') \cup \sfc_{k,X}(f_*(\gamma)) = \sfc_{i,X}(\beta') + \sum_{k=1}^i \ \sfc_{i-k,X}(\beta') \cup f_\push\sfP_{k-1,D/X}(\gamma) \]
by Corollary \ref{cor4-1}\,(1), (3) and Corollary \ref{cor9-1} (see also the assumption ($*_3$) at the beginning of this section). In particular, we have
\[ \sfc_{i,X}(\beta)=\sfc_{i,X}(\beta')+f_\push(x) \quad \hbox{ for } \;\;^{\exists} x \in \tH^{2i-2}(D,\sfGa(i-1)) \]
by the projection formula, and thus $j^\back(\sfc_{i,X}(\beta))=j^\back(\sfc_{i,X}(\beta'))$. Note that the resulting total Chern class map $\sfc=(\sfc_i)_{i \ge 0}$ satisfies the properties (L1)\bb(L4) by the surjectivity of the map $\tK_0(X) \to \tK_0(U)$ and Corollary \ref{cor4-1} for $X$.

Next suppose that $r \ge 2$, and fix an irreducible component $Y$ of $D$. Put
\[ D':= \ol{D \ssm Y} \;\; \hbox{(Zariski closure in $X$)},\qquad  E := Y \times_X D' \]
for each $p \ge 1$. By the  induction hypothesis, we are given total Chern class maps
\begin{align*}
 \sfc_{(X,D')} : \tK_0(U') \lra \tH^{2*}(X,D';\sfGa(*)) \qquad &\hbox{($U':=X \ssm D'$)} \\
 \sfc_{(Y,E)} : \tK_0(V) \lra \tH^{2*}(Y,E;\sfGa(*)) \;\;\, \qquad &\hbox{($V:=Y \ssm E$)},
\end{align*}
which satisfy the properties (L1)\bb(L4). Let $j : (X,D) \hra (X,D')$ (resp.\ $f : (Y,E) \hra (X,D')$) be the natural open (resp.\ closed) immersion of simplicial log pairs, and let $k : U \hra U'$ (resp.\ $g : V \hra U'$) be the natural open (resp.\ closed) immersion of simplicial schemes. In view of the diagram with exact rows
\[\xymatrix{ \tK_0(V) \ar[r]^-{g_*} & \tK_0(U') \ar[r]^-{k^\back} \ar[d]_{\sfc_{(X,D')}} & \tK_0(U) \ar[r] & 0 \\
\tH^{2*-2}(Y,E;\sfGa(*-1)) \ar[r]^-{f_\push} & \tH^{2*}(X,D';\sfGa(*)) \ar[r]^-{j^\back} & \tH^{2*}(X,D;\sfGa(*))}\]
and the construction in the case $r=1$, it is enough to check the following Riemann-Roch theorem for our induction step:
\addtocounter{thm}{-1}
\begin{thm}\label{thm13-1}
The following diagram commutes for each $i \ge 1${\rm:}
\[\xymatrix{ \tK_0(V) \ar[r]^-{g_*} \ar[d]_{\sfp_{i-1,(Y,E)/(X,D')}} & \tK_0(U') \ar[d]^{\sfc_{i,(X,D')}} \\
\tH^{2i-2}(Y,E;\sfGa(i-1)) \ar[r]^-{f_\push} & \tH^{2i}(X,D';\sfGa(i)),}\]
where $\sfp_{i-1,(Y,E)/(X,D')}$ denotes the polynomial in Chern class maps
\[ \sfp_{i-1,1}(\rk,\sfc_{1,(Y,E)},\sfc_{2,(Y,E)},\dotsc,\sfc_{i-1,(Y,E)};\,\sfc_{1,(Y,E)}(N_{V/U'})), \]
and $\sfp_{i-1,1}$ denotes the universal polynomial constructed in \S\ref{sect8}.
\end{thm}
We first note the following lemma, which is an immediate consequence of the assumption ($*_3$) at the beginning of this section:
\begin{lem}\label{lem13-1}
Let $\varrho : \sfM \to X \times \bP^1$ be the blow-up of $X \times \bP^1:=X \times_{\Spec(\bZ)} \bP^1_{\bZ}$ along $Y \times \{\infty\}$, and put ${\sf D}:=\varrho^{-1}(D' \times \bP^1)$. Then $(\sfM,{\sf D})$ is a log pair in $\cC$, and the morphism $\varrho : (\sfM,{\sf D}) \to (X \times \bP^1,D' \times \bP^1)$ belongs to $\Mor(\cC^\log)$.
\end{lem}
\begin{pf*}{\it Proof of Theorem \ref{thm13-1}}
In view of Lemma \ref{lem13-1} and the push-forward map in Proposition \ref{prop12-2}, we see that the argument in Step 2 of the proof of Theorem \ref{thm9-1} works in this situation and reduces Theorem \ref{thm13-1} to the following case:
\begin{itemize}
\item {\it $f$ is isomorphic to the zero section of the projective completion $\pi : \bP(E \oplus \bs{1}) \to Y$ of a vector bundle $E \to Y$ of rank $r$, and we have $D'=\pi^{-1}(E)$. Here $\bs{1}$ denotes the trivial line bundle over $Y$.}
\end{itemize}
Let $\varpi$ be the natural projection $U'=\pi^{-1}(V) \to V$. We use the following commutative diagram instead of \eqref{eq9-3}:
\begin{equation}\label{eq9-3'}
 \xymatrix{\tK_0(V) \ar[rrd]_-{g_*} \ar[rr]^-{\varpi^*} && \tK_0(U') \ar[d]^-{?\, \otimes \,g_*(1_{\tK})}\ar[rr]^-{\wt{\sfc}_{(X,D')}} && \wtH^{2*}(X,D';\sfGa(*)) \ar[d]^{? \tstar\, \wt{\sfc}_{(X,D')}(g_*(1_{\tK}))} \\
 && \tK_0(U')   \ar[rr]^-{\wt{\sfc}_{(X,D')}} && \wtH^{2*}(X,D';\sfGa(*))\,,}
\end{equation}
where $1_{\tK}$ denotes the unity of $\tK_0(V)$, and the square commutes by (L4) for $(X,D')$. The triangle commutes by the projection formula for Grothendieck groups. By \eqref{eq9-3'} and \eqref{eq9-6} for $g$, the composite map $\wt{\sfc}_{(X,D')} \circ g_*$ agrees with the composite map
\[\xymatrix{\tK_0(V) \ar[r]^-{\pi^\back} & \tK_0(U') \ar[r]^-{\wt{\sfc}_{(X,D')}} & \wtH^{2*}(X,D';\sfGa(*)) \ar[rrr]^-{\tstar\, \wt{\sfc}_{(X,D')}(\lam_{-1}(Q'^\vee))} &&& \wtH^{2*}(X,D';\sfGa(*)) \,,  }\]
where $Q'$ denotes the pullback, onto $U'$, of the universal quotient bundle $Q$ on $X=\bP(E \oplus \bs{1})$. Therefore we obtain the assertion by the same computations as in Step 1 of the proof of Theorem \ref{thm9-1}.
\end{pf*}
\par
By Theorem \ref{thm13-1}, we obtain the desired total Chern class map $\sfc=\sfc_{(X,D)}$ in \eqref{eq13-0}. Note that this map satisfies the properties (L1)\bb(L4) and does not depend on the choice of the irreducible component $Y \subset D$, by the surjectivity of the map $\tK_0(X) \to \tK_0(U)$. 

\begin{exmp}
{\rm The categories $\cC$ considered in Example \ref{exmp12-1}\,(1)\bb(3) satisfy the condition ($*_3$) at the beginning of this section. Consequently we obtain total Chern class map \eqref{eq13-0} for $\sfGa(*)^\log$ in those examples.}
\end{exmp}

\begin{rem}
One can pursue a simplicial analogue of the map \eqref{eq13-0}. However the authors do not know if such a Chern class in the universal situation provides us with a universal Chern class of the from
\[ \sfC_{i,(X,D)} : Rj_*(\OBQPs_U) \lra  \cK(\sfGa(i)_{(X,D)},2i) \qquad \hbox{{\rm(}$j : U=X \ssm D \hra X${\rm)}} \]
in the homotopy category $\cHs(X)$.
\end{rem}

\newpage
\appendix

\section{Motivic complex is an admissible cohomology theory}\label{appA}
\long\def\remind#1{\textcolor[named]{Peach}{}
\textcolor[named]{Blue}{\bf {#1}}}
\numberwithin{equation}{subsection}

Let the notation ($k$, $\cC$, $\bZ(i)$ for $i \geq 0$) be as in \S\ref{ex:motivic}.
We put $\bZ(i):=0$ for $i < 0$.
For $X \in \Ob(\cC)$, the motivic cohomology $\tH^*_{\!\!\cM}(X,\bZ(i))$ with $\bZ(i)$-coefficients is defined as
\[ \tH^*_{\!\!\cM}(X,\bZ(i)):=\tH^*(X_{\Zar},\bZ(i)). \]
We have
\begin{equation}\label{eqA-0}
 \tH^*_{\!\!\cM}(X,\bZ(i)) \cong \tH^*(X_{\Nis},\bZ(i)) \cong \Hom_{D^-(\NSwT(k))}(M(X),\bZ(i)[*])
\end{equation}
by \cite{SV} Corollaries 1.1.1 and 1.11.2, Proposition 1.8, where $\NSwT(k)$ denotes the category of Nisnevich sheaves with transfers over $k$ and $M(X)$ denotes $\rC\hspace{1pt}^\bullet(\bZ_\tr(X))$.
In this appendix, we prove
\begin{thm}\label{thmA-1}
Assume that $k$ admits the resolution of singularities in the sense of {\rm\cite{SV}} Definition {\rm0.1}.
Then the motivic complex $\bZ(*)=\{\bZ(i)\}_{i \in \bZ}$ is an admissible cohomology theory on $\cC$.
\end{thm}

It is easy to see that $\bZ(*)$ is a graded cohomology theory on $\cC$.
We define $\varrho$ of the axiom \ref{axiom1}\,(1) as the natural quasi-isomorphism $\bZ(1) \cong \cO^\times[-1]$.
One obtains the axiom {\rm\ref{axiom1}\,(2)} from \cite{SV} Theorem 4.5 and the fact that we have
\begin{equation}\label{eqA-1}\tag{A.0.2}
 \Hom_{D^-(\NSwT(k))}(M(X)(j),\bZ(i)[q]) = 0 \quad \hbox{for} \;\; ^\forall X \in \Ob(\cC),\, ^\forall j>i,\, ^\forall q\in \bZ
\end{equation}
under the assumption on the resolution of singularities.
One can check \eqref{eqA-1} as follows. By \eqref{eqA-0} above and loc.\ cit.\ Theorem 4.12 and Corollary 4.12.1, the problem is reduced to showing that
\[ \tH^q(X_{\Zar},\bZ(0)) \cong \tH^q((X \times \bP^1)_{\Zar},\bZ(0)) \quad \hbox{for}  \;\;^\forall q\in \bZ \]
which follows from the fact that the constant sheaf $\bZ (\cong \bZ(0))$ is flasque on $X_\Zar$ and $(X\times \bP^1)_\Zar$.
We check the axiom \ref{axiom1}\,(3) in what follows, which will be complete in \S\ref{sectA.2'} below.

\subsection{The category of unbounded cochain complexes}\label{sectA.0}
For an abelian category $\cA$, we will often write $\rC(\cA)$ for the category of unbounded cochain complexes of objects of $\cA$.
We will use the following facts freely in this appendix A,
where a {\it Grothendieck category} means an abelian category which satisfies AB5 and has a generator.
\begin{itemize}
\item
{\it $\sfS^\ab(\cC_\Zar)$ and $\NSwT(k)$ are Grothendieck categories.}
\item
{\it If $\cA$ is a Grothendieck category, then $\vD^\op\cA$ is again a Grothendieck category.
In particular, $\vD^\op\NSwT(k)$ is a Grothendieck category.}
\end{itemize}

When $\cA$ is a Grothendieck category, there exists a model category structure on $\rC(\cA)$ whose cofibrations are monomorphisms (cf.\ Lemma \ref{lemB+2} below) and whose weak equivalences are quasi-isomorphisms \cite{Bek} Proposition 3.13.
We will refer this model structure as {\it injective model structure},
and we will mean `fibrant with respect to the injective model structure' by {\it injectively fibrant}.
For a complex $K^\bullet \in \rC(\cA)$, an {\it injectively fibrant resolution} of $K^\bullet$ is a trivial cofibration
$K^\bullet \to I^\bullet$ with $I^\bullet$ injectively fibrant.
See Theorem \ref{thmB-1} below for the relation between injectively fibrant complexes and \tK-injective complexes.


For a Grothendieck category $\cA$, we define the derived category $\sfD(\cA)$ in the usual way (see e.g.\ \cite{Ha0} p.\ 37). It will be discussed in \S\ref{sectB-2} below that $\sfD(\cA)$ has small hom-sets (Corollary \ref{corB-1}) and isomorphic to $\cHsn(\rC(\cA))$, the homotopy category of $\rC(\cA)$ with respect to the injective model structure (Corollary \ref{corB-2}).

\subsection{A zig-zag representing the Gysin morphism}\label{sectA.1}
Let $f: Y \hra X$ be a closed immersion of pure codimension $r$ in $\cC$.
For $i \geq 0$, we fix an injectively fibrant resolution $\ep^i: \bZ(i) \to I^{i,\bullet}$ in $\rC(\NSwT(k))$.
We also fix a product map
\[ \alpha^i :  (I^{1,\bullet})^{\otimes i} \lra I^{i,\bullet} \quad \hbox{ in $\rC(\NSwT(k))$ } \]
which lifts, up to homotopy, the canonical isomorphism $\bZ(1)^{\otimes i} \cong \bZ(i)$ of complexes (cf.\ \cite{SV} Lemma 3.2\,(3)).
See Corollary \ref{corB-1} below for the existence of such a lift, and see also Proposition \ref{propB-1} below for the compatibility of the notions of homotopy in the sense of homotopical and homological algebra.
We construct here a (non-canonical but functorial) zig-zag of maps in $\rC(\NSwT(k))$
representing the Gysin morphism of Suslin-Voevodsky (\cite{SV} \S4)
\begin{equation}\label{eqA-sv}
 \xymatrix{ G_f=G_{Y}^{X} : M(X) \ar[r] & M(Y)(r)[2r]} \quad \hbox{ in } \;\; \sfD(\NSwT(k)).
\end{equation}

\begin{defn}\label{defnA-1}
{\rm
Let $p : \tX \to X$ (resp.\ $q : \tX' \to X \times \bA^1$) be the blow-up of $X$ along $Y$ (resp.\ the blow-up of $X \times \bA^1$ along $Y \times \{0\}$), and let $\pi : N_{Y/X} \to Y$ be the normal bundle of $f$.
Let $L^\taut$ be the tautological line bundle on $\bP(N_{Y/X} \oplus \bs{1}) \cong q^{-1}(Y \times \{0\})$, and
we fix a lift
\[ \Xi \in \Hom_{\rC^-(\NSwT(k))}(M(\bP(N_{Y/X} \oplus \bs{1})),I^{1,\bullet}[2]) \]
of $\xi=\sfc_1(L^\taut) \in \tH^2_{\!\!\cM}(\bP(N_{Y/X} \oplus \bs{1}),\bZ(1))$, cf.\ \eqref{eqA-0}, Corollary \ref{corB-1}.
We construct a zig-zag
\[\xymatrix{ \vG^X_{Y,\Xi} : M(X) \ar@{.>}[r] & M(Y)(r)[2r] \quad \hbox{ in } \;\; \rC(\NSwT(k)) }\]
representing the Gysin morphism \eqref{eqA-sv}, from the data $((I^{i,\bullet})_{i \geq 1}, (\alpha^i)_{i \geq 2},\Xi)$.
\begin{enumerate}
\item[(1)]
We define a map of complexes
\[ \sigma^i : M(\bP(N_{Y/X} \oplus \bs{1})) \lra M(Y)\otimes I^{i,\bullet}[2i] \]
as the composite
\begin{multline*}
\xymatrix{ M(\bP(N_{Y/X} \oplus \bs{1})) \ar[r]^-{\vD} & M(\bP(N_{Y/X} \oplus \bs{1}))^{\otimes i+1} } \\
\xymatrix{ \ar[rr]^-{M(q') \otimes \Xi^{\otimes i}}&& M(Y) \otimes (I^{1,\bullet}[2])^{\otimes i} \ar[r]^-{\id \otimes \alpha^i} & M(Y) \otimes I^{i,\bullet}[2i],}
\end{multline*}
where $q'$ denotes the structure morphism $\bP(N_{Y/X} \oplus \bs{1}) \to Y$.
We write $\tM(\bP(N_{Y/X} \oplus \bs{1}))$ for the mapping fiber of $M(q') : M(\bP(N_{Y/X} \oplus \bs{1})) \to M(Y)$.
The map $\sigma^r$ and the map
\[ \Sigma : \tM(\bP(N_{Y/X} \oplus \bs{1})) \lra \bigoplus_{i=1}^{r-1} \ M(Y)\otimes I^{i,\bullet}[2i] \]
defined as $\Sigma:=(\sigma^1,\sigma^2,\dotsc,\sigma^{r-1})$ will be important in what follows.
The map $\Sigma$ is a quasi-isomorphism by \cite{SV} Theorem 4.5.
\item[(2)]
For $\ep = 0,1$, let $i_\ep : X \times \{\ep \} \hra X \times \bA^1$ be the natural closed immersion.
We denote the closed immersion $\tX=\tX \times \{0\} \hra \tX'$ (resp.\ $X=X \times \{1\} \hra \tX'$) induced by $i_0$ (resp.\ $i_1$), again by $i_0$ (resp.\ $i_1$).
We put
\[\xymatrix{
 s := M(i_1 \circ \pr_1) : M(X \times \bA^1) \ar[r]^-{M(\pr_1)} & M(X)
\ar[r]^-{M(i_1)} & M(\tX'), }\]
which gives a section to $M(q) : M(\tX') \to M(X \times \bA^1)$ (in $\sfD(\NSwT(k))$).
The following zig-zag plays a key role, which represents $\beta^X_Y$ of \cite{SV} \S4:
\begin{multline*}
\xymatrix{
 t_+ : M(\tX) \ar[r]^-{M(i_0)} & M(\tX') \ar[r]^-{u} & \cone(s) & \ar[l]_-{v}^-\qis \tM(q^{-1}(Y\times \{0\}))} \\
\xymatrix{ \ar[r]^-{\Sigma} & \displaystyle \bigoplus_{i=1}^{r-1} \ M(Y)\otimes I^{i,\bullet}[2i].}
\end{multline*}
Here the arrow $u$ denotes the canonical map and $v$ denotes the composite map
\[ v : \tM(q^{-1}(Y\times \{0\})) \lra M(q^{-1}(Y\times \{0\})) \lra M(\tX') \os{u}\lra \cone(s). \]
We have used the assumption on the resolution of singularities to verify that $v$ is a quasi-isomorphism (loc.\ cit.\ Theorem 4.7). See (1) for the definition of $\Sigma$.
\item[(3)]
Because we cannot take a mapping fiber of $t_+$ constructed in (2), we introduce an auxiliary complex
\[ C_1:=\cone\bigg(\tM(q^{-1}(Y\times \{0\})) \os{(v,\Sigma)}\lra \cone(s) \oplus \bigoplus_{i=1}^{r-1} \ M(Y)\otimes I^{i,\bullet}[2i]\bigg) \]
and replace $t_+$ with the following composite map:
\[ \xymatrix{
 t : M(\tX) \ar[rr]^-{ u \circ M(i_0)} && \cone(s) \ar[r] & C_1. } \]
Since the canonical map $\bigoplus_{i=1}^{r-1} \ M(Y)\otimes I^{i,\bullet}[2i] \to C_1$ is a quasi-isomorphism, the mapping fiber $\cone(t)[-1]$ is quasi-isomorphic to $M(X)$ via the composite map
\[ b : \cone(t)[-1] \os{w}\lra M(\tX) \os{M(p)}\lra M(X) \]
by loc.\ cit.\ Theorem 4.8, where $w$ denotes the canonical map.
In other words, $w$ gives a section to $M(p)$ in $\sfD(\NSwT(k))$.
\item[(4)]
Finally, we obtain a zig-zag of maps of complexes
\begin{multline*}
\xymatrix{
 \vG^X_{Y,\Xi} : M(X) &  \ar[l]^-{\qis}_-b \cone(t)[-1] \ar[rr]^-{u \circ M(i_0) \circ w} &&  \cone(s) } \\
\xymatrix{ & \ar[l]_-v^-\qis \tM(q^{-1}(Y \times \{0\})) \ar[r]^-{\sigma^r} & M(Y)\otimes I^{r,\bullet}[2r]
 & \ar[l]^-\qis_-{\id \otimes \ep^r} M(Y)(r)[2r],}
\end{multline*}
which represents the Gysin morphism $G^X_Y$.
\end{enumerate}
}
\end{defn}


\begin{defn}\label{defnA-2}
{\rm
Suppose we are given a diagram in $\rC(\NSwT(k))$ (or in $\rC(\Ab)$)
\begin{equation}\label{eqA-3}\tag{A.1.2}
 \xymatrix{ A \ar@{.>}[r]^g \ar[d]_p & B \ar[d]^q \\
 A' \ar@{.>}[r]^{h} & B', }
\end{equation}
where $p$ and $q$ are maps of complexes and $g$ and $h$ are zig-zags of maps of complexes
\begin{align*}
g &: A=A_1 \os{g_1}\lra A_2 \os{g_2}{\us{\qis}\lla} A_3 \os{g_3}\lra \dotsb \os{g_{m-2}}{\us{\qis}\lla} A_{m-1} \os{g_{m-1}}\lra A_m=B \\
h &: A'=A'_1 \os{h_1}\lra A'_2 \os{h_2}{\us{\qis}\lla} A'_3 \os{h_3}\lra \dotsb \os{h_{n-2}}{\us{\qis}\lla} A'_{n-1} \os{h_{n-1}}\lra A'_n=B',
\end{align*}
respectively.
Then we will say that the diagram \eqref{eqA-3} {\it is commutative} or {\it commutes}, iff $m=n$ and there exist commutative squares of complexes
\[ \xymatrix{ A_i \ar[r]^-{g_i} \ar[d] & A_{i+1} \ar[d]^{\qquad\;\;\; \hbox{or}} & & A_i  \ar[d] & \ar[l]^-{\qis}_-{g_i} A_{i+1} \ar[d] \\
 A'_i \ar[r]^-{h_i} & A'_{i+1} && A'_i  & \ar[l]^-{\qis}_-{h_i} A'_{i+1}} \]
for all $i=1,2,\dotsc,m-1$.
}
\end{defn}

\begin{lem}[cf.\ \cite{SV} Lemma 4.9\,(1), (2)]\label{lemA-1}
Let $f : Y \to X$ and $\Xi$ be as in Definition \ref{defnA-1}.
\begin{enumerate}
\item[{\rm(1)}]
Let $f' : Y' \to X'$ be another closed immersion of pure codimension $r$ in $\cC$, and
suppose we are given a cartesian square
\[ \xymatrix{ Y' \, \ar@{}[rd]|\square  \ar@<-1pt>@{^{(}->}[r]^{f'} \ar@<-1pt>[d]_{p'} &  X' \ar@<-1pt>[d]^{p} \\
  Y \, \ar@<-1pt>@{^{(}->}[r]^{f} & X,} \]
where $p$ is an arbitrary morphism in $\cC$.
Then the following diagram commutes in $\rC(\NSwT(k))${\rm:}
\[\xymatrix{ M(X') \ar@{.>}[rr]^-{\vG^{X'}_{Y',p^*\Xi}} \ar[d]_{M(p)} && M(Y')(r)[2r] \ar[d]^{M(p') \otimes \id} \\
  M(X) \ar@{.>}[rr]^-{\vG^X_{Y,\Xi}} && M(Y)(r)[2r], }\]
where $p^*\Xi$ denotes the pull-back of $\Xi$ by $\bP(N_{Y'/X'}\oplus \bs{1}) \to \bP(N_{Y/X}\oplus \bs{1})$.
\item[{\rm(2)}]
For $Z \in \Ob(\cC)$, the diagram
\[\xymatrix{
M(Z \times X) \ar@{=}[d] \ar@{.>}[rr]^-{\vG^{Z \times X}_{Z\times Y,\pr_2^* \Xi}} &&
 M(Z \times Y)(r)[2r]\ar@<-10pt>@{=}[d] \\
M(Z) \otimes M(X) \ar@{.>}[rr]^-{\id \otimes \vG^X_{Y,\Xi}}  &&  M(Z) \otimes M(Y)(r)[2r]}\]
commutes in $\rC(\NSwT(k))$, where $\pr_2$ denotes the second projection $Z \times X \to X$.
\end{enumerate}
\end{lem}
\begin{pf}
(1)\, Since $N_{Y'/X'} \cong N_{Y/X} \times_Y Y'$,
 the pull-back of the tautological line bundle on $\bP(N_{Y/X}\oplus \bs{1})$ is isomorphic to that on $\bP(N_{Y'/X'}\oplus \bs{1})$,
 which implies that the zig-zag $\vG^{X'}_{Y',p^* \Xi}$ using the pull-back $p^*\Xi$ makes sense. The assertion follows from this fact and the construction of the zig-zags.
 \par
The assertion (2) follows from similar arguments as for (1), whose details are left to the reader.
\end{pf}

The following lemma is a slight generalization of Lemma \ref{lemA-1}\,(1), whose proof is also straight-forward and left to the reader:
\begin{lem}\label{lemA-2}
Under the same setting as in Lemma \ref{lemA-1}\,{\rm(}1{\rm)},
we fix another injectively fibrant resolution $\bZ(i) \to J^{i,\bullet}$ in $\rC(\NSwT(k))$ and a product map
\[ \beta^i :  (J^{1,\bullet})^{\otimes i} \lra J^{i,\bullet} \quad \hbox{ in \;\; $\rC(\NSwT(k))$ } \]
lifting, up to homotopy, the canonical isomorphism $\bZ(1)^{\otimes i} \cong \bZ(i)$, for each $i \geq 0$.
Let
\[ \Xi' : M(\bP(N_{Y'/X'} \oplus \bs{1})) \lra J^{1,\bullet}[2]\]
be a map of complexes representing the first Chern class of the tautological line bundle on $\bP(N_{Y'/X'} \oplus \bs{1})$,
and suppose that we are given maps $\gamma^i : I^{i,\bullet} \to J^{i,\bullet}$ of injectively fibrant resolutions of $\bZ(i)$ for $i \geq 1$, fitting into the following commutative diagrams in $\rC(\NSwT(k))${\rm:}
\[
\xymatrix{ M(\bP(N_{Y'/X'} \oplus \bs{1})) \ar[r]^-{\Xi'} \ar[d] & I^{1,\bullet} \ar[d]^{\gamma^1} \\
M(\bP(N_{Y/X} \oplus \bs{1})) \ar[r]^-{\Xi} & J^{1,\bullet},} \qquad\qquad
\xymatrix{(I^{1,\bullet})^{\otimes i} \ar[r]^-{\alpha^i} \ar[d]_{(\gamma^1)^{\otimes i}} & I^{i,\bullet} \ar[d]^{\gamma^i} \\
(J^{1,\bullet})^{\otimes i} \ar[r]^-{\beta^i} & J^{i,\bullet}.}\]
Then the diagram
\[\xymatrix{ M(X') \ar@{.>}[rr]^-{\vG^{X'}_{Y',\Xi'}} \ar[d]_{M(p)} && M(Y')(r)[2r] \ar[d]^{M(p')\otimes \id} \\
  M(X) \ar@{.>}[rr]^-{\vG^X_{Y,\Xi}} && M(Y)(r)[2r]}\]
commutes in $\rC(\NSwT(k))$, where $\vG^{X'}_{Y',\Xi'}$ is defined by the data $((J^{i,\bullet})_{i \geq 1}, (\beta^i)_{i \geq 2},\Xi')$.
\end{lem}

\subsection{Lifts of the first Chern class of line bundles}\label{sectA.1'}
The aim of this subsection is to construct a morphism \eqref{eqA+6} below from the first Chern class of a line bundle $\Ls$ over a simplicial variety $\Zs \in \vD^\op\cC$, which will play an important role in our construction of Gysin morphisms for simplicial varieties.
We work with the category $\vD^\op\NSwT(k)$, the category of simplicial objects in $\NSwT(k)$. Let
\[ \tau : \NSwT(k) \lra \vD^\op\NSwT(k) \]
be the functor which assigns $A \in \Ob(\NSwT(k))$ to the constant simplicial object $A$.

For an object $\As \in \Ob(\vD^\op\NSwT(k))$,
let $\rC\hspace{1pt}^\bullet(\As)$ be the Suslin complex of $\As$, i.e., the term $\rC\hspace{1pt}^{-q}(\As)$ in degree $\bullet=-q$ is the simplicial Nisnevich sheaf with transfers
\[ U \in \Ob(\cC) \longmapsto \As(U \times \bA^q). \]
This construction naturally extends to complexes $\As^\bullet \in \Ob(\rC(\vD^\op\NSwT(k)))$.
\begin{prop}[\cite{SV} Corollary 1.11.2]\label{lemA-3}
Let $\As^\bullet$ and $\Bs^\bullet$ be objects of $\rC\hspace{1pt}^-(\vD^\op\NSwT(k))$,
and assume that $A_p^\bullet$ has homotopy invariant cohomology sheaves for each $p \geq 0$.
Then the natural embedding $\Bs^\bullet \to \rC\hspace{1pt}^\bullet(\Bs^\bullet)$ induces an isomorphism
\[ \Hom_{\sfD(\vD^\op\NSwT(k))}(\Bs^\bullet,\As^\bullet) \cong \Hom_{\sfD(\vD^\op\NSwT(k))}(\rC\hspace{1pt}^\bullet(\Bs^\bullet),\As^\bullet). \]
\end{prop}
\begin{pf}
One obtains the assertion by the same arguments as in \cite{SV} Corollary 1.10.2--Corollary 1.11.2,
starting from loc.\ cit.\ Corollary 1.10.2 for usual Nisnevich sheaves with transfers.
\end{pf}

\begin{prop}\label{lemA+3}
For $\Zs \in \Ob(\vD^\op \cC)$, there exists a natural equivalence of functors
\[ \vG(Z_\star,-) \cong \Hom_{\vD^\op\NSwT(k)}(\bZ_\tr(\Zs),\tau(-)) : \NSwT(k) \lra \Ab, \]
where $\bZ_\tr(\Zs)$ denotes the simplicial Nisnevich sheaves with transfers consisting of the data $((\bZ_\tr(Z_p))_{[p] \in \Ob(\vD)},(\bZ_\tr(a^Z))_{a \in \Mor(\vD)})$.
\end{prop}
\begin{pf}
For $0 \leq i \leq p$, let $d_i: [0] \to [p]$ be the map sending $0$ to $i$.
Let $\cF \in \Ob(\NSwT(k))$ be an arbitrary sheaf.
By the definition of $\vG(Z_\star,\cF)$ and Yoneda's lemma, $\vG(Z_\star,\cF)$ agrees with the kernel of the map
\[ \bZ_\tr(d_0^Z)^*-\bZ_\tr(d_1^Z)^* :
 \Hom_{\NSwT(k)}(\bZ_\tr(Z_0),\cF) \lra \Hom_{\NSwT(k)}(\bZ_\tr(Z_1),\cF), \]
which show that
\[ \vG(Z_\star,\cF) \supset \Hom_{\vD^\op\NSwT(k)}(\bZ_\tr(\Zs),\tau(\cF)). \]
To show the inclusion in the other direction, it is enough to show that for any $s \in \vG(Z_\star,\cF)$ and $p \geq 1$, the value of $s$ under the map
\[  \bZ_\tr(d_i^Z)^* :
 \Hom_{\NSwT(k)}(\bZ_\tr(Z_0),\cF) \lra \Hom_{\NSwT(k)}(\bZ_\tr(Z_p),\cF) \]
is independent of $0 \leq i \leq p$. 
One can easily check this claim by induction on $p$ and simplicial identities, whose details are left to the reader.
\end{pf}

Let $\Zs$ be a simplicial object in $\cC$, and let $\Ls$ be a line bundle on $\Zs$.
Let $\ep : \bZ(1) \to I^\bullet$ be an injectively fibrant resolution of $\bZ(1)$ in $\rC(\NSwT(k))$.
We consider a composite map
{\allowdisplaybreaks
\begin{align}
\tH^2(Z_{\star,\Nis},\bZ(1)) & \; = \;
  \frac{\Ker(d: \vG(Z_\star,I^2)\to\vG(Z_\star,I^3))}{\Image(d: \vG(Z_\star,I^1)\to\vG(Z_\star,I^2))} \notag \\
& \; \cong \; \Hom_{\tK(\vD^\op\NSwT(k))}(\bZ_\tr(\Zs),\tau(I^\bullet)[2]) \notag \\
& \hspace{.8pt}\to \; \Hom_{\sfD(\vD^\op\NSwT(k))}(\bZ_\tr(\Zs),\tau(I^\bullet)[2]) \notag \\
& \; \cong \; \Hom_{\sfD(\vD^\op\NSwT(k))}(M(\Zs),\tau(\bZ(1))[2]), \label{eqA-7}
\end{align}
}where $\tK(\vD^\op\NSwT(k))$ denotes the homotopy category of $\rC(\vD^\op\NSwT(k))$ (cf.\ \S\ref{sectB-2} below), and the second isomorphism is obtained from Proposition \ref{lemA+3}; the last isomorphism follows from Proposition \ref{lemA-3} for $\Bs=\bZ_\tr(\Zs)$.
Consequently, the first Chern class $\sfc_1(\Ls)$ is defined in the last group of \eqref{eqA-7}.
Now we fix an injectively fibrant resolution $\ep_\star^1 : \tau(\bZ(1)) \to \Is^{1,\bullet}$ in $\rC(\vD^\op\NSwT(k))$, 
and a morphism
\begin{equation}\label{eqA+6}\tag{A.3.2}
 \zeta_\star : M(\Zs) \lra \Is^{1,\bullet}[2] \quad \hbox{in} \;\; \rC(\vD^\op\NSwT(k))
\end{equation}
which lifts $\sfc_1(\Ls)$, cf.\ Corollary \ref{corB-1}.
For a morphism $a : [p] \to [q]$ in $\vD$, there is a commutative diagram
\begin{equation}\label{eqA-6}\tag{A.3.3}
\xymatrix{ M(Z_q) \ar[d]_{M(a^Z)} \ar[r]^{\zeta_q} & I_q^{1,\bullet}[2] \ar[d]_{a^{I^{1,\bullet}}} & \ar[l]_{\ep_q^1[2]}^{\qis} \bZ(1) [2] \ar@{=}[d] \\
 M(Z_p) \ar[r]^{\zeta_p} & I_p^{1,\bullet}[2] &  \ar[l]_{\ep_p^1[2]}^{\qis}  \bZ(1)[2] }
\end{equation}
in $\rC(\NSwT(k))$.
We note here that the arrows $\ep_p^1$ and $\ep_q^1$ are fibrant resolutions by Lemma \ref{lemA-4} below, and that the zig-zag in the lower row (resp.\ the upper row) represents the first Chern class $\sfc_1(L_p) \in  \tH^2(Z_{p,\Nis},\bZ(1))$ (resp.\ $\sfc_1(L_q) \in  \tH^2(Z_{q,\Nis},\bZ(1))$), cf.\ \eqref{eqA-0}.
In the following lemma, an {\it injective fibration} means a fibration with respect to the injective model structure:

\begin{lem}\label{lemA-4}
Let $g_\star : \As^\bullet \to \Bs^\bullet$ be an injective fibration in $\rC(\vD^\op\NSwT(k))$.
Then $g_p : A_p^\bullet \to B_p^\bullet$ is an injective fibration in $\rC(\NSwT(k))$ for each $p \geq 0$.
\end{lem}
\begin{pf}
We endow the category $\vD^\op\rC(\NSwT(k))=\rC(\vD^\op\NSwT(k))$ with another model structure, the Reedy model structure (cf.\ \cite{GJ} Chapter \VII) associated with the injective model structure on $\rC(\NSwT(k))$,
and consider the following diagram of model categories (in fact, a Quillen equivalence):
\[\xymatrix{ (\vD^\op\rC(\NSwT(k)),\text{Reedy structure}) \ar@<2pt>[r]^-{\id} & \ar@<2pt>[l]^-{\id} (\rC(\vD^\op\NSwT(k)),\text{injective structure}) }\]
Both hand sides share the class of weak equivalences (i.e., level-wise quasi-isomorphisms), and
 the upper arrow sends Reedy cofibrations to injective cofibrations (i.e., injective maps), cf.\ \cite{GJ} Chapter \VII, Corollary 2.6\,(2).
Hence the lower arrow sends injective fibrations to Reedy fibrations,
 and the assertion follows from the fact that Reedy fibrations give injective fibrations at each level, loc.\ cit.\ VII Corollary 2.6\,(1).
\end{pf}

\subsection{Construction of push-forward maps}\label{sectA.2}
Let $\tau : \NSwT(k) \to \vD^\op\NSwT(k)$ be as in \S\ref{sectA.1'}.
For each $i \geq 0$, we fix an injectively fibrant resolution $\ep_\star^i : \tau(\bZ(i)) \to \Is^{i,\bullet}$ in $\rC(\vD^\op\NSwT(k))$.
For $A^\bullet \in \rC(\NSwT(k))$ and $p \geq 0$,
put
\[ \tR\Hom(A^\bullet,\bZ(i))_p:= \Hom_{\rC(\NSwT(k))}^\bullet(A^\bullet,I_p^{i,\bullet})
 \quad \hbox{(hom-complex)}. \]
Each $\ep_p^i : \bZ(i) \to I_p^{i,\bullet}$ is an injectively fibrant resolution in $\rC(\NSwT(k))$ by Lemma \ref{lemA-4}.
\par
Let $\nu_* : \NSwT(k) \to \sfS^\ab(\cC_\Zar)$ be the natural restriction (forgetful) functor.
For each $i \geq 0$, we fix an injectively fibrant resolution $\kappa^i : \bZ(i) \to J^{i,\bullet}$ of $\bZ(i)=\nu_*\bZ(i)$ in $\rC(\sfS^\ab(\cC_\Zar))$.
By the left hand isomorphism of \eqref{eqA-0} for all $X \in \Ob(\cC)$,
 the map $\bZ(i) \to \nu_*I_0^{i,\bullet}$ is a trivial cofibration in $\rC(\sfS^\ab(\cC_\Zar))$, and there exists a quasi-isomorphism of complexes
$\lam^i : \nu_*I_0^{i,\bullet} \to J^{i,\bullet}$
extending $\kappa^i : \bZ(i) \to J^{i,\bullet}$.
Under these settings,
one can rephrase the isomorphisms in \eqref{eqA-0} as
\begin{align}
\tR\vG(X,\bZ(i)) & \us{\qis}{\os{\,\lam^i}{\lla}} \tR\vG(X_\Nis,\bZ(i))_0 \notag \\
 & \us{\qis}{\os{c^{I^{i,\bullet}}}\lra} \tR\vG(X_\Nis,\bZ(i))_p  \us{\qis}\lla \tR\Hom(M(X),\bZ(i))_p
\label{eqA-01}
\end{align}
for arbitrary $p \geq 0$ and $X \in \Ob(\cC)$, where $\tR\vG(X,\bZ(i))$ (resp.\ $\tR\vG(X_\Nis,\bZ(i))_p$) denotes $\vG(X,J^{i,\bullet})$ (resp.\ $\vG(X,I_p^{i,\bullet})$) and $c$ denotes the canonical map $[p] \to [0]$ in $\vD$.
\par
For $i \geq 2$, we fix a map in $\rC(\vD^\op\NSwT(k))$
\[ \alpha_\star^i : (\Is^{1,\bullet})^{\otimes i} \lra \Is^{i,\bullet} \]
which represents the canonical isomorphism $\tau(\bZ(1))^{\otimes i} \cong \tau(\bZ(i))$.
For $i,j \geq 0$, we also fix a map in $\rC(\vD^\op\NSwT(k))$
\[ \alpha_\star^{i,j} : \Is^{i,\bullet} \otimes \Is^{j,\bullet} \lra \Is^{i+j,\bullet} \]
which represents the canonical isomorphism $\tau(\bZ(i)) \otimes \tau(\bZ(j)) \cong \tau(\bZ(i+j))$.
(We do not suppose that $\alpha_\star^{1,1}=\alpha_\star^2$.)
For $i,j,r \geq 0$, there is a square which is commutative up to homotopy in $\rC(\vD^\op\NSwT(k))$
\begin{equation}\label{eqA-5}
\xymatrix{ \Is^{i,\bullet} \otimes \Is^{j,\bullet} \otimes \Is^{r,\bullet}
 \ar[rr]^-{\alpha^{i,j} \otimes \id} \ar[d]_{\id \otimes \alpha^{j,r}} && \Is^{i+j,\bullet} \otimes \Is^{r,\bullet}
 \ar[d]^{\alpha^{i+j,r}} \\
 \Is^{i,\bullet} \otimes \Is^{j+r,\bullet} \ar[rr]^-{\alpha^{i,j+r}} && \Is^{i+j+r,\bullet}.}
\end{equation}
Let $f : Y_\star \hra X_\star$ be a strict closed immersion of pure codimension $r$ of simplicial schemes in $\cC$.
Let $N_\star$ be the normal bundle of $f : \Ys \hra \Xs$, and let 
 \[ \Xi_\star \in  \Hom_{\rC^-(\vD^\op\NSwT(k))}(M(\bP(\Ns \oplus \bs{1})),\Is^{1,\bullet}[2]) \]
be a lift of $\sfc_1(\Ls^\taut)$, the first Chern class of the tautological line bundle on $\bP(\Ns \oplus \bs{1})$, cf.\ \S\ref{sectA.1'}.
For $U \in \Ob((X_\star)_\Zar)$, we put $Y_U := Y_p \times_{X_p} U$, where $p \ge 0$ denotes the integer such that $U$ is an open subset of $X_p$, and consider a zig-zag in $\rC(\Ab)$
\begin{multline*}
\!\!\!\!\xymatrix{
 (f_U)_{\push,\Xi_p} : \tR\vG(Y_{U},\bZ(i)) \ar@{.>}[r]^-{\eqref{eqA-01}} &
\tR\Hom(M(Y_U),\bZ(i))_p
 \, \ar[r]^-{\delta_p^{i,r}} & \, \tR\Hom(M(Y_U)(r),\bZ(i+r))_p} \\
 \xymatrix{\ar@{.>}[rr]^-{-(\vG^U_{Y_U,\Xi_p})^*} && \tR \Hom(M(U)[-2r],\bZ(i+r))_p \ar@{.>}[r]^-{\eqref{eqA-01}} &\tR\vG(U,\bZ(i+r))[2r].}
\end{multline*}
Here $\vG^U_{Y_U,\Xi_p}$ denotes the zig-zag defined by the data $((I_p^{i,\bullet})_{i \geq 1}, (\alpha_p^i)_{i \geq 2},\Xi_p)$ (cf.\ \S\ref{sectA.1}), and $\delta_p^{i,r}$ denotes the composite map
\begin{multline*}
 \xymatrix{\delta_p^{i,r} : \tR\Hom(M(Y_U),\bZ(i))_p \ar[rr]^-{g \mapsto g \otimes \ep^r} &&
 \Hom_{\rC(\NSwT(k))}^\bullet(M(Y_U)(r),I_p^{i,\bullet}\otimes I_p^{r,\bullet}) } \\
 \xymatrix{ \ar[r]^-{\alpha_p^{i,r}} &
 \Hom_{\rC(\NSwT(k))}^\bullet(M(Y_U)(r), I_p^{i+r,\bullet}) = \tR\Hom(M(Y_U)(r),\bZ(i+r))_p.}
\end{multline*}
The reason of the sign of $-(\vG^U_{Y_U,\Xi_p})^*$ will be explained in Lemmas \ref{propA-0} and \ref{lem3-1}\,(3) below.
The zig-zag $(f_U)_{\push,\Xi_p}$ is contravariantly functorial in $U \in \Ob((X_\star)_\Zar)$ by \eqref{eqA-6} and  Lemma \ref{lemA-2}, and hence, yields a morphism
\[ f_{\push,\Xi_\star} : f_*\bZ(i)_{Y_\star} \lra \bZ(i+r)_{X_\star}[2r] \quad \hbox{ in } \;\; \sfD((\Xs)_\Zar).\]
One can easily check that $f_{\push,\Xi_\star}$ is independent of the choice of $\Xi_\star$ lifting $\sfc_1(\Ls^\taut)$.
For this reason, we will write $f_\push$ for $f_{\push,\Xi_\star}$, in what follows.
We will often write $f_\push$ for the following morphism induced by the above $f_\push$:
\[  f_\push : \bZ(i)_{Y_\star} \lra Rf^!\bZ(i+r)_{X_\star}[2r] \quad \hbox{ in } \;\; \sfD((\Ys)_\Zar).\]
\begin{lem}\label{propA-0}
$f_\push$ satisfies {\rm(}3a{\rm)} and {\rm(}3d{\rm)} of Definition \rm\ref{axiom1}.
\end{lem}
\def\tXps{\tX_{\hspace{-.5pt}\star}\hspace{-3.5pt}{}'}
\def\tXs{\tX_{\hspace{-.5pt}\star}}
\begin{pf}
(3d) immediately follows from Lemma \ref{lemA-1}\,(1).
To check (3a), suppose that $f$ has pure codimension $1$, i.e., $\Ys$ is an effective Cartier divisor on $\Xs$ via $f$.
Let $\tXps$ be the blow-up of $\Xs \times \bA^1$ along $\Ys \times \{0\}$, and let $\Es$ be the exceptional divisor on $\tXps$.
Let $\Ls^\taut$ be the tautological line bundle on $\Es \cong \bP(\Ns \oplus \bs{1})$.
Note that the strict closed immersion $i_0 : \Xs=\Xs \times \{0\} \hra \Xs \times \bA^1$ induces a strict closed immersion $\Xs \hra \tXps$, which we denote again by $i_0$.
Let $L_{\Es}$ be the line bundle on $\tXps$ associated with the effective Cartier divisor $\Es$, whose sheaf of sections gives the invertible sheaf $\cO(\Es)$.
The pull-back of $L_{\Es}$ by $i_0 : \Xs \hra \tXps$ is isomorphic to the line bundle $L_{\Ys}$ on $\Xs$ associated with $\Ys$, which implies that the map
\[\xymatrix{ i_0^\back : \tH^{2}_{\Es}(\tXps,\bZ(1)) \ar[r] & \tH^{2}_{\Ys}(\Xs,\bZ(1))} \]
sends $\sfc_1(L_{\Es})$ to $\sfc_1(L_{\Ys})$. On the other hand, the dual $(L_{\Es})^\vee$ restricts to $\Ls^\taut$ under $\Es \hra \tXps$. By this fact and our construction of $f_!$, we have $i_0^\back(\sfc_1(L_{\Es}))=f_\push(1)$ in $\tH^{2}_{\Ys}(\Xs,\bZ(1))$. We obtain (3a) from these facts.
\end{pf}

\begin{prop}\label{propA-1}
$f_\push$ satisfies the projection formula {\rm(}3b{\rm)} of Definition \ref{axiom1}.
\end{prop}
\begin{pf}
For $i,j \geq 0$ and $U \in \Ob((\Xs)_\Zar)$ with $U \subset X_p$, we put
\[ \tR^{i,\bullet}(U):=\tR\Hom(M(U),\bZ(i))_p \quad \hbox{and}\quad \tR^{i,\bullet}(Y_U(j)):=\tR\Hom(M(Y_U)(j),\bZ(i))_p \]
for simplicity. We consider the following diagram in $\rC(\Ab)$:
\begin{equation}\label{eqA-4}\tag{A.4.3}
\xymatrix{
\tR^{i,\bullet}(U)  \os{\bL}\otimes \tR^{j,\bullet}(Y_U)
 \ar[rr]^-{\alpha_p^{i,j}\,\circ(f_U^\back \otimes \id)} \ar[d]_{\id \otimes \delta_p^{j,r}}
  && \tR^{i+j,\bullet}(Y_U) \ar[d]^{\delta_p^{i+j,r}} \\
\tR^{i,\bullet}(U)  \os{\bL}\otimes \tR^{j+r,\bullet}(Y_U(r))
 \ar[rr]^-{\alpha_p^{i,j+r}\,\circ(f_U^\back \otimes \id)} \ar@{.>}[d]_{\id \otimes (-\vG^U_{Y_U,\Xi_p})^*}
  && \tR^{i+j+r,\bullet}(Y_U(r)) \ar@{.>}[d]^{(-\vG^U_{Y_U,\Xi_p})^*} \\
 \tR^{i,\bullet}(U) \os{\bL}\otimes \tR^{j+r,\bullet}(U) [2r] \ar[rr]^-{\alpha_p^{i,j+r}}
 && \tR^{i+j+r,\bullet}(U)[2r].}
\end{equation}
whose upper square commutes up to a homotopy defined globally on $(\Xs)_\Zar$ by \eqref{eqA-5}. 
The lower square commutes in the sense of Definition \ref{defnA-2} by the commutative squares
{\small\[\xymatrix{
M(U) \otimes M(Y_U)(r)[2r] \ar@{=}[r] & M(U \times Y_U)(r)[2r]
 &&& \ar[lll]_-{M(f_U \times \id)(r)[2r]} M(Y_U)(r)[2r] \\
M(U) \otimes M(U) \ar@{=}[r] \ar@{.>}[u]^-{\id \otimes (-\vG^U_{Y_U,\Xi_p})} 
 & M(U \times U) \ar@{.>}[u]_-{-\vG^{U \times U}_{U\times Y_U,\Xi_p}}
 &&& \ar[lll]_-{M(\vD)} \ar@{.>}[u]_-{-\vG^U_{Y_U,\Xi_p}}  M(U)  }\]
}in $\rC(\NSwT(k))$, cf.\ Lemma \ref{lemA-1}\,(1), (2).
One can easily derive the assertion from the commutativity (up to a functorial homotopy) of the above diagram.
\end{pf}

\subsection{Proof of (3c)}\label{sectA.2'}
\def\Es{E}
\def\Xs{X}
\def\Ys{Y}
\def\Zs{Z}
\def\tXs{\tX}
\def\tXps{\tX{}'}

Let $f : \Ys \hra \Xs$ be a closed immersion of pure codimension $r$ in $\cC$,
the above construction of $f_\push$ works in the Nisnevich topology as well, and gives rise to a morphism
\[ f_\push : f_*\bZ(i)_{\Ys} \lra \bZ(i+r)_{\Xs}[2r] \]
both in $\sfD(\Xs_\Zar)$ and $\sfD(\Xs_\Nis)$.
Our task is to check (3c) for $f_\push$ in the Zariski topology, which is reduced to showing the same properties in the Nisnevich topology by \cite{SV} Corollary 1.1.1.
\begin{lem}\label{lem3-1}
Let $\tau$ be either $\Zar$ or $\Nis$, and put $\cl_{\Xs}(\Ys):=f_!(1) \in \tH^{2r}_{\Ys}(\Xs_\tau,\bZ(r))$. Then the following holds{\rm:}
\begin{enumerate}
\item[{\rm(1)}]
Let $p : \tXs \to \Xs$ be the blow-up along $\Ys$, and let $\Es \subset \tXs$ be the exceptional fiber.
Let $\xi  \in \tH^2_{\Es}(\tXs_\tau,\bZ(1))$ be the localized first Chern class of the invertible sheaf $\cO(1)=\cO(-E)$ on $\tX$.
Then we have
\[ \xi^r \equiv  -p^\back \cl_{\Xs}(\Ys) \; \mod \; \bigoplus_{i=1}^{r-1}\, \tH^{2i}(\Ys_\tau,\bZ(i))\]
under the isomorphism
\[ \tH^{2r}_{\Es}(\tXs_\tau,\bZ(r)) \cong \tH^{2r}_{\Ys}(\Xs_\tau,\bZ(r)) \oplus \bigoplus_{i=1}^{r-1} \ \tH^{2i}(\Ys_\tau,\bZ(i)). \]
For this isomorphism, see {\rm \cite{SV}} Theorem {\rm 4.8} and the construction after loc.\ cit.\ Theorem {\rm 4.8}.
\item[{\rm(2)}]
If $f$ has a retraction $\pi : X \to Y$, then $f_\push : f_*\bZ(i)_{Y} \to \bZ(i+r)_{X}[2r]$ factors as
\[ \xymatrix{
f_*\bZ(i)_{Y}=f_*f^*\pi^*\bZ(i)_{Y} \ar[r]^-{f_*(\pi^\back)} & f_*f^*\bZ(i)_{X} \ar[rr]^-{\cup \, \cl_X(Y)} && \bZ(i+r)_{X}[2r]
}\quad \hbox{ in } \;\; \sfD(X_\tau). \]
\item[{\rm(3)}]
If there exists a simple normal crossing divisor $D=\bigcup_{j=1}^r \, D_j$ on $X$ with\, $\bigcap_{j=1}^r \, D_j = Y$ and each $D_j$ smooth, then we have
\[ \cl_X(Y) = \cl_X(D_1) \cup \cl_X(D_2) \cup \dotsb \cup \cl_X(D_r) \quad \hbox{ in } \;\;
 \tH^{2r}_Y(X_\tau,\bZ(r)). \]
\end{enumerate}
\end{lem}
We first prove (3c) in $\sfD(X_\Nis)$ admitting this lemma.
By an argument of Voevodsky in \cite{V1} Proof of Theorem 4.14,
for a point $x \in X_p$ contained in $Z$, there exist a Nisnevich neighborhood $x \to U \os{w}\to X$
and an \'etale map $v : U \to Z \times \bA^{r+s}$ such that
 $w^{-1}(Y)=v^{-1}(Y \times 0_r)$,
 $w^{-1}(Z)=v^{-1}(Z \times (0_s,0_r))$,
 where $0_r$ (resp.\ $0_s$) denotes the origin of $\bA^r=\bA^r_k$ (resp.\ $\bA^s$).
Hence one can reduce (3c) to the following sufficiently local situation (L) by taking a Nisnevich hypercovering:
\begin{itemize}
\item[(L)]
{\it We have $\Xs \cong \Zs \times \bA^{r+s}$, and the closed immersions $f$, $g$ and $g \circ f$ are isomorphic to the zero sections
 $\Ys \to \Ys \times \bA^r$, $\Zs \to \Zs \times \bA^s$ and $\Zs \to \Zs \times \bA^{r+s}$, respectively.}
\end{itemize}
In this case, (3c) follows from Lemma \ref{lem3-1}\,(2), (3), which is parallel to the arguments in \cite{FG} Proof of Proposition 1.2.1.
\par\medskip
\begin{pf*}{Proof of Lemma \ref{lem3-1}}
(1) Let $p' : \tX' \to X \times \bA^1$ be the blow-up along $Y \times \{0\}$.
Let $E' \subset \tX'$ be the exceptional fiber, and let $\zeta  \in \tH^2_{E'}(\tX'_\tau,\bZ(1))$ be the localized first Chern class of the invertible sheaf $\cO(1)=\cO(-E')$ on $\tX'$.
Let $i_0 : \tX \hra \tX'$ be the natural closed immersion induced by the closed immersion $X \times \{0\} \hra X \times \bA^1$.
Consider the following composite map:
\[\xymatrix{
\bZ \ar[rr]^-{1 \mapsto \zeta^r} && \tH^{2r}_{E'}(\tX'_\tau,\bZ(r))
 \ar[r]^-{i_0^\back} & \tH^{2r}_{E}(\tX_\tau,\bZ(r)) \ar[r] & \tH^{2r}_Y(X_\tau,\bZ(r))},\]
where the last arrow is the retraction to $p^\sharp$.
This composite map sends $1$ to $-\cl_X(Y)$ by the construction of the Gysin morphism $G_{Y}^{X} : M(X) \to M(Y)(r)[2r]$ and our definition of $f_\push$.
Now the assertion follows from the fact that $i_0^\back(\zeta^r)=\xi^r$.\par
(2) By Proposition \ref{propA-1}, the following triangle commutes in $\sfD(X_\tau)$:
\[ \xymatrix{ f_*f^*\bZ(i)_{X} \ar[rd]^-{\cup \, \cl_X(Y)} \ar[d]_{f_*(f^\back)} \\
 f_*\bZ(i)_{Y} \ar[r]^-{f_\push} & \bZ(i+r)_{X}[2r]. }\]
One obtains the assertion immediately from this fact.
\par
(3) Since we have $\cl_X(D_j)=\sfc_1(\cO_X(D_j))$ by Lemma \ref{propA-0}, the assertion follows from (1) and the same arguments as in \cite{FG} Proposition 1.1.4.
\end{pf*}
\noindent
This completes the proof of Theorem \ref{thmA-1}.

\section{Two remarks on homotopical and homological algebra}\label{sectB}
\begin{center}
{ by\, Kei Hagihara}
\end{center}
Let $\cA$ be an abelian category.
We say that a morphism $f : K \ra L$ in $\rC(\cA)$ is a {\it monic quasi-isomorphism} if it is a monomorphism and a quasi-isomorphism.

\subsection{\tK-injectivity and right lifting property}

The following fact is stated in \S2 of an earlier version of \cite{Ho2}, where the credit is due to Avramov, Foxby and Halperin. Unfortunately their preprint (referred as [AFH97] there) has not been available so far.
\begin{thm}\label{thmB-1}
For a complex $X=X^\bullet \in \rC(\cA)$, the following two conditions are equivalent to each other{\rm:}
\begin{itemize}
\item[{\rm (i)}]
The morphism $X \to 0$ has the right lifting property with respect to all monic quasi-isomorphisms $f : K \to L$,
that is, an arbitrary morphism $g : K \to X$ in $\rC(\cA)$ is extended to a morphism $h : L \to X$ in $\rC(\cA)${\rm:}
\[ \xymatrix{K \ar[d]_f \ar[r]^g &  X \ar[d] \\ L \ar@{.>}[ru]_{\!\!h} \ar[r] & 0. \hspace{-4pt}}\]
\item[{\rm (ii)}]
$X$ is \tK-injective in the sense of Spaltenstein {\rm \cite{Spa}}, and $X^q$ is injective for all $q \in \bZ$.
\end{itemize}
\end{thm}
\noindent
We give a proof of this fact in this subsection.
To prove the theorem, we prepare Lemmas \ref{lemB+1}--\ref{lemB+3} below.
We omit the proof of the first lemma:
\begin{lem}\label{lemB+1}
Let $T$ be an object of $\cA$, and let $K=K^\bullet$ be an object of $\rC(\cA)$.
Then we have
\[ \Hom_{\rC(\cA)}(\cone(\id_T)[-n],K) \cong \Hom_{\cA}(T,K^{n-1}) \]
for any $n \in \bZ$.
\end{lem}


\begin{lem}\label{lemB+2}
A morphism $f : K \to L$ in $\rC(\cA)$ is a monomorphism, if and only if $f^n : K^n \to L^n$ is a monomorphism in $\cA$ for any $n \in \bZ$.
\end{lem}
\begin{pf}
Assume first that $f$ is a monomorphism. Then for arbitrary $T \in \Ob(\cA)$ and $n \in \bZ$, the map
\[ f_* : \Hom_{\rC(\cA)}(\cone(\id_T)[-n-1],K) \lra \Hom_{\rC(\cA)}(\cone(\id_T)[-n-1],L) \]
is injective, which together with Lemma \ref{lemB+1} implies that $f^n : K^n \to L^n$ is a monomorphism.
The converse immediately follows from the fact that 
\[ \Hom_{\rC(\cA)}(T,K) \subset \prod_{n \in \bZ} \ \Hom_{\cA}(T^n,K^n). \]
This completes the proof of the lemma.
\end{pf}

\begin{lem}\label{lemB+3}
For $T, X \in \Ob(\rC(\cA))$ and $n \in \bZ$, there is a natural isomorphism
\[ \Hom^{n-1}(T,X)=\prod_{m \in \bZ} \ \Hom_{\cA}(T^{m-n+1},X^{m})
 \cong \Hom_{\rC(\cA)}(\cone(\id_T)[-n],X) \]
\end{lem}
\begin{pf}
We define a map $F:  \Hom^{n-1}(T,X) \to \Hom_{\rC(\cA)}(\cone(\id_T)[-n],X)$ by sending
\begin{align*}
 &(f_m:T^{m+1-n} \to X^m)_m \longmapsto \\
 &\qquad \quad ((d_X^{m-1}\circ f_{m-1}+f_{m}\circ d_T^{m+1-n},f_m): T^{m-n} \oplus T^{m+1-n} \to X^{m})_m,
\end{align*}
which is well-defined (i.e., the data on the right hand side gives a morphism in $\rC(\cA)$)
by a straight-forward computation (see \cite{GM} p.\ 154 for the differential of $\cone(\id_T)$).
On the other hand, we define a map $G :\Hom_{\rC(\cA)}(\cone(\id_T)[-n],X) \to \Hom^{n-1}(T,X)$ by sending
\[ ((a_m ,b_m): T^{m-n} \oplus T^{m+1-n} \to X^m)_m \longmapsto (b_m:T^{m+1-n} \to X^{m})_m. \]
It is obvious that $G \circ F = \id_{\Hom^{n-1}}$. Moreover, we obtain
\[ (F \circ G)((a_m ,b_m)_m)
 = (d_X^{m-1}\circ b_{m-1}+b_m \circ d_T^{m+1-n},b_m)_m = (a_m,b_m)_m, \]
by the condition that $(a_m ,b_m)_m$ is a morphism in $\rC(\cA)$. Thus we obtain the lemma.
\end{pf}

\begin{pf*}{Proof of Theorem \ref{thmB-1}}
{\bf (i)$\Ra$(ii).}\,
We first show that $X^n$ is injective for any $n \in \bZ$.
Let $\varphi : S \hra T$ be an arbitrary monomorphism in $\cA$.
Then the map $\wt{\varphi} : \cone(\id_S) \to \cone(\id_T)$ induced by $\varphi$ is a monomorphism by Lemma \ref{lemB+2}, and moreover it is a quasi-isomorphism, because $\cone(\id_S)$ and $\cone(\id_T)$ are acyclic. Hence by the lifting property of $X \to 0$, the map
\[  \Hom_{\rC(\cA)}(\cone(\id_T)[-n-1],X) \lra \Hom_{\rC(\cA)}(\cone(\id_S)[-n-1],X) \]
induced by $\wt{\varphi}[-n-1]$ is surjective, which implies that the map
\[ \varphi^* : \Hom_{\cA}(T,X^n) \lra \Hom_{\cA}(S,X^n) \]
is surjective by Lemma \ref{lemB+1}. Thus $X^n$ is injective.
We next show that $X$ is \tK-injective, i.e., the hom-complex $\Hom^\bullet(T,X)$ is acyclic for any acyclic $T \in \Ob(\rC(\cA))$.
Our task is to show that $\H^n(\Hom^\bullet(T,X))=0$ for any $n \in \bZ$.
Note that we have
\begin{align*}
 & \H^n(\Hom^\bullet(T,X)) =
\frac{\,\Ker(d_{\Hom}^n:\Hom^n(T,X)\to \Hom^{n+1}(T,X))\,}
{\Image(d_{\Hom}^{n-1}:\Hom^{n-1}(T,X)\to \Hom^n(T,X))}, \\
 &\Ker(d_{\Hom}^n:\Hom^n(T,X)\to \Hom^{n+1}(T,X)) = \Hom_{\rC(\cA)}(T[-n],X).
\end{align*}
Note also that we have
\[ \Hom^{n-1}(T,X)=\Hom_{\rC(\cA)}(\cone(\id_T)[-n],X) \]
by Lemma \ref{lemB+3}.
The natural map $i : T \to \cone(\id_T)$ is a monic quasi-isomorphism by Lemma \ref{lemB+2} and the acyclicity of  $T$. Hence the lifting property of $X \to 0$ implies that the map
\[ (i[-n])^* : \Hom_{\rC(\cA)}(\cone(\id_T)[-n],X) \lra \Hom_{\rC(\cA)}(T[-n],X) \]
is surjective. Noting that this map corresponds to the differential map of $\Hom^\bullet(T,X)$,
 one obtains that $\H^n(\Hom^\bullet(T,X)) = 0$. Thus $X$ is $\tK$-injective.

\par\medskip
{\bf (ii)$\Ra$(i).}\,
Let $f : K \hra L$ be a monic quasi-isomorphism in $\rC(\cA)$. Our task is to show that the map
\[ f^* : \Hom_{\rC(\cA)}(L,X) \lra \Hom_{\rC(\cA)}(K,X) \]
is surjective. Since the quotient complex $L/K$ is acyclic, the hom-complex $\Hom^\bullet(L/K,X)$ is acyclic by the assumption that $X$ is \tK-injective, and the map
\begin{equation}\label{eqB+1}
 f^* : \Hom^\bullet(L,X) \lra \Hom^\bullet(K,X)
\end{equation}
is a quasi-isomorphism in $\rC(\cA)$. Now noting that $Z^0(\Hom^\bullet(-,X))=\Hom_{\rC(\cA)}(-,X)$,
 consider a commutative diagram with exact rows in $\cA$
\[\xymatrix{
\Hom^{-1}(L,X) \ar[r]^-d \ar[d]^{f^*} & \Hom_{\rC(\cA)}(L,X) \ar[d]^{f^*} \ar[r]
 & \H^0(\Hom^\bullet(L,X)) \ar[d]^{f^*}_{\wr \!} \ar[r] & 0 \\
\Hom^{-1}(K,X) \ar[r]^-d & \Hom_{\rC(\cA)}(K,X) \ar[r] & \H^0(\Hom^\bullet(K,X)) \ar[r] & 0,\hspace{-3pt}
}\]
where the right vertical arrow is bijective by the fact that \eqref{eqB+1} is a quasi-isomorphism.
On the other hand, the left vertical arrow in the above diagram is surjective, because we have
\[ \Hom^{-1}(L,X)=\prod_{m \in \bZ} \ \Hom_{\cA}(L^m,X^{m-1}), \quad
\Hom^{-1}(K,X)=\prod_{m \in \bZ} \ \Hom_{\cA}(K^m,X^{m-1}), \]
and $f^m: K^m \to L^m$ is a monomorphism (resp.\ $X^{m-1}$ is injective) for any $m$ by Lemma \ref{lemB+2}
 (resp.\ by assumption).
Hence the middle vertical arrow is surjective as well, which shows the lifting property.
\end{pf*}

\subsection{Homotopy categories and derived category}\label{sectB-2}
We assume here that $\cA$ is a Grothendieck category.
Recall that the category $\rC(\cA)$ is endowed with the injective model structure, whose cofibrations are injective maps and whose weak equivalences are quasi-isomorphisms, cf.\ \cite{Bek} Proposition 3.13.
Let $\bI$ be the full subcategory of $\rC(\cA)$ consisting of injectively fibrant objects, and
let $\cHsn(\rC(\cA))$ be the homotopy category of $\rC(\cA)$ with respect to the injective model category structure.
On the other hand, let $\tK(\cA)$ be the homotopy category of $\rC(\cA)$ in the sense of homological algebra, i.e.,
the category whose objects are those of  $\rC(\cA)$ and whose morphisms are homotopy equivalence classes of morphisms in $\rC(\cA)$, cf.\ \cite{Ha0} Chapter I, \S2. 
The following fact for unbounded complexes is standard and was useful in Appendix A, but the authors do not know a written literature for it:

\begin{prop}\label{propB-1}
For $X \in \Ob(\rC(\cA))$ and $Y \in \Ob(\bI)$, there exists a canonical bijection
\[ \Hom_{\tK(\cA)}(X,Y) \cong \Hom_{\cHsn(\rC(\cA))}(X,Y). \]
\end{prop}
\noindent
To prove this proposition, the following standard construction will be useful:
\begin{defn}\label{defnB-1}
{\rm
For $X=(X^\bullet,d_X^\bullet) \in \Ob(\rC(\cA))$, let $\tX=(\tX^\bullet,d_{\tX}^\bullet) \in \Ob(\rC(\cA))$ be the mapping cylinder of the identity map $\id_X : X \to X$ (cf.\ {\rm \cite{GM}} p.\ {\rm 154}):
\begin{align*}
& \tX^n:=X^n \oplus X^{n+1}\oplus  X^n, \\
& d_{\tX}^n : \tX^n \lra \tX^{n+1},\qquad
(x,y,z) \mapsto  (d_X^n(x)-y,-d_X^{n+1}(y),y+d_X^n(z)).
\end{align*}
We define a morphism $i : X \oplus X \to \tX$ in $\rC(\cA)$ as $i(u,v) := (u,0,v)$, and define $j : \tX \to X$ in $\rC(\cA)$ as $j(x,y,z):=x+z$.
Then $i$ is a monomorphism by Lemma \ref{lemB+2}, and $j$ is a quasi-isomorphism by loc.\ cit., p.\ 155 Lemma 3. In other words, the diagram
\[ X \oplus X \os{i}\lra \tX \os{j}\lra X \]
is a cylinder object of $X$ with respect to the injective model structure, cf.\ \cite{GJ} p.\ 68.
}
\end{defn}

\begin{pf*}{Proof of Proposition \ref{propB-1}}
Let $X$ and $Y$ be as in Proposition \ref{propB-1}, and let $f, g : X \to Y$ be two morphisms in $\rC(\cA)$.
Our task to show that $f$ and $g$ are homotopic in the sense of homological algebra if and only if $f$ and $g$ are homotopic with respect to the injective model structure.
By the assumption that $Y \in \Ob(\bI)$ and \cite{GJ} Chapter \II, Corollary 1.9, it is enough to prove the following claim ($\spadesuit$) for $f, g : X \to Y$ in $\rC(\cA)$, where we do not assume that $Y \in \Ob(\bI)$:
\begin{enumerate}
\item[($\spadesuit$)]
{\it There exists a collection $(H^n : X^n \to Y^{n-1})_{n \in \bZ}$ of morphisms in $\cA$ satisfying the relation
\begin{equation}\label{eqB-1}
 d_Y^{n-1} \circ H^n + H^{n+1} \circ d_X^n = f^n - g^n \qquad ({}^\forall n \in \bZ),
\end{equation}
if and only if there exists a map $h : \tX \to Y$ in $\rC(\cA)$ fitting into the following commutative diagram}:
\begin{equation}\label{eqB-2}
\xymatrix{ X \oplus X \ar[d]_{i} \ar[rd]^{(f,g)} \\
 \tX \ar[r]_-h & Y.
}\end{equation}
\end{enumerate}
We prove ($\spadesuit$).
Assume that there exists a collection $(H^n : X^n \to Y^{n-1})_{n \in \bZ}$ of morphisms satisfying \eqref{eqB-1}.
Then we define the morphism
\[ h^n: \tX^n=X^n \oplus X^{n+1} \oplus X^n \lra Y^n \qquad (n \in \bZ) \]
by $h^n(x,y,z):=f^n(x)+H^{n+1}(y)+g^n(z)$. It is easy to check that $h=(h^n)_{n \in \bZ} : \tX \to Y$ is a morphism of complexes and fits into the commutative diagram \eqref{eqB-2}.
Conversely, suppose that we are given a morphism $h : \tX \to Y$ of complexes fitting into \eqref{eqB-2}.
We then define $H^n  : X^n \to Y^{n-1}$ by $H^n(x):=h^{n-1}(0,x,0)$.
One can easily check that the collection $(H^n : X^n \to Y^{n-1})_{n \in \bZ}$ satisfies the relation \eqref{eqB-1}.
\end{pf*}

By Proposition \ref{propB-1}, there is no fear of confusion in saying that two maps $f,g : X \to Y$ in $\rC(\cA)$ with $Y \in \Ob(\bI)$ are {\it homotopic}, either in the sense of homotopical or homological algebra.
We next show the following lemma, which extends the facts in \cite{Ha0} Chapter I, Lemmas 4.4 and 4.5 to unbounded complexes:

\begin{lem}\label{lemB-1}
\begin{enumerate}
\item[{\rm(1)}]
Let $f : X \to Y$ be a morphism in $\rC(\cA)$ with $Y \in \Ob(\bI)$, and assume that $X$ is acyclic.
Then $f$ is homotopic to zero.
\item[{\rm(2)}]
Let $g : Y \to Z$ be a quasi-isomorphism in $\rC(\cA)$ with $Y \in \Ob(\bI)$.
Then there exists a morphism $s : Z \to Y$ such that $s \circ g$ is homotopic to $\id_Y$.
\end{enumerate}
\end{lem}

\begin{pf}
(1)
Since $Y$ is \tK-injective by Theorem \ref{thmB-1} and $X$ is acyclic,
$\Hom^\bullet(X,Y)$ is acyclic. Hence we have
\[ \Hom_{\tK(\cA)}(X,Y) \cong \H^0(\Hom^\bullet(X,Y)) = 0, \]
which implies the assertion.
\par
(2) One can deduce the assertion from (1) by the same arguments as in the proof of  \cite{Ha0} Chapter I, Lemma 4.5.
\end{pf}
\smallskip
\noindent
The following corollary of Lemma \ref{lemB-1}\,(2) verifies that the derived category $\sfD(\cA)$ defined in the usual way has small hom-sets, under the assumption that $\cA$ is a Grothendieck category.
Compare with \cite{Bek} Remark 3.15,  and see also \cite{Ha0} Chapter I, Proposition 3.1 and the proof of loc.\ cit.\ Proposition 3.2 for the description of hom-sets of localized categories.

\begin{cor}\label{corB-1}
For $X \in \rC(\cA)$ and $Y \in \bI$, the inductive limit with respect to all quasi-isomorphisms $Y \os{\qis}\to Y'$ in $\rC(\cA)$
\[ \varinjlim_{Y \os{\qis}\to Y'} \ \Hom_{\tK(\cA)}(X,Y') \]
is small, and canonically isomorphic to $\Hom_{\tK(\cA)}(X,Y)$.
\end{cor}
\noindent
By Proposition \ref{propB-1} and Corollary \ref{corB-1}, we obtain the following corollary, which shows that the definition of the derived category of $\rC(\cA)$ given in \cite{Bek} Corollary 3.14 is consistent with the usual construction:
\begin{cor}\label{corB-2}
For $X,Y \in \rC(\cA)$, there exists a canonical bijection
\[ \Hom_{\cHsn(\rC(\cA))}(X,Y) \cong \Hom_{\sfD(\cA)}(X,Y). \]
Consequently, $\sfD(\cA)$ is canonically isomorphic to $\cHsn(\rC(\cA))$.
\end{cor}

\newpage

\noindent
Department of Mathematics, Hokkaido University
\par\noindent
Sapporo 060-0810,
JAPAN

\smallskip

\noindent
{\it E-mail} : \textbf{asakura@math.sci.hokudai.ac.jp}

\bigskip

\noindent
Department of Mathematics, Chuo University
\par\noindent
1-13-27 Kasuga, Bunkyo-ku,
Tokyo 112-8551,
JAPAN

\smallskip

\noindent
{\it Email} : \textbf{kanetomo\_\hspace{.8pt}sato@hotmail.com}
\bigskip

\noindent
Mathematical Science Team, RIKEN Center for Advanced Intelligence Project (AIP)
\par\noindent
1-4-1 Nihonbashi, Chuo-ku, Tokyo 103-0027, JAPAN
\par\smallskip\noindent
Department of Mathematics, Keio University
\par\noindent
3-14-1 Hiyoshi, Kohoku-ku,
Yokohama 223-8522,
JAPAN

\smallskip

\noindent
{\it Email} : \textbf{kei.hagihara@gmail.com}
\end{document}